\documentclass{article}
\usepackage{amsfonts}
\usepackage[section]{placeins}
\def\Erdos{Erd\H{o}s}

\def\Q{{\mathbb Q}}

\def\A{{\mathcal A}}
\def\dd{{\bf d}}

\def\sqN{{\sqrt{\frac N 2}}}  

\def\vector#1#2#3{\left(
\begin{array}{l} 
#1\\ 
#2\\ 
#3 
\end{array}\right)}

\def\matrix#1#2#3#4#5#6#7#8#9{{
 \left( \begin{array}{ccc}
#1 & #2 & #3 \\
#4 & #5 & #6 \\
#7 & #8 & #9 
\end{array} \right)}}

\def\dettwo#1#2#3#4{{
 \left\vert \begin{array}{ccc}
#1 & #2  \\
#3 & #4
\end{array} \right\vert}}

\usepackage{pstricks}
\begin{document}

\psset{unit=3cm}
\def\ThreeTiling{%
\pspicture(1, 1)
\qline(0,0)(1,0)            
\qline(0,0)(0.5,0.866025404)    
\qline(0.5,0.866025404)(1,0)    
\qline(0,0)(0.5,0.288675135)     
\qline(0.50,0.866025404)(0.50,0.288675135)  
\qline(0.50,0.288675135)(1,0)     
\endpspicture}

\def\ThreeTilingA{%
\pspicture(2, 1)
\qline(0,0)(2,0)            
\qline(0,0)(0.5,0.866025404)    
\qline(0.5,0.866025404)(2,0)    
\qline(0,0)(1,0.577350269)     
\qline(1,0.577350269)(1,0)  
\endpspicture}

\def\FourTilingA{%
\psset{unit=2cm}
\pspicture(2,1.2)
\qline(0,0)(1.73205081,0)              
\qline(1.73205081,0)(1.73205081,1)    
\qline(0.866025404,0.5)( 0.866025404,0)   
\qline( 0.866025404,0.5)(1.73205081,0)   
\qline( 0,0)(1.73205081,1) 
\qline( 0.866025404,0.5)(1.73205081,0.5)  
\endpspicture
\psset{unit=3cm}
}

\def\FourTilingB{%
\psset{unit=2cm}
\pspicture(2,1.2)
\qline(0,0)(1.73205081,0)              
\qline(1.73205081,0)(1.73205081,1)    
\qline(0.866025404,0.5)( 0.866025404,0)   
\qline( 0.866025404,0.5)(1.73205081,0)   
\qline( 0,0)(1.73205081,1) 
\qline(1.29903811,0.25)(1.73205081,1)
\endpspicture
\psset{unit=3cm}
}

\def\FourTilingC{%
\psset{unit=2cm}
\pspicture(2,1.2)
\qline(0,0)(1.73205081,0)              
\qline(1.73205081,0)(1.73205081,1)    
\qline(0.866025404,0.5)( 0.866025404,0)   
\qline( 0.866025404,0.5)(1.73205081,0)   
\qline( 0,0)(1.73205081,1) 
\qline(1.73205081,0)(1.29903811,0.75)
\endpspicture
\psset{unit=3cm}
}

\def\FourTiling{%
\pspicture(1,1)
\qline(0,0)(1,0)            
\qline(0,0)(0.5,0.866025404)    
\qline(0.5,0.866025404)(1,0)    
\qline(0.25,0.433012602)(0.75,0.433012602)  
\qline(0.25,0.433012602)(0.5,0)  
\qline(0.5,0)(0.75,0.433012602)   
\endpspicture}

\def\FiveTiling{%
\psset{unit=2cm}
\pspicture(2.5,1.2)
\qline(0,0)(2.5,0)              
\qline(0,0)(2,1)    
\qline(2,1)(2.5,0)   
\qline(2,1)(2,0)     
\qline(1,0.5)(2,0.5)  
\qline(2,0.5)(1,0)     
\qline(1,0)(1,0.5)    
\psset{unit=3cm}
\endpspicture}

\def\FiveTilingB{%
\psset{unit=2cm}
\pspicture(2.5,1.2)
\qline(0,0)(2.5,0)              
\qline(0,0)(2,1)    
\qline(2,1)(2.5,0)   
\qline(2,1)(2,0)     
\qline(1,0.5)(2,0.5)  
\qline(1,0.5)(2,0)     
\qline(1,0)(1,0.5)    
\psset{unit=3cm}
\endpspicture}

\def\SixTiling{%
\pspicture(1, 1)
\qline(0,0)(1,0)            
\qline(0,0)(0.5,0.866025404)    
\qline(0.5,0.866025404)(1,0)    
\qline(0,0)(0.75, 0.433012702)     
\qline(0.50,0.866025404)(0.50,0)  
\qline(0.25, 0.433012702)(1,0)     
\endpspicture}

\def\SixTilingA{%
\psset{unit=2cm}
\pspicture(1, 1)
\qline(0,0)( 3.46410162,0)            
\qline(0,0)(1.73205081,1)        
\qline(1.73205081,1)( 3.46410162,0)    
\qline(1.73205081,1)(1.73205081,0)  
\qline( 0.866025404,0.5)(1.15470054,0)   
\qline(1.15470054,0)(1.73205081,1)
\qline(1.73205081,1)(2.30940108,0)      
\qline(2.30940108,0)(2.59807621,0.5)    
\endpspicture}

\def\EightTiling{%
\pspicture(1,1)
\qline(0,0)(1,0)            
\qline(0,0)(0.5,0.866025404)    
\qline(0.5,0.866025404)(1,0)    
\qline(0.25,0.433012602)(0.75,0.433012602)  
\qline(0.25,0.433012602)(0.25,0)  
\qline(0.75,0.433012602)(0.75,0)  
\qline(0.5,0.866025404)(0.5,0)   
\qline(0.25,0.433012602)(0.5,0)  
\qline(0.5,0)(0.75,0.433012602)   
\endpspicture}

\def\NineTiling{%
\pspicture(1,1)
\qline(0,0)(1,0)            
\qline(0,0)(0.5,0.866025404)    
\qline(0.5,0.866025404)(1,0)    
\qline(0.167777,0.288675135)(0.833333,0.288675135)  
\qline(0.333333,0.577350269)(0.6777777, 0.577350269) 
\qline(0.1677777,0.288675135)(0.333333,0)
\qline(0.333333,0.577350269)(0.6777777,0)
\qline(0.333333,0)(0.677777,0.577350269)
\qline(0.677777,0)(0.833333,0.288675135)
\endpspicture}

\def\NineTilingA{%
\psset{unit=1cm}
\pspicture(6,3.3)
\qline(0,0)(6,0)            
\qline(0,0)(6,3)            
\qline(6,0)(6,3)            
\qline(2,0)(2,1)         
\qline(4,0)(4,2)          
\qline(5,0)(5,2)
\qline(4,2)(6,2)        
\qline(2,1)(4,1)
\qline(2,0)(4,1)        
\qline(4,0)(5,2)
\qline(5,0)(6,2)
\endpspicture
\psset{unit=3cm}
}

\def\TwelveTiling{%
\pspicture(1,1)
\qline(0,0)(1,0)            
\qline(0,0)(0.5,0.866025404)    
\qline(0.5,0.866025404)(1,0)    
\qline(0.25,0.433012602)(0.75,0.433012602)  
\qline(0.25,0.433012602)(0.5,0)  
\qline(0.5,0)(0.75,0.433012602)   
\qline(0,0)(0.25,0.144337567)
\qline(0.25,0.144337567)(0.25,0.433012602)
\qline(0.25,0.144337567)(0.5,0)
\qline(0.5,0)(0.75,0.144337567)
\qline(0.75,0.144337567)(0.75,0.433012602)
\qline(0.75,0.144337567)(1,0)
\qline(0.5,0.288675035)(0.5,0)
\qline(0.5,0.288675035)(0.25,0.433012602)
\qline(0.5,0.288675035)(0.75,0.433012602)
\qline(0.5,0.577350169)(0.5,0.866025404)
\qline(0.5,0.577350169)(0.25,0.433012602)
\qline(0.5,0.577350169)(0.75,0.433012602)
\endpspicture}

\def\ThirteenTiling{%
\hskip0.8in
\pspicture(2.9,1.2)
\psset{unit=0.5cm}
\newrgbcolor{lightblue}{0.8 0.8 1}
\newrgbcolor{pink}{1 0.8 0.8}
\newrgbcolor{lightgreen}{0.8 1 0.8}
\pspolygon[fillstyle=solid,linewidth=1pt,fillcolor=lightgreen](0.00,0.00)(9.00,6.00)(9.00,0.00)
\psline(9.00,6.00)(9.00,0.00)
\psline(6.00,4.00)(6.00,0.00)
\psline(3.00,2.00)(3.00,0.00)
\psline(0.00,0.00)(0.00,0.00)
\psline(9.00,0.00)(0.00,0.00)
\psline(9.00,2.00)(3.00,2.00)
\psline(9.00,4.00)(6.00,4.00)
\psline(9.00,6.00)(9.00,6.00)
\psline(0.00,0.00)(9.00,6.00)
\psline(3.00,0.00)(9.00,4.00)
\psline(6.00,0.00)(9.00,2.00)
\psline(9.00,0.00)(9.00,0.00)
\pspolygon[fillstyle=solid,linewidth=1pt,fillcolor=lightblue](9.00,6.00)(9.00,0.00)(13.00,0.00)
\psline(9.00,0.00)(13.00,0.00)
\psline(9.00,3.00)(11.00,3.00)
\psline(9.00,6.00)(9.00,6.00)
\psline(13.00,0.00)(9.00,6.00)
\psline(11.00,0.00)(9.00,3.00)
\psline(9.00,0.00)(9.00,0.00)
\psline(9.00,6.00)(9.00,0.00)
\psline(11.00,3.00)(11.00,0.00)
\psline(13.00,0.00)(13.00,0.00)
\endpspicture}

\def\SixteenTiling{%
\pspicture(1,1)
\qline(0,0)(1,0)            
\qline(0,0)(0.5,0.866025404)    
\qline(0.5,0.866025404)(1,0)    
\qline(0.125,0.216506351)(0.875,0.216506351)  
\qline(0.25,0.433012702)(0.75,  0.433012702) 
\qline(0.375,0.649519052)(0.625,0.649519052 ) 
\qline(0.125,0.216506351)(0.25,0)
\qline(0.25,0.433012702)(0.5,0)
\qline(0.375,0.649519052)(0.75,0)
\qline(0.25,0)(0.625,0.649519052 )
\qline(0.5,0)(0.75,  0.433012702)
\qline(0.75,0)(0.875,0.216506351)
\endpspicture}


\def\TwentySevenTiling{%
\psset{unit=0.5cm}
\pspicture(10.3923048,9.4)
\qline(0,0)(10.3923048,0)             
\qline(0,0)(5.19615242,9)             
\qline(5.19615242,9)(10.3923048,0)    
\qline(0,0)(1.73205081,1)             
\qline(1.73205081,1)(3.46410162,0)
\qline(3.46410162,0)(5.19615242,1)
\qline(5.19615242,1)(6.92820323,0)
\qline(6.92820323,0)(8.66025404,1)
\qline(8.66025404,1)(10.3923048,0)
\qline(1.73205081,1)(8.66025404,1)  
\qline(1.73205081,1)(1.73205081,3)  
\qline(5.19615242,1)(5.19615242,3)
\qline(8.66025404,1)(8.66025404,3)
\qline(1.73205081,3)(3.46410162,4)  
\qline(3.46410162,4)(5.19615242,3)
\qline(5.19615242,3)(6.92820323,4)
\qline(6.92820323,4)(8.66025404,3)  
\qline(3.46410162,4)(3.46410162,6)  
\qline(3.46410162,6)(5.19615242,7)  
\qline(5.19615242,7)(6.92820323,6)  
\qline(6.92820323,6)(6.92820323,4)  
\qline(5.19615242,7)(5.19615242,9)  
\qline(1.73205081,1)(5.19615242,7)  
\qline(5.19615242,7)(8.66025404,1)
\qline(3.46410162,4)(6.92820323,4)
\qline(3.46410162,4)(5.19615242,1)
\qline(5.19615242,1)(6.92820323,4)
\qline(3.46410162,2)(1.73205081,1)  
\qline(3.46410162,2)(3.46410162,4)
\qline(3.46410162,2)(5.19615242,1)
\qline(6.92820323,2)(5.19615242,1)  
\qline(6.92820323,2)(6.92820323,4)
\qline(6.92820323,2)(8.66025404,1)
\qline(5.19615242,5)(3.46410162,4)  
\qline(5.19615242,5)(6.92820323,4)
\qline(5.19615242,5)(5.19615242,7)
\endpspicture}

\def\TwelveTilingA{
\pspicture(1,1)
\psset{unit=0.005cm} 
\qline(409.5,10)(609.25,125.3257)
\qline(609.25,125.3257)(675.8334,10)
\qline(675.8334,10)(409.5,10)
\qline(809,240.6514)(809,10)
\qline(809,10)(675.8334,10)
\qline(675.8334,10)(809,240.6514)
\qline(809,240.6514)(609.25,125.3257)
\qline(609.25,125.3257)(675.8334,10)
\qline(675.8334,10)(809,240.6514)
\qline(809,240.6514)(609.25,125.3257)
\qline(609.25,125.3257)(542.6667,240.6514)
\qline(542.6667,240.6514)(809,240.6514)
\qline(409.5,10)(409.5,240.6514)
\qline(409.5,240.6514)(542.6667,240.6514)
\qline(542.6667,240.6514)(409.5,10)
\qline(409.5,10)(609.25,125.3257)
\qline(609.25,125.3257)(542.6667,240.6514)
\qline(542.6667,240.6514)(409.5,10)
\qline(9.999969,10)(209.75,125.3257)
\qline(209.75,125.3257)(276.3333,10)
\qline(276.3333,10)(9.999969,10)
\qline(409.5,240.6514)(409.5,10)
\qline(409.5,10)(276.3333,10)
\qline(276.3333,10)(409.5,240.6514)
\qline(409.5,240.6514)(209.75,125.3257)
\qline(209.75,125.3257)(276.3333,10)
\qline(276.3333,10)(409.5,240.6514)
\qline(409.5,240.6514)(609.25,355.9771)
\qline(609.25,355.9771)(675.8334,240.6514)
\qline(675.8334,240.6514)(409.5,240.6514)
\qline(809,471.3029)(809,240.6514)
\qline(809,240.6514)(675.8334,240.6514)
\qline(675.8334,240.6514)(809,471.3029)
\qline(809,471.3029)(609.25,355.9771)
\qline(609.25,355.9771)(675.8334,240.6514)
\qline(675.8334,240.6514)(809,471.3029)
\endpspicture}

\def\TwelveTilingB{
\pspicture(1,1)
\psset{unit=0.005cm} 
\qline(276.3333,10)(476.0833,125.3257)
\qline(476.0833,125.3257)(542.6667,10)
\qline(542.6667,10)(276.3333,10)
\qline(675.8333,240.6514)(675.8333,10)
\qline(675.8333,10)(542.6667,10)
\qline(542.6667,10)(675.8333,240.6514)
\qline(675.8333,240.6514)(476.0833,125.3257)
\qline(476.0833,125.3257)(542.6667,10)
\qline(542.6667,10)(675.8333,240.6514)
\qline(9.999969,10)(276.3333,10)
\qline(276.3333,10)(209.75,125.3257)
\qline(209.75,125.3257)(9.999969,10)
\qline(409.5,240.6514)(276.3333,10)
\qline(276.3333,10)(209.75,125.3257)
\qline(209.75,125.3257)(409.5,240.6514)
\qline(409.5,240.6514)(276.3333,10)
\qline(276.3333,10)(476.0833,125.3257)
\qline(476.0833,125.3257)(409.5,240.6514)
\qline(409.5,240.6514)(675.8334,240.6514)
\qline(675.8334,240.6514)(609.2501,355.9771)
\qline(609.2501,355.9771)(409.5,240.6514)
\qline(809,471.3029)(675.8334,240.6514)
\qline(675.8334,240.6514)(609.2501,355.9771)
\qline(609.2501,355.9771)(809,471.3029)
\qline(675.8333,10)(809,10)
\qline(809,10)(675.8333,240.6514)
\qline(675.8333,240.6514)(675.8333,10)
\qline(809,240.6514)(809,10)
\qline(809,10)(675.8333,240.6514)
\qline(675.8333,240.6514)(809,240.6514)
\qline(675.8333,240.6514)(809,240.6514)
\qline(809,240.6514)(809,471.3029)
\qline(809,471.3029)(675.8333,240.6514)
\endpspicture}
\def\TwentySevenTilingA{
\pspicture(1.0,1.0)
\psset{unit=0.005cm} 
\qline(542.6667,10)(675.8334,86.88382)
\qline(675.8334,86.88382)(720.2222,10)
\qline(720.2222,10)(542.6667,10)
\qline(809,163.7676)(809,10)
\qline(809,10)(720.2222,10)
\qline(720.2222,10)(809,163.7676)
\qline(809,163.7676)(675.8334,86.88382)
\qline(675.8334,86.88382)(720.2222,10)
\qline(720.2222,10)(809,163.7676)
\qline(809,163.7676)(675.8334,86.88382)
\qline(675.8334,86.88382)(631.4445,163.7676)
\qline(631.4445,163.7676)(809,163.7676)
\qline(542.6667,10)(542.6667,163.7676)
\qline(542.6667,163.7676)(631.4445,163.7676)
\qline(631.4445,163.7676)(542.6667,10)
\qline(542.6667,10)(675.8334,86.88382)
\qline(675.8334,86.88382)(631.4445,163.7676)
\qline(631.4445,163.7676)(542.6667,10)
\qline(276.3333,10)(409.5,86.88382)
\qline(409.5,86.88382)(453.8889,10)
\qline(453.8889,10)(276.3333,10)
\qline(542.6667,163.7676)(542.6667,10)
\qline(542.6667,10)(453.8889,10)
\qline(453.8889,10)(542.6667,163.7676)
\qline(542.6667,163.7676)(409.5,86.88382)
\qline(409.5,86.88382)(453.8889,10)
\qline(453.8889,10)(542.6667,163.7676)
\qline(542.6667,163.7676)(409.5,86.88382)
\qline(409.5,86.88382)(365.1111,163.7676)
\qline(365.1111,163.7676)(542.6667,163.7676)
\qline(276.3333,10)(276.3333,163.7676)
\qline(276.3333,163.7676)(365.1111,163.7676)
\qline(365.1111,163.7676)(276.3333,10)
\qline(276.3333,10)(409.5,86.88382)
\qline(409.5,86.88382)(365.1111,163.7676)
\qline(365.1111,163.7676)(276.3333,10)
\qline(9.999992,10)(143.1667,86.88382)
\qline(143.1667,86.88382)(187.5555,10)
\qline(187.5555,10)(9.999992,10)
\qline(276.3333,163.7676)(276.3333,10)
\qline(276.3333,10)(187.5555,10)
\qline(187.5555,10)(276.3333,163.7676)
\qline(276.3333,163.7676)(143.1667,86.88382)
\qline(143.1667,86.88382)(187.5555,10)
\qline(187.5555,10)(276.3333,163.7676)
\qline(542.6667,163.7676)(675.8334,240.6514)
\qline(675.8334,240.6514)(720.2222,163.7676)
\qline(720.2222,163.7676)(542.6667,163.7676)
\qline(809,317.5352)(809,163.7676)
\qline(809,163.7676)(720.2222,163.7676)
\qline(720.2222,163.7676)(809,317.5352)
\qline(809,317.5352)(675.8334,240.6514)
\qline(675.8334,240.6514)(720.2222,163.7676)
\qline(720.2222,163.7676)(809,317.5352)
\qline(809,317.5352)(675.8334,240.6514)
\qline(675.8334,240.6514)(631.4445,317.5352)
\qline(631.4445,317.5352)(809,317.5352)
\qline(542.6667,163.7676)(542.6667,317.5353)
\qline(542.6667,317.5353)(631.4445,317.5352)
\qline(631.4445,317.5352)(542.6667,163.7676)
\qline(542.6667,163.7676)(675.8334,240.6514)
\qline(675.8334,240.6514)(631.4445,317.5352)
\qline(631.4445,317.5352)(542.6667,163.7676)
\qline(276.3333,163.7676)(409.5,240.6514)
\qline(409.5,240.6514)(453.8889,163.7676)
\qline(453.8889,163.7676)(276.3333,163.7676)
\qline(542.6667,317.5352)(542.6667,163.7676)
\qline(542.6667,163.7676)(453.8889,163.7676)
\qline(453.8889,163.7676)(542.6667,317.5352)
\qline(542.6667,317.5352)(409.5,240.6514)
\qline(409.5,240.6514)(453.8889,163.7676)
\qline(453.8889,163.7676)(542.6667,317.5352)
\qline(542.6667,317.5352)(675.8334,394.4191)
\qline(675.8334,394.4191)(720.2222,317.5352)
\qline(720.2222,317.5352)(542.6667,317.5352)
\qline(809,471.3029)(809,317.5352)
\qline(809,317.5352)(720.2222,317.5352)
\qline(720.2222,317.5352)(809,471.3029)
\qline(809,471.3029)(675.8334,394.4191)
\qline(675.8334,394.4191)(720.2222,317.5352)
\qline(720.2222,317.5352)(809,471.3029)
\endpspicture}
\def\TwentySevenTilingB{
\pspicture(1.0,1.0)
\psset{unit=0.005cm} 
\qline(453.8889,10)(587.0555,86.88382)
\qline(587.0555,86.88382)(631.4445,10)
\qline(631.4445,10)(453.8889,10)
\qline(720.2222,163.7676)(720.2222,10)
\qline(720.2222,10)(631.4445,10)
\qline(631.4445,10)(720.2222,163.7676)
\qline(720.2222,163.7676)(587.0555,86.88382)
\qline(587.0555,86.88382)(631.4445,10)
\qline(631.4445,10)(720.2222,163.7676)
\qline(720.2222,163.7676)(587.0555,86.88382)
\qline(587.0555,86.88382)(542.6667,163.7676)
\qline(542.6667,163.7676)(720.2222,163.7676)
\qline(453.8889,10)(453.8889,163.7676)
\qline(453.8889,163.7676)(542.6667,163.7676)
\qline(542.6667,163.7676)(453.8889,10)
\qline(453.8889,10)(587.0555,86.88382)
\qline(587.0555,86.88382)(542.6667,163.7676)
\qline(542.6667,163.7676)(453.8889,10)
\qline(187.5555,10)(320.7222,86.88382)
\qline(320.7222,86.88382)(365.1111,10)
\qline(365.1111,10)(187.5555,10)
\qline(453.8889,163.7676)(453.8889,10)
\qline(453.8889,10)(365.1111,10)
\qline(365.1111,10)(453.8889,163.7676)
\qline(453.8889,163.7676)(320.7222,86.88382)
\qline(320.7222,86.88382)(365.1111,10)
\qline(365.1111,10)(453.8889,163.7676)
\qline(453.8889,163.7676)(587.0555,240.6514)
\qline(587.0555,240.6514)(631.4445,163.7676)
\qline(631.4445,163.7676)(453.8889,163.7676)
\qline(720.2222,317.5352)(720.2222,163.7676)
\qline(720.2222,163.7676)(631.4445,163.7676)
\qline(631.4445,163.7676)(720.2222,317.5352)
\qline(720.2222,317.5352)(587.0555,240.6514)
\qline(587.0555,240.6514)(631.4445,163.7676)
\qline(631.4445,163.7676)(720.2222,317.5352)
\qline(9.999992,10)(187.5555,10)
\qline(187.5555,10)(143.1667,86.88382)
\qline(143.1667,86.88382)(9.999992,10)
\qline(276.3333,163.7676)(187.5555,10)
\qline(187.5555,10)(143.1667,86.88382)
\qline(143.1667,86.88382)(276.3333,163.7676)
\qline(276.3333,163.7676)(187.5555,10)
\qline(187.5555,10)(320.7222,86.88382)
\qline(320.7222,86.88382)(276.3333,163.7676)
\qline(276.3333,163.7676)(453.8889,163.7676)
\qline(453.8889,163.7676)(409.5,240.6515)
\qline(409.5,240.6515)(276.3333,163.7676)
\qline(542.6666,317.5353)(453.8889,163.7676)
\qline(453.8889,163.7676)(409.5,240.6515)
\qline(409.5,240.6515)(542.6666,317.5353)
\qline(542.6666,317.5353)(453.8889,163.7676)
\qline(453.8889,163.7676)(587.0555,240.6515)
\qline(587.0555,240.6515)(542.6666,317.5353)
\qline(542.6667,317.5352)(720.2222,317.5352)
\qline(720.2222,317.5352)(675.8333,394.4191)
\qline(675.8333,394.4191)(542.6667,317.5352)
\qline(809,471.3029)(720.2222,317.5352)
\qline(720.2222,317.5352)(675.8333,394.4191)
\qline(675.8333,394.4191)(809,471.3029)
\qline(720.2222,10)(809,10)
\qline(809,10)(809,163.7676)
\qline(809,163.7676)(720.2222,10)
\qline(720.2222,163.7676)(720.2222,10)
\qline(720.2222,10)(809,163.7676)
\qline(809,163.7676)(720.2222,163.7676)
\qline(720.2222,163.7676)(809,163.7676)
\qline(809,163.7676)(720.2222,317.5353)
\qline(720.2222,317.5353)(720.2222,163.7676)
\qline(809,317.5353)(809,163.7676)
\qline(809,163.7676)(720.2222,317.5353)
\qline(720.2222,317.5353)(809,317.5353)
\qline(720.2222,317.5352)(809,317.5352)
\qline(809,317.5352)(809,471.3029)
\qline(809,471.3029)(720.2222,317.5352)
\endpspicture}

\def\FortyEightTiling{
\pspicture(2.0,2.0)
\psset{unit=0.01cm} 
\qline(10,10)(148.75,10)
\qline(148.75,10)(79.375,50.05365)
\qline(79.375,50.05365)(10,10)
\qline(148.75,10)(287.5,10)
\qline(287.5,10)(218.125,50.05365)
\qline(218.125,50.05365)(148.75,10)
\qline(287.5,10)(426.25,10)
\qline(426.25,10)(356.875,50.05365)
\qline(356.875,50.05365)(287.5,10)
\qline(426.25,10)(565,10)
\qline(565,10)(495.625,50.05365)
\qline(495.625,50.05365)(426.25,10)
\qline(10,10)(79.375,50.05365)
\qline(79.375,50.05365)(79.375,130.161)
\qline(79.375,130.161)(10,10)
\qline(79.375,130.161)(148.75,170.2147)
\qline(148.75,170.2147)(148.75,250.3221)
\qline(148.75,250.3221)(79.375,130.161)
\qline(148.75,250.3221)(218.125,290.3757)
\qline(218.125,290.3757)(218.125,370.4831)
\qline(218.125,370.4831)(148.75,250.3221)
\qline(218.125,370.4831)(287.5,410.5368)
\qline(287.5,410.5368)(287.5,490.6441)
\qline(287.5,490.6441)(218.125,370.4831)
\qline(565,10)(495.625,50.05365)
\qline(495.625,50.05365)(495.625,130.161)
\qline(495.625,130.161)(565,10)
\qline(495.625,130.161)(426.25,170.2147)
\qline(426.25,170.2147)(426.25,250.3221)
\qline(426.25,250.3221)(495.625,130.161)
\qline(426.25,250.3221)(356.875,290.3757)
\qline(356.875,290.3757)(356.875,370.4831)
\qline(356.875,370.4831)(426.25,250.3221)
\qline(356.875,370.4831)(287.5,410.5368)
\qline(287.5,410.5368)(287.5,490.6441)
\qline(287.5,490.6441)(356.875,370.4831)
\qline(287.5,410.5368)(218.125,290.3757)
\qline(218.125,290.3757)(218.125,370.4831)
\qline(218.125,370.4831)(287.5,410.5368)
\qline(287.5,410.5368)(218.125,290.3757)
\qline(218.125,290.3757)(287.5,330.4294)
\qline(287.5,330.4294)(287.5,410.5368)
\qline(218.125,290.3757)(356.875,290.3757)
\qline(356.875,290.3757)(287.5,250.3221)
\qline(287.5,250.3221)(218.125,290.3757)
\qline(218.125,290.3757)(356.875,290.3757)
\qline(356.875,290.3757)(287.5,330.4294)
\qline(287.5,330.4294)(218.125,290.3757)
\qline(356.875,290.3757)(287.5,410.5368)
\qline(287.5,410.5368)(356.875,370.4831)
\qline(356.875,370.4831)(356.875,290.3757)
\qline(356.875,290.3757)(287.5,410.5368)
\qline(287.5,410.5368)(287.5,330.4294)
\qline(287.5,330.4294)(356.875,290.3757)
\qline(218.125,290.3757)(148.75,170.2147)
\qline(148.75,170.2147)(148.75,250.3221)
\qline(148.75,250.3221)(218.125,290.3757)
\qline(218.125,290.3757)(148.75,170.2147)
\qline(148.75,170.2147)(218.125,210.2684)
\qline(218.125,210.2684)(218.125,290.3757)
\qline(148.75,170.2147)(287.5,170.2147)
\qline(287.5,170.2147)(218.125,130.161)
\qline(218.125,130.161)(148.75,170.2147)
\qline(148.75,170.2147)(287.5,170.2147)
\qline(287.5,170.2147)(218.125,210.2684)
\qline(218.125,210.2684)(148.75,170.2147)
\qline(287.5,170.2147)(218.125,290.3757)
\qline(218.125,290.3757)(287.5,250.3221)
\qline(287.5,250.3221)(287.5,170.2147)
\qline(287.5,170.2147)(218.125,290.3757)
\qline(218.125,290.3757)(218.125,210.2684)
\qline(218.125,210.2684)(287.5,170.2147)
\qline(356.875,290.3757)(287.5,170.2147)
\qline(287.5,170.2147)(287.5,250.3221)
\qline(287.5,250.3221)(356.875,290.3757)
\qline(356.875,290.3757)(287.5,170.2147)
\qline(287.5,170.2147)(356.875,210.2684)
\qline(356.875,210.2684)(356.875,290.3757)
\qline(287.5,170.2147)(426.25,170.2147)
\qline(426.25,170.2147)(356.875,130.161)
\qline(356.875,130.161)(287.5,170.2147)
\qline(287.5,170.2147)(426.25,170.2147)
\qline(426.25,170.2147)(356.875,210.2684)
\qline(356.875,210.2684)(287.5,170.2147)
\qline(426.25,170.2147)(356.875,290.3757)
\qline(356.875,290.3757)(426.25,250.3221)
\qline(426.25,250.3221)(426.25,170.2147)
\qline(426.25,170.2147)(356.875,290.3757)
\qline(356.875,290.3757)(356.875,210.2684)
\qline(356.875,210.2684)(426.25,170.2147)
\qline(148.75,170.2147)(79.375,50.05371)
\qline(79.375,50.05371)(79.375,130.161)
\qline(79.375,130.161)(148.75,170.2147)
\qline(148.75,170.2147)(79.375,50.05371)
\qline(79.375,50.05371)(148.75,90.10736)
\qline(148.75,90.10736)(148.75,170.2147)
\qline(79.375,50.05371)(218.125,50.05371)
\qline(218.125,50.05371)(148.75,10)
\qline(148.75,10)(79.375,50.05371)
\qline(79.375,50.05371)(218.125,50.05371)
\qline(218.125,50.05371)(148.75,90.10736)
\qline(148.75,90.10736)(79.375,50.05371)
\qline(218.125,50.05371)(148.75,170.2147)
\qline(148.75,170.2147)(218.125,130.161)
\qline(218.125,130.161)(218.125,50.05371)
\qline(218.125,50.05371)(148.75,170.2147)
\qline(148.75,170.2147)(148.75,90.10736)
\qline(148.75,90.10736)(218.125,50.05371)
\qline(287.5,170.2147)(218.125,50.05371)
\qline(218.125,50.05371)(218.125,130.161)
\qline(218.125,130.161)(287.5,170.2147)
\qline(287.5,170.2147)(218.125,50.05371)
\qline(218.125,50.05371)(287.5,90.10736)
\qline(287.5,90.10736)(287.5,170.2147)
\qline(218.125,50.05371)(356.875,50.05371)
\qline(356.875,50.05371)(287.5,10)
\qline(287.5,10)(218.125,50.05371)
\qline(218.125,50.05371)(356.875,50.05371)
\qline(356.875,50.05371)(287.5,90.10736)
\qline(287.5,90.10736)(218.125,50.05371)
\qline(356.875,50.05371)(287.5,170.2147)
\qline(287.5,170.2147)(356.875,130.161)
\qline(356.875,130.161)(356.875,50.05371)
\qline(356.875,50.05371)(287.5,170.2147)
\qline(287.5,170.2147)(287.5,90.10736)
\qline(287.5,90.10736)(356.875,50.05371)
\qline(426.25,170.2147)(356.875,50.05371)
\qline(356.875,50.05371)(356.875,130.161)
\qline(356.875,130.161)(426.25,170.2147)
\qline(426.25,170.2147)(356.875,50.05371)
\qline(356.875,50.05371)(426.25,90.10736)
\qline(426.25,90.10736)(426.25,170.2147)
\qline(356.875,50.05371)(495.625,50.05371)
\qline(495.625,50.05371)(426.25,10)
\qline(426.25,10)(356.875,50.05371)
\qline(356.875,50.05371)(495.625,50.05371)
\qline(495.625,50.05371)(426.25,90.10736)
\qline(426.25,90.10736)(356.875,50.05371)
\qline(495.625,50.05371)(426.25,170.2147)
\qline(426.25,170.2147)(495.625,130.161)
\qline(495.625,130.161)(495.625,50.05371)
\qline(495.625,50.05371)(426.25,170.2147)
\qline(426.25,170.2147)(426.25,90.10736)
\qline(426.25,90.10736)(495.625,50.05371)
\endpspicture}

\def\OneHundredTwentyFiveTiling{
\pspicture(2.0,2.0)
\psset{unit=0.01cm} 
\qline(10,10)(121,10)
\qline(121,10)(65.5,42.04297)
\qline(65.5,42.04297)(10,10)
\qline(121,10)(232,10)
\qline(232,10)(176.5,42.04297)
\qline(176.5,42.04297)(121,10)
\qline(232,10)(343,10)
\qline(343,10)(287.5,42.04297)
\qline(287.5,42.04297)(232,10)
\qline(343,10)(454,10)
\qline(454,10)(398.5,42.04297)
\qline(398.5,42.04297)(343,10)
\qline(454,10)(565,10)
\qline(565,10)(509.5,42.04297)
\qline(509.5,42.04297)(454,10)
\qline(10,10)(65.5,42.04297)
\qline(65.5,42.04297)(65.5,106.1288)
\qline(65.5,106.1288)(10,10)
\qline(65.5,106.1288)(121,138.1718)
\qline(121,138.1718)(121,202.2576)
\qline(121,202.2576)(65.5,106.1288)
\qline(121,202.2577)(176.5,234.3006)
\qline(176.5,234.3006)(176.5,298.3865)
\qline(176.5,298.3865)(121,202.2577)
\qline(176.5,298.3865)(232,330.4294)
\qline(232,330.4294)(232,394.5153)
\qline(232,394.5153)(176.5,298.3865)
\qline(232,394.5153)(287.5,426.5582)
\qline(287.5,426.5582)(287.5,490.6441)
\qline(287.5,490.6441)(232,394.5153)
\qline(565,10)(509.5,42.04297)
\qline(509.5,42.04297)(509.5,106.1288)
\qline(509.5,106.1288)(565,10)
\qline(509.5,106.1288)(454,138.1718)
\qline(454,138.1718)(454,202.2576)
\qline(454,202.2576)(509.5,106.1288)
\qline(454,202.2577)(398.5,234.3006)
\qline(398.5,234.3006)(398.5,298.3865)
\qline(398.5,298.3865)(454,202.2577)
\qline(398.5,298.3865)(343,330.4294)
\qline(343,330.4294)(343,394.5153)
\qline(343,394.5153)(398.5,298.3865)
\qline(343,394.5153)(287.5,426.5582)
\qline(287.5,426.5582)(287.5,490.6441)
\qline(287.5,490.6441)(343,394.5153)
\qline(287.5,426.5583)(232,330.4294)
\qline(232,330.4294)(232,394.5153)
\qline(232,394.5153)(287.5,426.5583)
\qline(287.5,426.5583)(232,330.4294)
\qline(232,330.4294)(287.5,362.4724)
\qline(287.5,362.4724)(287.5,426.5583)
\qline(232,330.4294)(343,330.4294)
\qline(343,330.4294)(287.5,298.3865)
\qline(287.5,298.3865)(232,330.4294)
\qline(232,330.4294)(343,330.4294)
\qline(343,330.4294)(287.5,362.4724)
\qline(287.5,362.4724)(232,330.4294)
\qline(343,330.4294)(287.5,426.5583)
\qline(287.5,426.5583)(343,394.5153)
\qline(343,394.5153)(343,330.4294)
\qline(343,330.4294)(287.5,426.5583)
\qline(287.5,426.5583)(287.5,362.4724)
\qline(287.5,362.4724)(343,330.4294)
\qline(232,330.4294)(176.5,234.3006)
\qline(176.5,234.3006)(176.5,298.3865)
\qline(176.5,298.3865)(232,330.4294)
\qline(232,330.4294)(176.5,234.3006)
\qline(176.5,234.3006)(232,266.3435)
\qline(232,266.3435)(232,330.4294)
\qline(176.5,234.3006)(287.5,234.3006)
\qline(287.5,234.3006)(232,202.2577)
\qline(232,202.2577)(176.5,234.3006)
\qline(176.5,234.3006)(287.5,234.3006)
\qline(287.5,234.3006)(232,266.3435)
\qline(232,266.3435)(176.5,234.3006)
\qline(287.5,234.3006)(232,330.4294)
\qline(232,330.4294)(287.5,298.3865)
\qline(287.5,298.3865)(287.5,234.3006)
\qline(287.5,234.3006)(232,330.4294)
\qline(232,330.4294)(232,266.3435)
\qline(232,266.3435)(287.5,234.3006)
\qline(343,330.4294)(287.5,234.3006)
\qline(287.5,234.3006)(287.5,298.3865)
\qline(287.5,298.3865)(343,330.4294)
\qline(343,330.4294)(287.5,234.3006)
\qline(287.5,234.3006)(343,266.3435)
\qline(343,266.3435)(343,330.4294)
\qline(287.5,234.3006)(398.5,234.3006)
\qline(398.5,234.3006)(343,202.2577)
\qline(343,202.2577)(287.5,234.3006)
\qline(287.5,234.3006)(398.5,234.3006)
\qline(398.5,234.3006)(343,266.3435)
\qline(343,266.3435)(287.5,234.3006)
\qline(398.5,234.3006)(343,330.4294)
\qline(343,330.4294)(398.5,298.3865)
\qline(398.5,298.3865)(398.5,234.3006)
\qline(398.5,234.3006)(343,330.4294)
\qline(343,330.4294)(343,266.3435)
\qline(343,266.3435)(398.5,234.3006)
\qline(176.5,234.3006)(121,138.1718)
\qline(121,138.1718)(121,202.2577)
\qline(121,202.2577)(176.5,234.3006)
\qline(176.5,234.3006)(121,138.1718)
\qline(121,138.1718)(176.5,170.2147)
\qline(176.5,170.2147)(176.5,234.3006)
\qline(121,138.1718)(232,138.1718)
\qline(232,138.1718)(176.5,106.1288)
\qline(176.5,106.1288)(121,138.1718)
\qline(121,138.1718)(232,138.1718)
\qline(232,138.1718)(176.5,170.2147)
\qline(176.5,170.2147)(121,138.1718)
\qline(232,138.1718)(176.5,234.3006)
\qline(176.5,234.3006)(232,202.2577)
\qline(232,202.2577)(232,138.1718)
\qline(232,138.1718)(176.5,234.3006)
\qline(176.5,234.3006)(176.5,170.2147)
\qline(176.5,170.2147)(232,138.1718)
\qline(287.5,234.3006)(232,138.1718)
\qline(232,138.1718)(232,202.2577)
\qline(232,202.2577)(287.5,234.3006)
\qline(287.5,234.3006)(232,138.1718)
\qline(232,138.1718)(287.5,170.2147)
\qline(287.5,170.2147)(287.5,234.3006)
\qline(232,138.1718)(343,138.1718)
\qline(343,138.1718)(287.5,106.1288)
\qline(287.5,106.1288)(232,138.1718)
\qline(232,138.1718)(343,138.1718)
\qline(343,138.1718)(287.5,170.2147)
\qline(287.5,170.2147)(232,138.1718)
\qline(343,138.1718)(287.5,234.3006)
\qline(287.5,234.3006)(343,202.2577)
\qline(343,202.2577)(343,138.1718)
\qline(343,138.1718)(287.5,234.3006)
\qline(287.5,234.3006)(287.5,170.2147)
\qline(287.5,170.2147)(343,138.1718)
\qline(398.5,234.3006)(343,138.1718)
\qline(343,138.1718)(343,202.2577)
\qline(343,202.2577)(398.5,234.3006)
\qline(398.5,234.3006)(343,138.1718)
\qline(343,138.1718)(398.5,170.2147)
\qline(398.5,170.2147)(398.5,234.3006)
\qline(343,138.1718)(454,138.1718)
\qline(454,138.1718)(398.5,106.1288)
\qline(398.5,106.1288)(343,138.1718)
\qline(343,138.1718)(454,138.1718)
\qline(454,138.1718)(398.5,170.2147)
\qline(398.5,170.2147)(343,138.1718)
\qline(454,138.1718)(398.5,234.3006)
\qline(398.5,234.3006)(454,202.2577)
\qline(454,202.2577)(454,138.1718)
\qline(454,138.1718)(398.5,234.3006)
\qline(398.5,234.3006)(398.5,170.2147)
\qline(398.5,170.2147)(454,138.1718)
\qline(121,138.1718)(65.5,42.04291)
\qline(65.5,42.04291)(65.5,106.1288)
\qline(65.5,106.1288)(121,138.1718)
\qline(121,138.1718)(65.5,42.04291)
\qline(65.5,42.04291)(121,74.08588)
\qline(121,74.08588)(121,138.1718)
\qline(65.5,42.04291)(176.5,42.04291)
\qline(176.5,42.04291)(121,10)
\qline(121,10)(65.5,42.04291)
\qline(65.5,42.04291)(176.5,42.04291)
\qline(176.5,42.04291)(121,74.08588)
\qline(121,74.08588)(65.5,42.04291)
\qline(176.5,42.04291)(121,138.1718)
\qline(121,138.1718)(176.5,106.1288)
\qline(176.5,106.1288)(176.5,42.04291)
\qline(176.5,42.04291)(121,138.1718)
\qline(121,138.1718)(121,74.08588)
\qline(121,74.08588)(176.5,42.04291)
\qline(232,138.1718)(176.5,42.04291)
\qline(176.5,42.04291)(176.5,106.1288)
\qline(176.5,106.1288)(232,138.1718)
\qline(232,138.1718)(176.5,42.04291)
\qline(176.5,42.04291)(232,74.08588)
\qline(232,74.08588)(232,138.1718)
\qline(176.5,42.04291)(287.5,42.04291)
\qline(287.5,42.04291)(232,10)
\qline(232,10)(176.5,42.04291)
\qline(176.5,42.04291)(287.5,42.04291)
\qline(287.5,42.04291)(232,74.08588)
\qline(232,74.08588)(176.5,42.04291)
\qline(287.5,42.04291)(232,138.1718)
\qline(232,138.1718)(287.5,106.1288)
\qline(287.5,106.1288)(287.5,42.04291)
\qline(287.5,42.04291)(232,138.1718)
\qline(232,138.1718)(232,74.08588)
\qline(232,74.08588)(287.5,42.04291)
\qline(343,138.1718)(287.5,42.04291)
\qline(287.5,42.04291)(287.5,106.1288)
\qline(287.5,106.1288)(343,138.1718)
\qline(343,138.1718)(287.5,42.04291)
\qline(287.5,42.04291)(343,74.08588)
\qline(343,74.08588)(343,138.1718)
\qline(287.5,42.04291)(398.5,42.04291)
\qline(398.5,42.04291)(343,10)
\qline(343,10)(287.5,42.04291)
\qline(287.5,42.04291)(398.5,42.04291)
\qline(398.5,42.04291)(343,74.08588)
\qline(343,74.08588)(287.5,42.04291)
\qline(398.5,42.04291)(343,138.1718)
\qline(343,138.1718)(398.5,106.1288)
\qline(398.5,106.1288)(398.5,42.04291)
\qline(398.5,42.04291)(343,138.1718)
\qline(343,138.1718)(343,74.08588)
\qline(343,74.08588)(398.5,42.04291)
\qline(454,138.1718)(398.5,42.04291)
\qline(398.5,42.04291)(398.5,106.1288)
\qline(398.5,106.1288)(454,138.1718)
\qline(454,138.1718)(398.5,42.04291)
\qline(398.5,42.04291)(454,74.08588)
\qline(454,74.08588)(454,138.1718)
\qline(398.5,42.04291)(509.5,42.04291)
\qline(509.5,42.04291)(454,10)
\qline(454,10)(398.5,42.04291)
\qline(398.5,42.04291)(509.5,42.04291)
\qline(509.5,42.04291)(454,74.08588)
\qline(454,74.08588)(398.5,42.04291)
\qline(509.5,42.04291)(454,138.1718)
\qline(454,138.1718)(509.5,106.1288)
\qline(509.5,106.1288)(509.5,42.04291)
\qline(509.5,42.04291)(454,138.1718)
\qline(454,138.1718)(454,74.08588)
\qline(454,74.08588)(509.5,42.04291)
\endpspicture}
 
\def\FigureBiquadratic{
\hskip 1.3in
\pspicture(2.0,1.3)
\psset{unit=0.1cm}
\newrgbcolor{lightblue}{0.8 0.8 1}
\newrgbcolor{pink}{1 0.8 0.8}
\newrgbcolor{lightgreen}{0.8 1 0.8}
\pspolygon[fillstyle=solid,linewidth=1pt,fillcolor=lightgreen](0.00,0.00)(49.00,35.00)(49.00,0.00)
\psline(49.00,35.00)(49.00,0.00)
\psline(42.00,30.00)(42.00,0.00)
\psline(35.00,25.00)(35.00,0.00)
\psline(28.00,20.00)(28.00,0.00)
\psline(21.00,15.00)(21.00,0.00)
\psline(14.00,10.00)(14.00,0.00)
\psline(7.00,5.00)(7.00,0.00)
\psline(0.00,0.00)(0.00,0.00)
\psline(49.00,0.00)(0.00,0.00)
\psline(49.00,5.00)(7.00,5.00)
\psline(49.00,10.00)(14.00,10.00)
\psline(49.00,15.00)(21.00,15.00)
\psline(49.00,20.00)(28.00,20.00)
\psline(49.00,25.00)(35.00,25.00)
\psline(49.00,30.00)(42.00,30.00)
\psline(49.00,35.00)(49.00,35.00)
\psline(0.00,0.00)(49.00,35.00)
\psline(7.00,0.00)(49.00,30.00)
\psline(14.00,0.00)(49.00,25.00)
\psline(21.00,0.00)(49.00,20.00)
\psline(28.00,0.00)(49.00,15.00)
\psline(35.00,0.00)(49.00,10.00)
\psline(42.00,0.00)(49.00,5.00)
\psline(49.00,0.00)(49.00,0.00)
\pspolygon[fillstyle=solid,linewidth=1pt,fillcolor=lightblue](49.00,35.00)(49.00,0.00)(74.00,0.00)
\psline(49.00,0.00)(74.00,0.00)
\psline(49.00,7.00)(69.00,7.00)
\psline(49.00,14.00)(64.00,14.00)
\psline(49.00,21.00)(59.00,21.00)
\psline(49.00,28.00)(54.00,28.00)
\psline(49.00,35.00)(49.00,35.00)
\psline(74.00,0.00)(49.00,35.00)
\psline(69.00,0.00)(49.00,28.00)
\psline(64.00,0.00)(49.00,21.00)
\psline(59.00,0.00)(49.00,14.00)
\psline(54.00,0.00)(49.00,7.00)
\psline(49.00,0.00)(49.00,0.00)
\psline(49.00,35.00)(49.00,0.00)
\psline(54.00,28.00)(54.00,0.00)
\psline(59.00,21.00)(59.00,0.00)
\psline(64.00,14.00)(64.00,0.00)
\psline(69.00,7.00)(69.00,0.00)
\psline(74.00,0.00)(74.00,0.00)
\endpspicture}

\def\FigurePythagorean{
\hskip 2in
\pspicture(2.0,1.5)
\psset{unit=0.2cm}
\newrgbcolor{lightblue}{0.8 0.8 1}
\newrgbcolor{pink}{1 0.8 0.8}
\newrgbcolor{lightgreen}{0.8 1 0.8}
\newrgbcolor{lightyellow}{1 1 0.8}
\pspolygon[fillstyle=solid,linewidth=1pt,fillcolor=lightyellow](0.00,0.00)(15.00,20.00)(15.00,0.00)
\psline(15.00,20.00)(15.00,0.00)
\psline(12.00,16.00)(12.00,0.00)
\psline(9.00,12.00)(9.00,0.00)
\psline(6.00,8.00)(6.00,0.00)
\psline(3.00,4.00)(3.00,0.00)
\psline(0.00,0.00)(0.00,0.00)
\psline(15.00,0.00)(0.00,0.00)
\psline(15.00,4.00)(3.00,4.00)
\psline(15.00,8.00)(6.00,8.00)
\psline(15.00,12.00)(9.00,12.00)
\psline(15.00,16.00)(12.00,16.00)
\psline(15.00,20.00)(15.00,20.00)
\psline(0.00,0.00)(15.00,20.00)
\psline(3.00,0.00)(15.00,16.00)
\psline(6.00,0.00)(15.00,12.00)
\psline(9.00,0.00)(15.00,8.00)
\psline(12.00,0.00)(15.00,4.00)
\psline(15.00,0.00)(15.00,0.00)
\pspolygon[fillstyle=solid,linewidth=1pt,fillcolor=lightblue](15.00,20.00)(15.00,0.00)(24.60,7.20)
\psline(15.00,0.00)(24.60,7.20)
\psline(15.00,5.00)(22.20,10.40)
\psline(15.00,10.00)(19.80,13.60)
\psline(15.00,15.00)(17.40,16.80)
\psline(15.00,20.00)(15.00,20.00)
\psline(24.60,7.20)(15.00,20.00)
\psline(22.20,5.40)(15.00,15.00)
\psline(19.80,3.60)(15.00,10.00)
\psline(17.40,1.80)(15.00,5.00)
\psline(15.00,0.00)(15.00,0.00)
\psline(15.00,20.00)(15.00,0.00)
\psline(17.40,16.80)(17.40,1.80)
\psline(19.80,13.60)(19.80,3.60)
\psline(22.20,10.40)(22.20,5.40)
\psline(24.60,7.20)(24.60,7.20)
\pspolygon[fillstyle=solid,linewidth=1pt,fillcolor=lightgreen](15.00,0.00)(24.60,7.20)(30.00,0.00)
\psline(24.60,7.20)(30.00,0.00)
\psline(21.40,4.80)(25.00,0.00)
\psline(18.20,2.40)(20.00,0.00)
\psline(15.00,0.00)(15.00,0.00)
\psline(30.00,0.00)(15.00,0.00)
\psline(28.20,2.40)(18.20,2.40)
\psline(26.40,4.80)(21.40,4.80)
\psline(24.60,7.20)(24.60,7.20)
\psline(15.00,0.00)(24.60,7.20)
\psline(20.00,0.00)(26.40,4.80)
\psline(25.00,0.00)(28.20,2.40)
\psline(30.00,0.00)(30.00,0.00)
\endpspicture}

\def\FigureIsoscelesCounterexample{
\hskip2.4in
\pspicture(2.0,1.5)
\psset{unit=1cm}
\newrgbcolor{lightblue}{0.8 0.8 1}
\newrgbcolor{pink}{1 0.8 0.8}
\newrgbcolor{lightgreen}{0.8 1 0.8}
\newrgbcolor{lightyellow}{1 1 0.8}
\pspolygon[fillstyle=solid,linewidth=1pt,fillcolor=lightblue](0.00,0.00)(1.00,2.00)(4.00,0.00)\pspolygon[fillstyle=solid,linewidth=1pt,fillcolor=lightblue](1.00,2.00)(2.00,4.00)(2.00,2.00)\pspolygon[fillstyle=solid,linewidth=1pt,fillcolor=lightblue](2.00,4.00)(2.00,2.00)(3.00,2.00)\pspolygon[fillstyle=solid,linewidth=1pt,fillcolor=lightblue](3.00,2.00)(3.00,0.00)(4.00,0.00)\pspolygon[fillstyle=solid,linewidth=1pt,fillcolor=lightyellow](1.00,2.00)(3.00,2.00)(1.00,1.00)\pspolygon[fillstyle=solid,linewidth=1pt,fillcolor=lightyellow](1.00,1.00)(3.00,2.00)(3.00,1.00)\pspolygon[fillstyle=solid,linewidth=1pt,fillcolor=lightgreen](1.00,1.00)(3.00,1.00)(3.00,0.00)\pspolygon[fillstyle=solid,linewidth=1pt,fillcolor=lightgreen](1.00,1.00)(1.00,0.00)(3.00,0.00)
\endpspicture
}

\def\FigureBigQuadratic{
\hskip0.5in
\pspicture(2.0,1.8)
\psset{unit=1cm}
\newrgbcolor{lightblue}{0.8 0.8 1}
\newrgbcolor{pink}{1 0.8 0.8}
\newrgbcolor{lightgreen}{0.8 1 0.8}
\newrgbcolor{lightyellow}{1 1 0.8}
\pspolygon[fillstyle=solid,linewidth=1pt,fillcolor=lightblue](1.00,0.00)(3.00,5.00)(12.00,0.00)
\psline(3.00,5.00)(12.00,0.00)
\psline(2.80,4.50)(10.90,0.00)
\psline(2.60,4.00)(9.80,0.00)
\psline(2.40,3.50)(8.70,0.00)
\psline(2.20,3.00)(7.60,0.00)
\psline(2.00,2.50)(6.50,0.00)
\psline(1.80,2.00)(5.40,0.00)
\psline(1.60,1.50)(4.30,0.00)
\psline(1.40,1.00)(3.20,0.00)
\psline(1.20,0.50)(2.10,0.00)
\psline(1.00,0.00)(1.00,0.00)
\psline(12.00,0.00)(1.00,0.00)
\psline(11.10,0.50)(1.20,0.50)
\psline(10.20,1.00)(1.40,1.00)
\psline(9.30,1.50)(1.60,1.50)
\psline(8.40,2.00)(1.80,2.00)
\psline(7.50,2.50)(2.00,2.50)
\psline(6.60,3.00)(2.20,3.00)
\psline(5.70,3.50)(2.40,3.50)
\psline(4.80,4.00)(2.60,4.00)
\psline(3.90,4.50)(2.80,4.50)
\psline(3.00,5.00)(3.00,5.00)
\psline(1.00,0.00)(3.00,5.00)
\psline(2.10,0.00)(3.90,4.50)
\psline(3.20,0.00)(4.80,4.00)
\psline(4.30,0.00)(5.70,3.50)
\psline(5.40,0.00)(6.60,3.00)
\psline(6.50,0.00)(7.50,2.50)
\psline(7.60,0.00)(8.40,2.00)
\psline(8.70,0.00)(9.30,1.50)
\psline(9.80,0.00)(10.20,1.00)
\psline(10.90,0.00)(11.10,0.50)
\psline(12.00,0.00)(12.00,0.00)
\endpspicture
}

\def\FigureNonQuadratic{
\hskip1.5in
\pspicture(2.0,2.0)
\psset{unit=0.4cm}
\newrgbcolor{lightblue}{0.8 0.8 1}
\newrgbcolor{pink}{1 0.8 0.8}
\newrgbcolor{lightgreen}{0.8 1 0.8}
\newrgbcolor{lightyellow}{1 1 0.8}
\pspolygon[fillstyle=solid,linewidth=1pt,fillcolor=lightblue](1.00,0.00)(8.50,12.99)(11.00,0.00)
\psline(8.50,12.99)(11.00,0.00)
\psline(7.00,10.39)(9.00,0.00)
\psline(5.50,7.79)(7.00,0.00)
\psline(4.00,5.20)(5.00,0.00)
\psline(2.50,2.60)(3.00,0.00)
\psline(1.00,0.00)(1.00,0.00)
\psline(11.00,0.00)(1.00,0.00)
\psline(10.50,2.60)(2.50,2.60)
\psline(10.00,5.20)(4.00,5.20)
\psline(9.50,7.79)(5.50,7.79)
\psline(9.00,10.39)(7.00,10.39)
\psline(8.50,12.99)(8.50,12.99)
\psline(1.00,0.00)(8.50,12.99)
\psline(3.00,0.00)(9.00,10.39)
\psline(5.00,0.00)(9.50,7.79)
\psline(7.00,0.00)(10.00,5.20)
\psline(9.00,0.00)(10.50,2.60)
\psline(11.00,0.00)(11.00,0.00)
\pspolygon[fillstyle=solid,linewidth=1pt, fillcolor=lightyellow](1.00,0.00)(7.00,0.00)(10.00,5.20)(4.00,5.20)\pspolygon[fillstyle=solid,linewidth=1pt, fillcolor=lightyellow](1.00,0.00)(4.00,0.00)(2.00,1.73)\pspolygon[fillstyle=solid,linewidth=1pt, fillcolor=lightyellow](2.00,1.73)(5.00,1.73)(3.00,3.46)\pspolygon[fillstyle=solid,linewidth=1pt, fillcolor=lightyellow](3.00,3.46)(6.00,3.46)(4.00,5.20)\pspolygon[fillstyle=solid,linewidth=1pt, fillcolor=lightyellow](4.00,0.00)(7.00,0.00)(5.00,1.73)\pspolygon[fillstyle=solid,linewidth=1pt, fillcolor=lightyellow](5.00,1.73)(8.00,1.73)(6.00,3.46)\pspolygon[fillstyle=solid,linewidth=1pt, fillcolor=lightyellow](6.00,3.46)(9.00,3.46)(7.00,5.20)
\endpspicture
}

\def\FigureThreeZeroTwoFour{
\hskip 1in
\pspicture(2.0,2.8)
\psset{unit=0.8cm}
\newrgbcolor{lightblue}{0.8 0.8 1}
\newrgbcolor{pink}{1 0.8 0.8}
\newrgbcolor{lightgreen}{0.8 1 0.8}
\newrgbcolor{lightyellow}{1 1 0.8}
\pspolygon[fillstyle=solid,linewidth=1pt,fillcolor=lightblue](1.63,2.96)(0.00,0.00)(7.00,0.00)
\psline(0.00,0.00)(7.00,0.00)
\psline(0.23,0.42)(6.23,0.42)
\psline(0.47,0.85)(5.47,0.85)
\psline(0.70,1.27)(4.70,1.27)
\psline(0.93,1.69)(3.93,1.69)
\psline(1.17,2.12)(3.17,2.12)
\psline(1.40,2.54)(2.40,2.54)
\psline(1.63,2.96)(1.63,2.96)
\psline(7.00,0.00)(1.63,2.96)
\psline(6.00,0.00)(1.40,2.54)
\psline(5.00,0.00)(1.17,2.12)
\psline(4.00,0.00)(0.93,1.69)
\psline(3.00,0.00)(0.70,1.27)
\psline(2.00,0.00)(0.47,0.85)
\psline(1.00,0.00)(0.23,0.42)
\psline(0.00,0.00)(0.00,0.00)
\psline(1.63,2.96)(0.00,0.00)
\psline(2.40,2.54)(1.00,0.00)
\psline(3.17,2.12)(2.00,0.00)
\psline(3.93,1.69)(3.00,0.00)
\psline(4.70,1.27)(4.00,0.00)
\psline(5.47,0.85)(5.00,0.00)
\psline(6.23,0.42)(6.00,0.00)
\psline(7.00,0.00)(7.00,0.00)
\pspolygon[fillstyle=solid,linewidth=1pt,fillcolor=pink](0.93,1.69)(1.70,1.27)(1.47,0.85)
\psline(1.70,1.27)(1.47,0.85)
\psline(0.93,1.69)(0.93,1.69)
\psline(1.47,0.85)(0.93,1.69)
\psline(1.70,1.27)(1.70,1.27)
\psline(0.93,1.69)(1.70,1.27)
\psline(1.47,0.85)(1.47,0.85)
\pspolygon[fillstyle=solid,linewidth=1pt,fillcolor=pink](0.93,1.69)(1.47,0.85)(0.70,1.27)
\psline(1.47,0.85)(0.70,1.27)
\psline(0.93,1.69)(0.93,1.69)
\psline(0.70,1.27)(0.93,1.69)
\psline(1.47,0.85)(1.47,0.85)
\psline(0.93,1.69)(1.47,0.85)
\psline(0.70,1.27)(0.70,1.27)
\pspolygon[fillstyle=solid,linewidth=1pt,fillcolor=pink](1.93,1.69)(2.70,1.27)(2.47,0.85)
\psline(2.70,1.27)(2.47,0.85)
\psline(1.93,1.69)(1.93,1.69)
\psline(2.47,0.85)(1.93,1.69)
\psline(2.70,1.27)(2.70,1.27)
\psline(1.93,1.69)(2.70,1.27)
\psline(2.47,0.85)(2.47,0.85)
\pspolygon[fillstyle=solid,linewidth=1pt,fillcolor=pink](1.93,1.69)(2.47,0.85)(1.70,1.27)
\psline(2.47,0.85)(1.70,1.27)
\psline(1.93,1.69)(1.93,1.69)
\psline(1.70,1.27)(1.93,1.69)
\psline(2.47,0.85)(2.47,0.85)
\psline(1.93,1.69)(2.47,0.85)
\psline(1.70,1.27)(1.70,1.27)
\pspolygon[fillstyle=solid,linewidth=1pt,fillcolor=pink](2.93,1.69)(3.70,1.27)(3.47,0.85)
\psline(3.70,1.27)(3.47,0.85)
\psline(2.93,1.69)(2.93,1.69)
\psline(3.47,0.85)(2.93,1.69)
\psline(3.70,1.27)(3.70,1.27)
\psline(2.93,1.69)(3.70,1.27)
\psline(3.47,0.85)(3.47,0.85)
\pspolygon[fillstyle=solid,linewidth=1pt,fillcolor=pink](2.93,1.69)(3.47,0.85)(2.70,1.27)
\psline(3.47,0.85)(2.70,1.27)
\psline(2.93,1.69)(2.93,1.69)
\psline(2.70,1.27)(2.93,1.69)
\psline(3.47,0.85)(3.47,0.85)
\psline(2.93,1.69)(3.47,0.85)
\psline(2.70,1.27)(2.70,1.27)
\pspolygon[fillstyle=solid,linewidth=1pt,fillcolor=pink](3.47,0.85)(4.23,0.42)(4.00,0.00)
\psline(4.23,0.42)(4.00,0.00)
\psline(3.47,0.85)(3.47,0.85)
\psline(4.00,0.00)(3.47,0.85)
\psline(4.23,0.42)(4.23,0.42)
\psline(3.47,0.85)(4.23,0.42)
\psline(4.00,0.00)(4.00,0.00)
\pspolygon[fillstyle=solid,linewidth=1pt,fillcolor=pink](3.47,0.85)(4.00,0.00)(3.23,0.42)
\psline(4.00,0.00)(3.23,0.42)
\psline(3.47,0.85)(3.47,0.85)
\psline(3.23,0.42)(3.47,0.85)
\psline(4.00,0.00)(4.00,0.00)
\psline(3.47,0.85)(4.00,0.00)
\psline(3.23,0.42)(3.23,0.42)
\psline[linestyle=dashed](1.63,2.96)(2.63,2.96)
\psline[linestyle=dashed](1.87,3.38)(2.40,2.54)
\psline(5.60,10.15)(0.00,0.00)
\psline(0.00,0.00)(11.21,0.00)
\psline(5.60,10.15)(11.21,0.00)
\psline(2.40,2.54)(2.63,2.96)
\psline(1.87,3.38)(2.63,2.96)
\psdot(5.60,10.15)\put(5.04,9.86){$A$}
\psdot(0.00,0.00)\put(0.00,-0.44){$B$}
\psdot(11.21,0.00)\put(11.21,-0.44){$C$}
\psdot(7.00,0.00)\put(7.00,-0.44){$U$}
\psdot(1.63,2.96)\put(1.03,2.95){$Q$}
\psdot(1.87,3.38)\put(1.27,3.37){$P$}
\psdot(2.40,2.54)\put(2.59,2.44){$R$}
\psdot(2.63,2.96)\put(2.83,2.86){$S$}
\endpspicture
}

\def\FigureThreeZeroFourEight{
\hskip 1in
\pspicture(2.0,2.8)
\psset{unit=0.8cm}
\newrgbcolor{lightblue}{0.8 0.8 1}
\newrgbcolor{pink}{1 0.8 0.8}
\newrgbcolor{lightgreen}{0.8 1 0.8}
\newrgbcolor{lightyellow}{1 1 0.8}
\pspolygon[fillstyle=solid,linewidth=1pt,fillcolor=lightblue](1.63,2.96)(0.00,0.00)(7.00,0.00)
\psline(0.00,0.00)(7.00,0.00)
\psline(0.23,0.42)(6.23,0.42)
\psline(0.47,0.85)(5.47,0.85)
\psline(0.70,1.27)(4.70,1.27)
\psline(0.93,1.69)(3.93,1.69)
\psline(1.17,2.12)(3.17,2.12)
\psline(1.40,2.54)(2.40,2.54)
\psline(1.63,2.96)(1.63,2.96)
\psline(7.00,0.00)(1.63,2.96)
\psline(6.00,0.00)(1.40,2.54)
\psline(5.00,0.00)(1.17,2.12)
\psline(4.00,0.00)(0.93,1.69)
\psline(3.00,0.00)(0.70,1.27)
\psline(2.00,0.00)(0.47,0.85)
\psline(1.00,0.00)(0.23,0.42)
\psline(0.00,0.00)(0.00,0.00)
\psline(1.63,2.96)(0.00,0.00)
\psline(2.40,2.54)(1.00,0.00)
\psline(3.17,2.12)(2.00,0.00)
\psline(3.93,1.69)(3.00,0.00)
\psline(4.70,1.27)(4.00,0.00)
\psline(5.47,0.85)(5.00,0.00)
\psline(6.23,0.42)(6.00,0.00)
\psline(7.00,0.00)(7.00,0.00)
\pspolygon[fillstyle=solid,linewidth=1pt,fillcolor=pink](0.93,1.69)(1.70,1.27)(1.47,0.85)
\psline(1.70,1.27)(1.47,0.85)
\psline(0.93,1.69)(0.93,1.69)
\psline(1.47,0.85)(0.93,1.69)
\psline(1.70,1.27)(1.70,1.27)
\psline(0.93,1.69)(1.70,1.27)
\psline(1.47,0.85)(1.47,0.85)
\pspolygon[fillstyle=solid,linewidth=1pt,fillcolor=pink](0.93,1.69)(1.47,0.85)(0.70,1.27)
\psline(1.47,0.85)(0.70,1.27)
\psline(0.93,1.69)(0.93,1.69)
\psline(0.70,1.27)(0.93,1.69)
\psline(1.47,0.85)(1.47,0.85)
\psline(0.93,1.69)(1.47,0.85)
\psline(0.70,1.27)(0.70,1.27)
\pspolygon[fillstyle=solid,linewidth=1pt,fillcolor=pink](1.93,1.69)(2.70,1.27)(2.47,0.85)
\psline(2.70,1.27)(2.47,0.85)
\psline(1.93,1.69)(1.93,1.69)
\psline(2.47,0.85)(1.93,1.69)
\psline(2.70,1.27)(2.70,1.27)
\psline(1.93,1.69)(2.70,1.27)
\psline(2.47,0.85)(2.47,0.85)
\pspolygon[fillstyle=solid,linewidth=1pt,fillcolor=pink](1.93,1.69)(2.47,0.85)(1.70,1.27)
\psline(2.47,0.85)(1.70,1.27)
\psline(1.93,1.69)(1.93,1.69)
\psline(1.70,1.27)(1.93,1.69)
\psline(2.47,0.85)(2.47,0.85)
\psline(1.93,1.69)(2.47,0.85)
\psline(1.70,1.27)(1.70,1.27)
\pspolygon[fillstyle=solid,linewidth=1pt,fillcolor=pink](2.93,1.69)(3.70,1.27)(3.47,0.85)
\psline(3.70,1.27)(3.47,0.85)
\psline(2.93,1.69)(2.93,1.69)
\psline(3.47,0.85)(2.93,1.69)
\psline(3.70,1.27)(3.70,1.27)
\psline(2.93,1.69)(3.70,1.27)
\psline(3.47,0.85)(3.47,0.85)
\pspolygon[fillstyle=solid,linewidth=1pt,fillcolor=pink](2.93,1.69)(3.47,0.85)(2.70,1.27)
\psline(3.47,0.85)(2.70,1.27)
\psline(2.93,1.69)(2.93,1.69)
\psline(2.70,1.27)(2.93,1.69)
\psline(3.47,0.85)(3.47,0.85)
\psline(2.93,1.69)(3.47,0.85)
\psline(2.70,1.27)(2.70,1.27)
\pspolygon[fillstyle=solid,linewidth=1pt,fillcolor=pink](3.47,0.85)(4.23,0.42)(4.00,0.00)
\psline(4.23,0.42)(4.00,0.00)
\psline(3.47,0.85)(3.47,0.85)
\psline(4.00,0.00)(3.47,0.85)
\psline(4.23,0.42)(4.23,0.42)
\psline(3.47,0.85)(4.23,0.42)
\psline(4.00,0.00)(4.00,0.00)
\pspolygon[fillstyle=solid,linewidth=1pt,fillcolor=pink](3.47,0.85)(4.00,0.00)(3.23,0.42)
\psline(4.00,0.00)(3.23,0.42)
\psline(3.47,0.85)(3.47,0.85)
\psline(3.23,0.42)(3.47,0.85)
\psline(4.00,0.00)(4.00,0.00)
\psline(3.47,0.85)(4.00,0.00)
\psline(3.23,0.42)(3.23,0.42)
\psline(5.60,10.15)(0.00,0.00)
\psline(0.00,0.00)(11.21,0.00)
\psline(5.60,10.15)(11.21,0.00)
\pspolygon[fillstyle=solid,linewidth=1pt,fillcolor=lightyellow](2.40,2.54)(2.63,2.96)(1.63,2.96)
\psline(2.63,2.96)(1.63,2.96)
\psline(2.40,2.54)(2.40,2.54)
\psline(1.63,2.96)(2.40,2.54)
\psline(2.63,2.96)(2.63,2.96)
\psline(2.40,2.54)(2.63,2.96)
\psline(1.63,2.96)(1.63,2.96)
\pspolygon[fillstyle=solid,linewidth=1pt,fillcolor=lightyellow](1.63,2.96)(1.87,3.38)(2.63,2.96)
\psline(1.87,3.38)(2.63,2.96)
\psline(1.63,2.96)(1.63,2.96)
\psline(2.63,2.96)(1.63,2.96)
\psline(1.87,3.38)(1.87,3.38)
\psline(1.63,2.96)(1.87,3.38)
\psline(2.63,2.96)(2.63,2.96)
\pspolygon[fillstyle=solid,linewidth=1pt,fillcolor=lightyellow](3.17,2.12)(3.40,2.54)(2.40,2.54)
\psline(3.40,2.54)(2.40,2.54)
\psline(3.17,2.12)(3.17,2.12)
\psline(2.40,2.54)(3.17,2.12)
\psline(3.40,2.54)(3.40,2.54)
\psline(3.17,2.12)(3.40,2.54)
\psline(2.40,2.54)(2.40,2.54)
\pspolygon[fillstyle=solid,linewidth=1pt,fillcolor=lightyellow](2.40,2.54)(2.63,2.96)(3.40,2.54)
\psline(2.63,2.96)(3.40,2.54)
\psline(2.40,2.54)(2.40,2.54)
\psline(3.40,2.54)(2.40,2.54)
\psline(2.63,2.96)(2.63,2.96)
\psline(2.40,2.54)(2.63,2.96)
\psline(3.40,2.54)(3.40,2.54)
\pspolygon[fillstyle=solid,linewidth=1pt,fillcolor=lightyellow](3.93,1.69)(4.17,2.12)(3.17,2.12)
\psline(4.17,2.12)(3.17,2.12)
\psline(3.93,1.69)(3.93,1.69)
\psline(3.17,2.12)(3.93,1.69)
\psline(4.17,2.12)(4.17,2.12)
\psline(3.93,1.69)(4.17,2.12)
\psline(3.17,2.12)(3.17,2.12)
\pspolygon[fillstyle=solid,linewidth=1pt,fillcolor=lightyellow](3.17,2.12)(3.40,2.54)(4.17,2.12)
\psline(3.40,2.54)(4.17,2.12)
\psline(3.17,2.12)(3.17,2.12)
\psline(4.17,2.12)(3.17,2.12)
\psline(3.40,2.54)(3.40,2.54)
\psline(3.17,2.12)(3.40,2.54)
\psline(4.17,2.12)(4.17,2.12)
\pspolygon[fillstyle=solid,linewidth=1pt,fillcolor=lightyellow](4.70,1.27)(4.93,1.69)(3.93,1.69)
\psline(4.93,1.69)(3.93,1.69)
\psline(4.70,1.27)(4.70,1.27)
\psline(3.93,1.69)(4.70,1.27)
\psline(4.93,1.69)(4.93,1.69)
\psline(4.70,1.27)(4.93,1.69)
\psline(3.93,1.69)(3.93,1.69)
\pspolygon[fillstyle=solid,linewidth=1pt,fillcolor=lightyellow](3.93,1.69)(4.17,2.12)(4.93,1.69)
\psline(4.17,2.12)(4.93,1.69)
\psline(3.93,1.69)(3.93,1.69)
\psline(4.93,1.69)(3.93,1.69)
\psline(4.17,2.12)(4.17,2.12)
\psline(3.93,1.69)(4.17,2.12)
\psline(4.93,1.69)(4.93,1.69)
\pspolygon[fillstyle=solid,linewidth=1pt,fillcolor=lightyellow](5.47,0.85)(5.70,1.27)(4.70,1.27)
\psline(5.70,1.27)(4.70,1.27)
\psline(5.47,0.85)(5.47,0.85)
\psline(4.70,1.27)(5.47,0.85)
\psline(5.70,1.27)(5.70,1.27)
\psline(5.47,0.85)(5.70,1.27)
\psline(4.70,1.27)(4.70,1.27)
\pspolygon[fillstyle=solid,linewidth=1pt,fillcolor=lightyellow](4.70,1.27)(4.93,1.69)(5.70,1.27)
\psline(4.93,1.69)(5.70,1.27)
\psline(4.70,1.27)(4.70,1.27)
\psline(5.70,1.27)(4.70,1.27)
\psline(4.93,1.69)(4.93,1.69)
\psline(4.70,1.27)(4.93,1.69)
\psline(5.70,1.27)(5.70,1.27)
\pspolygon[fillstyle=solid,linewidth=1pt,fillcolor=lightblue](5.47,0.85)(6.47,0.85)(6.23,0.42)
\psline(6.47,0.85)(6.23,0.42)
\psline(5.47,0.85)(5.47,0.85)
\psline(6.23,0.42)(5.47,0.85)
\psline(6.47,0.85)(6.47,0.85)
\psline(5.47,0.85)(6.47,0.85)
\psline(6.23,0.42)(6.23,0.42)
\put(6.06,0.51){$5$}
\pspolygon[fillstyle=solid,linewidth=1pt,fillcolor=lightgreen](5.47,0.85)(5.95,0.85)(5.95,1.72)
\psline(5.95,0.85)(5.95,1.72)
\psline(5.47,0.85)(5.47,0.85)
\psline(5.95,1.72)(5.47,0.85)
\psline(5.95,0.85)(5.95,0.85)
\psline(5.47,0.85)(5.95,0.85)
\psline(5.95,1.72)(5.95,1.72)
\put(5.69,0.90){$6$}
\psdot(5.60,10.15)\put(5.04,9.86){$A$}
\psdot(0.00,0.00)\put(0.00,-0.44){$B$}
\psdot(11.21,0.00)\put(11.21,-0.44){$C$}
\psdot(7.00,0.00)\put(7.00,-0.44){$U$}
\psdot(1.63,2.96)\put(1.03,2.95){$Q$}
\psdot(1.87,3.38)\put(1.27,3.37){$P$}
\psdot(5.70,1.27)\put(5.26,1.58){$S_k$}
\psdot(6.47,0.85)\put(6.66,0.74){$S_{k+1}$}
\psdot(6.23,0.42)\put(6.43,0.32){$R_{k+1}$}
\psline(5.70,1.27)(8.12,5.65)
\endpspicture
}

\def\FigureThreeZeroSixEight{
\hskip 1in
\pspicture(2.0,2.8)
\psset{unit=0.8cm}
\newrgbcolor{lightblue}{0.8 0.8 1}
\newrgbcolor{pink}{1 0.8 0.8}
\newrgbcolor{lightgreen}{0.8 1 0.8}
\newrgbcolor{lightyellow}{1 1 0.8}
\pspolygon[fillstyle=solid,linewidth=1pt,fillcolor=lightblue](1.63,2.96)(0.00,0.00)(7.00,0.00)
\psline(0.00,0.00)(7.00,0.00)
\psline(0.23,0.42)(6.23,0.42)
\psline(0.47,0.85)(5.47,0.85)
\psline(0.70,1.27)(4.70,1.27)
\psline(0.93,1.69)(3.93,1.69)
\psline(1.17,2.12)(3.17,2.12)
\psline(1.40,2.54)(2.40,2.54)
\psline(1.63,2.96)(1.63,2.96)
\psline(7.00,0.00)(1.63,2.96)
\psline(6.00,0.00)(1.40,2.54)
\psline(5.00,0.00)(1.17,2.12)
\psline(4.00,0.00)(0.93,1.69)
\psline(3.00,0.00)(0.70,1.27)
\psline(2.00,0.00)(0.47,0.85)
\psline(1.00,0.00)(0.23,0.42)
\psline(0.00,0.00)(0.00,0.00)
\psline(1.63,2.96)(0.00,0.00)
\psline(2.40,2.54)(1.00,0.00)
\psline(3.17,2.12)(2.00,0.00)
\psline(3.93,1.69)(3.00,0.00)
\psline(4.70,1.27)(4.00,0.00)
\psline(5.47,0.85)(5.00,0.00)
\psline(6.23,0.42)(6.00,0.00)
\psline(7.00,0.00)(7.00,0.00)
\pspolygon[fillstyle=solid,linewidth=1pt,fillcolor=pink](0.93,1.69)(1.70,1.27)(1.47,0.85)
\psline(1.70,1.27)(1.47,0.85)
\psline(0.93,1.69)(0.93,1.69)
\psline(1.47,0.85)(0.93,1.69)
\psline(1.70,1.27)(1.70,1.27)
\psline(0.93,1.69)(1.70,1.27)
\psline(1.47,0.85)(1.47,0.85)
\pspolygon[fillstyle=solid,linewidth=1pt,fillcolor=pink](0.93,1.69)(1.47,0.85)(0.70,1.27)
\psline(1.47,0.85)(0.70,1.27)
\psline(0.93,1.69)(0.93,1.69)
\psline(0.70,1.27)(0.93,1.69)
\psline(1.47,0.85)(1.47,0.85)
\psline(0.93,1.69)(1.47,0.85)
\psline(0.70,1.27)(0.70,1.27)
\pspolygon[fillstyle=solid,linewidth=1pt,fillcolor=pink](1.93,1.69)(2.70,1.27)(2.47,0.85)
\psline(2.70,1.27)(2.47,0.85)
\psline(1.93,1.69)(1.93,1.69)
\psline(2.47,0.85)(1.93,1.69)
\psline(2.70,1.27)(2.70,1.27)
\psline(1.93,1.69)(2.70,1.27)
\psline(2.47,0.85)(2.47,0.85)
\pspolygon[fillstyle=solid,linewidth=1pt,fillcolor=pink](1.93,1.69)(2.47,0.85)(1.70,1.27)
\psline(2.47,0.85)(1.70,1.27)
\psline(1.93,1.69)(1.93,1.69)
\psline(1.70,1.27)(1.93,1.69)
\psline(2.47,0.85)(2.47,0.85)
\psline(1.93,1.69)(2.47,0.85)
\psline(1.70,1.27)(1.70,1.27)
\pspolygon[fillstyle=solid,linewidth=1pt,fillcolor=pink](2.93,1.69)(3.70,1.27)(3.47,0.85)
\psline(3.70,1.27)(3.47,0.85)
\psline(2.93,1.69)(2.93,1.69)
\psline(3.47,0.85)(2.93,1.69)
\psline(3.70,1.27)(3.70,1.27)
\psline(2.93,1.69)(3.70,1.27)
\psline(3.47,0.85)(3.47,0.85)
\pspolygon[fillstyle=solid,linewidth=1pt,fillcolor=pink](2.93,1.69)(3.47,0.85)(2.70,1.27)
\psline(3.47,0.85)(2.70,1.27)
\psline(2.93,1.69)(2.93,1.69)
\psline(2.70,1.27)(2.93,1.69)
\psline(3.47,0.85)(3.47,0.85)
\psline(2.93,1.69)(3.47,0.85)
\psline(2.70,1.27)(2.70,1.27)
\pspolygon[fillstyle=solid,linewidth=1pt,fillcolor=pink](3.47,0.85)(4.23,0.42)(4.00,0.00)
\psline(4.23,0.42)(4.00,0.00)
\psline(3.47,0.85)(3.47,0.85)
\psline(4.00,0.00)(3.47,0.85)
\psline(4.23,0.42)(4.23,0.42)
\psline(3.47,0.85)(4.23,0.42)
\psline(4.00,0.00)(4.00,0.00)
\pspolygon[fillstyle=solid,linewidth=1pt,fillcolor=pink](3.47,0.85)(4.00,0.00)(3.23,0.42)
\psline(4.00,0.00)(3.23,0.42)
\psline(3.47,0.85)(3.47,0.85)
\psline(3.23,0.42)(3.47,0.85)
\psline(4.00,0.00)(4.00,0.00)
\psline(3.47,0.85)(4.00,0.00)
\psline(3.23,0.42)(3.23,0.42)
\psline(5.60,10.15)(0.00,0.00)
\psline(0.00,0.00)(11.21,0.00)
\psline(5.60,10.15)(11.21,0.00)
\pspolygon[fillstyle=solid,linewidth=1pt,fillcolor=lightyellow](2.40,2.54)(2.63,2.96)(1.63,2.96)
\psline(2.63,2.96)(1.63,2.96)
\psline(2.40,2.54)(2.40,2.54)
\psline(1.63,2.96)(2.40,2.54)
\psline(2.63,2.96)(2.63,2.96)
\psline(2.40,2.54)(2.63,2.96)
\psline(1.63,2.96)(1.63,2.96)
\pspolygon[fillstyle=solid,linewidth=1pt,fillcolor=lightyellow](1.63,2.96)(1.87,3.38)(2.63,2.96)
\psline(1.87,3.38)(2.63,2.96)
\psline(1.63,2.96)(1.63,2.96)
\psline(2.63,2.96)(1.63,2.96)
\psline(1.87,3.38)(1.87,3.38)
\psline(1.63,2.96)(1.87,3.38)
\psline(2.63,2.96)(2.63,2.96)
\pspolygon[fillstyle=solid,linewidth=1pt,fillcolor=lightyellow](3.17,2.12)(3.40,2.54)(2.40,2.54)
\psline(3.40,2.54)(2.40,2.54)
\psline(3.17,2.12)(3.17,2.12)
\psline(2.40,2.54)(3.17,2.12)
\psline(3.40,2.54)(3.40,2.54)
\psline(3.17,2.12)(3.40,2.54)
\psline(2.40,2.54)(2.40,2.54)
\pspolygon[fillstyle=solid,linewidth=1pt,fillcolor=lightyellow](2.40,2.54)(2.63,2.96)(3.40,2.54)
\psline(2.63,2.96)(3.40,2.54)
\psline(2.40,2.54)(2.40,2.54)
\psline(3.40,2.54)(2.40,2.54)
\psline(2.63,2.96)(2.63,2.96)
\psline(2.40,2.54)(2.63,2.96)
\psline(3.40,2.54)(3.40,2.54)
\pspolygon[fillstyle=solid,linewidth=1pt,fillcolor=lightyellow](3.93,1.69)(4.17,2.12)(3.17,2.12)
\psline(4.17,2.12)(3.17,2.12)
\psline(3.93,1.69)(3.93,1.69)
\psline(3.17,2.12)(3.93,1.69)
\psline(4.17,2.12)(4.17,2.12)
\psline(3.93,1.69)(4.17,2.12)
\psline(3.17,2.12)(3.17,2.12)
\pspolygon[fillstyle=solid,linewidth=1pt,fillcolor=lightyellow](3.17,2.12)(3.40,2.54)(4.17,2.12)
\psline(3.40,2.54)(4.17,2.12)
\psline(3.17,2.12)(3.17,2.12)
\psline(4.17,2.12)(3.17,2.12)
\psline(3.40,2.54)(3.40,2.54)
\psline(3.17,2.12)(3.40,2.54)
\psline(4.17,2.12)(4.17,2.12)
\pspolygon[fillstyle=solid,linewidth=1pt,fillcolor=lightyellow](4.70,1.27)(4.93,1.69)(3.93,1.69)
\psline(4.93,1.69)(3.93,1.69)
\psline(4.70,1.27)(4.70,1.27)
\psline(3.93,1.69)(4.70,1.27)
\psline(4.93,1.69)(4.93,1.69)
\psline(4.70,1.27)(4.93,1.69)
\psline(3.93,1.69)(3.93,1.69)
\pspolygon[fillstyle=solid,linewidth=1pt,fillcolor=lightyellow](3.93,1.69)(4.17,2.12)(4.93,1.69)
\psline(4.17,2.12)(4.93,1.69)
\psline(3.93,1.69)(3.93,1.69)
\psline(4.93,1.69)(3.93,1.69)
\psline(4.17,2.12)(4.17,2.12)
\psline(3.93,1.69)(4.17,2.12)
\psline(4.93,1.69)(4.93,1.69)
\pspolygon[fillstyle=solid,linewidth=1pt,fillcolor=lightyellow](5.47,0.85)(5.70,1.27)(4.70,1.27)
\psline(5.70,1.27)(4.70,1.27)
\psline(5.47,0.85)(5.47,0.85)
\psline(4.70,1.27)(5.47,0.85)
\psline(5.70,1.27)(5.70,1.27)
\psline(5.47,0.85)(5.70,1.27)
\psline(4.70,1.27)(4.70,1.27)
\pspolygon[fillstyle=solid,linewidth=1pt,fillcolor=lightyellow](4.70,1.27)(4.93,1.69)(5.70,1.27)
\psline(4.93,1.69)(5.70,1.27)
\psline(4.70,1.27)(4.70,1.27)
\psline(5.70,1.27)(4.70,1.27)
\psline(4.93,1.69)(4.93,1.69)
\psline(4.70,1.27)(4.93,1.69)
\psline(5.70,1.27)(5.70,1.27)
\pspolygon[fillstyle=solid,linewidth=1pt,fillcolor=lightyellow](5.47,0.85)(5.70,1.27)(6.23,0.42)
\psline(5.70,1.27)(6.23,0.42)
\psline(5.47,0.85)(5.47,0.85)
\psline(6.23,0.42)(5.47,0.85)
\psline(5.70,1.27)(5.70,1.27)
\psline(5.47,0.85)(5.70,1.27)
\psline(6.23,0.42)(6.23,0.42)
\put(5.66,0.77){$5$}
\pspolygon[fillstyle=solid,linewidth=1pt,fillcolor=lightgreen](6.23,0.42)(5.98,0.83)(6.72,1.30)
\psline(5.98,0.83)(6.72,1.30)
\psline(6.23,0.42)(6.23,0.42)
\psline(6.72,1.30)(6.23,0.42)
\psline(5.98,0.83)(5.98,0.83)
\psline(6.23,0.42)(5.98,0.83)
\psline(6.72,1.30)(6.72,1.30)
\put(6.10,0.65){$6$}
\psdot(5.60,10.15)\put(5.04,9.86){$A$}
\psdot(0.00,0.00)\put(0.00,-0.44){$B$}
\psdot(11.21,0.00)\put(11.21,-0.44){$C$}
\psdot(7.00,0.00)\put(7.00,-0.44){$U$}
\psdot(1.63,2.96)\put(1.03,2.95){$Q$}
\psdot(1.87,3.38)\put(1.27,3.37){$P$}
\psdot(5.70,1.27)\put(5.26,1.58){$S_k$}
\psdot(6.47,0.85)\put(6.66,0.74){$S_{k+1}$}
\psdot(6.23,0.42)\put(6.43,0.32){$R_{k+1}$}
\psline(5.70,1.27)(3.03,5.50)
\endpspicture
}

\def\FigureThreeFiveSevenEight{
\hskip 1in
\pspicture(2.0,2.8)
\psset{unit=0.8cm}
\newrgbcolor{lightblue}{0.8 0.8 1}
\newrgbcolor{pink}{1 0.8 0.8}
\newrgbcolor{lightgreen}{0.8 1 0.8}
\newrgbcolor{lightyellow}{1 1 0.8}
\pspolygon[fillstyle=solid,linewidth=1pt,fillcolor=lightblue](1.63,2.96)(0.00,0.00)(7.00,0.00)
\psline(0.00,0.00)(7.00,0.00)
\psline(0.23,0.42)(6.23,0.42)
\psline(0.47,0.85)(5.47,0.85)
\psline(0.70,1.27)(4.70,1.27)
\psline(0.93,1.69)(3.93,1.69)
\psline(1.17,2.12)(3.17,2.12)
\psline(1.40,2.54)(2.40,2.54)
\psline(1.63,2.96)(1.63,2.96)
\psline(7.00,0.00)(1.63,2.96)
\psline(6.00,0.00)(1.40,2.54)
\psline(5.00,0.00)(1.17,2.12)
\psline(4.00,0.00)(0.93,1.69)
\psline(3.00,0.00)(0.70,1.27)
\psline(2.00,0.00)(0.47,0.85)
\psline(1.00,0.00)(0.23,0.42)
\psline(0.00,0.00)(0.00,0.00)
\psline(1.63,2.96)(0.00,0.00)
\psline(2.40,2.54)(1.00,0.00)
\psline(3.17,2.12)(2.00,0.00)
\psline(3.93,1.69)(3.00,0.00)
\psline(4.70,1.27)(4.00,0.00)
\psline(5.47,0.85)(5.00,0.00)
\psline(6.23,0.42)(6.00,0.00)
\psline(7.00,0.00)(7.00,0.00)
\pspolygon[fillstyle=solid,linewidth=1pt,fillcolor=pink](0.93,1.69)(1.70,1.27)(1.47,0.85)
\psline(1.70,1.27)(1.47,0.85)
\psline(0.93,1.69)(0.93,1.69)
\psline(1.47,0.85)(0.93,1.69)
\psline(1.70,1.27)(1.70,1.27)
\psline(0.93,1.69)(1.70,1.27)
\psline(1.47,0.85)(1.47,0.85)
\pspolygon[fillstyle=solid,linewidth=1pt,fillcolor=pink](0.93,1.69)(1.47,0.85)(0.70,1.27)
\psline(1.47,0.85)(0.70,1.27)
\psline(0.93,1.69)(0.93,1.69)
\psline(0.70,1.27)(0.93,1.69)
\psline(1.47,0.85)(1.47,0.85)
\psline(0.93,1.69)(1.47,0.85)
\psline(0.70,1.27)(0.70,1.27)
\pspolygon[fillstyle=solid,linewidth=1pt,fillcolor=pink](1.93,1.69)(2.70,1.27)(2.47,0.85)
\psline(2.70,1.27)(2.47,0.85)
\psline(1.93,1.69)(1.93,1.69)
\psline(2.47,0.85)(1.93,1.69)
\psline(2.70,1.27)(2.70,1.27)
\psline(1.93,1.69)(2.70,1.27)
\psline(2.47,0.85)(2.47,0.85)
\pspolygon[fillstyle=solid,linewidth=1pt,fillcolor=pink](1.93,1.69)(2.47,0.85)(1.70,1.27)
\psline(2.47,0.85)(1.70,1.27)
\psline(1.93,1.69)(1.93,1.69)
\psline(1.70,1.27)(1.93,1.69)
\psline(2.47,0.85)(2.47,0.85)
\psline(1.93,1.69)(2.47,0.85)
\psline(1.70,1.27)(1.70,1.27)
\pspolygon[fillstyle=solid,linewidth=1pt,fillcolor=pink](2.93,1.69)(3.70,1.27)(3.47,0.85)
\psline(3.70,1.27)(3.47,0.85)
\psline(2.93,1.69)(2.93,1.69)
\psline(3.47,0.85)(2.93,1.69)
\psline(3.70,1.27)(3.70,1.27)
\psline(2.93,1.69)(3.70,1.27)
\psline(3.47,0.85)(3.47,0.85)
\pspolygon[fillstyle=solid,linewidth=1pt,fillcolor=pink](2.93,1.69)(3.47,0.85)(2.70,1.27)
\psline(3.47,0.85)(2.70,1.27)
\psline(2.93,1.69)(2.93,1.69)
\psline(2.70,1.27)(2.93,1.69)
\psline(3.47,0.85)(3.47,0.85)
\psline(2.93,1.69)(3.47,0.85)
\psline(2.70,1.27)(2.70,1.27)
\pspolygon[fillstyle=solid,linewidth=1pt,fillcolor=pink](3.47,0.85)(4.23,0.42)(4.00,0.00)
\psline(4.23,0.42)(4.00,0.00)
\psline(3.47,0.85)(3.47,0.85)
\psline(4.00,0.00)(3.47,0.85)
\psline(4.23,0.42)(4.23,0.42)
\psline(3.47,0.85)(4.23,0.42)
\psline(4.00,0.00)(4.00,0.00)
\pspolygon[fillstyle=solid,linewidth=1pt,fillcolor=pink](3.47,0.85)(4.00,0.00)(3.23,0.42)
\psline(4.00,0.00)(3.23,0.42)
\psline(3.47,0.85)(3.47,0.85)
\psline(3.23,0.42)(3.47,0.85)
\psline(4.00,0.00)(4.00,0.00)
\psline(3.47,0.85)(4.00,0.00)
\psline(3.23,0.42)(3.23,0.42)
\psline(5.60,10.15)(0.00,0.00)
\psline(0.00,0.00)(11.21,0.00)
\psline(5.60,10.15)(11.21,0.00)
\pspolygon[fillstyle=solid,linewidth=1pt,fillcolor=lightyellow](2.40,2.54)(2.63,2.96)(1.63,2.96)
\psline(2.63,2.96)(1.63,2.96)
\psline(2.40,2.54)(2.40,2.54)
\psline(1.63,2.96)(2.40,2.54)
\psline(2.63,2.96)(2.63,2.96)
\psline(2.40,2.54)(2.63,2.96)
\psline(1.63,2.96)(1.63,2.96)
\pspolygon[fillstyle=solid,linewidth=1pt,fillcolor=lightyellow](1.63,2.96)(1.87,3.38)(2.63,2.96)
\psline(1.87,3.38)(2.63,2.96)
\psline(1.63,2.96)(1.63,2.96)
\psline(2.63,2.96)(1.63,2.96)
\psline(1.87,3.38)(1.87,3.38)
\psline(1.63,2.96)(1.87,3.38)
\psline(2.63,2.96)(2.63,2.96)
\pspolygon[fillstyle=solid,linewidth=1pt,fillcolor=lightyellow](3.17,2.12)(3.40,2.54)(2.40,2.54)
\psline(3.40,2.54)(2.40,2.54)
\psline(3.17,2.12)(3.17,2.12)
\psline(2.40,2.54)(3.17,2.12)
\psline(3.40,2.54)(3.40,2.54)
\psline(3.17,2.12)(3.40,2.54)
\psline(2.40,2.54)(2.40,2.54)
\pspolygon[fillstyle=solid,linewidth=1pt,fillcolor=lightyellow](2.40,2.54)(2.63,2.96)(3.40,2.54)
\psline(2.63,2.96)(3.40,2.54)
\psline(2.40,2.54)(2.40,2.54)
\psline(3.40,2.54)(2.40,2.54)
\psline(2.63,2.96)(2.63,2.96)
\psline(2.40,2.54)(2.63,2.96)
\psline(3.40,2.54)(3.40,2.54)
\pspolygon[fillstyle=solid,linewidth=1pt,fillcolor=lightyellow](3.93,1.69)(4.17,2.12)(3.17,2.12)
\psline(4.17,2.12)(3.17,2.12)
\psline(3.93,1.69)(3.93,1.69)
\psline(3.17,2.12)(3.93,1.69)
\psline(4.17,2.12)(4.17,2.12)
\psline(3.93,1.69)(4.17,2.12)
\psline(3.17,2.12)(3.17,2.12)
\pspolygon[fillstyle=solid,linewidth=1pt,fillcolor=lightyellow](3.17,2.12)(3.40,2.54)(4.17,2.12)
\psline(3.40,2.54)(4.17,2.12)
\psline(3.17,2.12)(3.17,2.12)
\psline(4.17,2.12)(3.17,2.12)
\psline(3.40,2.54)(3.40,2.54)
\psline(3.17,2.12)(3.40,2.54)
\psline(4.17,2.12)(4.17,2.12)
\pspolygon[fillstyle=solid,linewidth=1pt,fillcolor=lightyellow](4.70,1.27)(4.93,1.69)(3.93,1.69)
\psline(4.93,1.69)(3.93,1.69)
\psline(4.70,1.27)(4.70,1.27)
\psline(3.93,1.69)(4.70,1.27)
\psline(4.93,1.69)(4.93,1.69)
\psline(4.70,1.27)(4.93,1.69)
\psline(3.93,1.69)(3.93,1.69)
\pspolygon[fillstyle=solid,linewidth=1pt,fillcolor=lightyellow](3.93,1.69)(4.17,2.12)(4.93,1.69)
\psline(4.17,2.12)(4.93,1.69)
\psline(3.93,1.69)(3.93,1.69)
\psline(4.93,1.69)(3.93,1.69)
\psline(4.17,2.12)(4.17,2.12)
\psline(3.93,1.69)(4.17,2.12)
\psline(4.93,1.69)(4.93,1.69)
\pspolygon[fillstyle=solid,linewidth=1pt,fillcolor=lightyellow](5.47,0.85)(5.70,1.27)(4.70,1.27)
\psline(5.70,1.27)(4.70,1.27)
\psline(5.47,0.85)(5.47,0.85)
\psline(4.70,1.27)(5.47,0.85)
\psline(5.70,1.27)(5.70,1.27)
\psline(5.47,0.85)(5.70,1.27)
\psline(4.70,1.27)(4.70,1.27)
\pspolygon[fillstyle=solid,linewidth=1pt,fillcolor=lightyellow](4.70,1.27)(4.93,1.69)(5.70,1.27)
\psline(4.93,1.69)(5.70,1.27)
\psline(4.70,1.27)(4.70,1.27)
\psline(5.70,1.27)(4.70,1.27)
\psline(4.93,1.69)(4.93,1.69)
\psline(4.70,1.27)(4.93,1.69)
\psline(5.70,1.27)(5.70,1.27)
\pspolygon[fillstyle=solid,linewidth=1pt,fillcolor=lightblue](5.47,0.85)(6.47,0.85)(6.23,0.42)
\psline(6.47,0.85)(6.23,0.42)
\psline(5.47,0.85)(5.47,0.85)
\psline(6.23,0.42)(5.47,0.85)
\psline(6.47,0.85)(6.47,0.85)
\psline(5.47,0.85)(6.47,0.85)
\psline(6.23,0.42)(6.23,0.42)
\put(6.06,0.51){$5$}
\pspolygon[fillstyle=solid,linewidth=1pt,fillcolor=lightgreen](5.47,0.85)(5.95,0.85)(5.95,1.72)
\psline(5.95,0.85)(5.95,1.72)
\psline(5.47,0.85)(5.47,0.85)
\psline(5.95,1.72)(5.47,0.85)
\psline(5.95,0.85)(5.95,0.85)
\psline(5.47,0.85)(5.95,0.85)
\psline(5.95,1.72)(5.95,1.72)
\put(5.69,0.90){$6$}
\psdot(1.63,2.96)\put(1.03,2.95){$Q$}
\psdot(1.87,3.38)\put(1.27,3.37){$P$}
\psdot(6.47,0.85)\put(6.66,0.74){$S_{k+1}$}
\psdot(6.23,0.42)\put(6.43,0.32){$R_{k+1}$}
\psline[linestyle=dashed](5.70,1.27)(8.12,5.65)
\pspolygon[fillstyle=solid,linewidth=1pt,fillcolor=lightyellow](2.63,2.96)(2.87,3.38)(1.87,3.38)
\psline(2.87,3.38)(1.87,3.38)
\psline(2.63,2.96)(2.63,2.96)
\psline(1.87,3.38)(2.63,2.96)
\psline(2.87,3.38)(2.87,3.38)
\psline(2.63,2.96)(2.87,3.38)
\psline(1.87,3.38)(1.87,3.38)
\pspolygon[fillstyle=solid,linewidth=1pt,fillcolor=lightyellow](1.87,3.38)(2.10,3.81)(2.87,3.38)
\psline(2.10,3.81)(2.87,3.38)
\psline(1.87,3.38)(1.87,3.38)
\psline(2.87,3.38)(1.87,3.38)
\psline(2.10,3.81)(2.10,3.81)
\psline(1.87,3.38)(2.10,3.81)
\psline(2.87,3.38)(2.87,3.38)
\pspolygon[fillstyle=solid,linewidth=1pt,fillcolor=lightyellow](3.40,2.54)(3.63,2.96)(2.63,2.96)
\psline(3.63,2.96)(2.63,2.96)
\psline(3.40,2.54)(3.40,2.54)
\psline(2.63,2.96)(3.40,2.54)
\psline(3.63,2.96)(3.63,2.96)
\psline(3.40,2.54)(3.63,2.96)
\psline(2.63,2.96)(2.63,2.96)
\pspolygon[fillstyle=solid,linewidth=1pt,fillcolor=lightyellow](2.63,2.96)(2.87,3.38)(3.63,2.96)
\psline(2.87,3.38)(3.63,2.96)
\psline(2.63,2.96)(2.63,2.96)
\psline(3.63,2.96)(2.63,2.96)
\psline(2.87,3.38)(2.87,3.38)
\psline(2.63,2.96)(2.87,3.38)
\psline(3.63,2.96)(3.63,2.96)
\pspolygon[fillstyle=solid,linewidth=1pt,fillcolor=lightyellow](4.17,2.12)(4.40,2.54)(3.40,2.54)
\psline(4.40,2.54)(3.40,2.54)
\psline(4.17,2.12)(4.17,2.12)
\psline(3.40,2.54)(4.17,2.12)
\psline(4.40,2.54)(4.40,2.54)
\psline(4.17,2.12)(4.40,2.54)
\psline(3.40,2.54)(3.40,2.54)
\pspolygon[fillstyle=solid,linewidth=1pt,fillcolor=lightyellow](3.40,2.54)(3.63,2.96)(4.40,2.54)
\psline(3.63,2.96)(4.40,2.54)
\psline(3.40,2.54)(3.40,2.54)
\psline(4.40,2.54)(3.40,2.54)
\psline(3.63,2.96)(3.63,2.96)
\psline(3.40,2.54)(3.63,2.96)
\psline(4.40,2.54)(4.40,2.54)
\pspolygon[fillstyle=solid,linewidth=1pt,fillcolor=lightyellow](4.93,1.69)(5.17,2.12)(4.17,2.12)
\psline(5.17,2.12)(4.17,2.12)
\psline(4.93,1.69)(4.93,1.69)
\psline(4.17,2.12)(4.93,1.69)
\psline(5.17,2.12)(5.17,2.12)
\psline(4.93,1.69)(5.17,2.12)
\psline(4.17,2.12)(4.17,2.12)
\pspolygon[fillstyle=solid,linewidth=1pt,fillcolor=lightyellow](4.17,2.12)(4.40,2.54)(5.17,2.12)
\psline(4.40,2.54)(5.17,2.12)
\psline(4.17,2.12)(4.17,2.12)
\psline(5.17,2.12)(4.17,2.12)
\psline(4.40,2.54)(4.40,2.54)
\psline(4.17,2.12)(4.40,2.54)
\psline(5.17,2.12)(5.17,2.12)
\pspolygon[fillstyle=solid,linewidth=1pt,fillcolor=lightyellow](5.70,1.27)(5.93,1.69)(4.93,1.69)
\psline(5.93,1.69)(4.93,1.69)
\psline(5.70,1.27)(5.70,1.27)
\psline(4.93,1.69)(5.70,1.27)
\psline(5.93,1.69)(5.93,1.69)
\psline(5.70,1.27)(5.93,1.69)
\psline(4.93,1.69)(4.93,1.69)
\pspolygon[fillstyle=solid,linewidth=1pt,fillcolor=lightyellow](4.93,1.69)(5.17,2.12)(5.93,1.69)
\psline(5.17,2.12)(5.93,1.69)
\psline(4.93,1.69)(4.93,1.69)
\psline(5.93,1.69)(4.93,1.69)
\psline(5.17,2.12)(5.17,2.12)
\psline(4.93,1.69)(5.17,2.12)
\psline(5.93,1.69)(5.93,1.69)
\pspolygon[fillstyle=solid,linewidth=1pt,fillcolor=lightyellow](2.87,3.38)(3.10,3.81)(2.10,3.81)
\psline(3.10,3.81)(2.10,3.81)
\psline(2.87,3.38)(2.87,3.38)
\psline(2.10,3.81)(2.87,3.38)
\psline(3.10,3.81)(3.10,3.81)
\psline(2.87,3.38)(3.10,3.81)
\psline(2.10,3.81)(2.10,3.81)
\pspolygon[fillstyle=solid,linewidth=1pt,fillcolor=lightyellow](2.10,3.81)(2.33,4.23)(3.10,3.81)
\psline(2.33,4.23)(3.10,3.81)
\psline(2.10,3.81)(2.10,3.81)
\psline(3.10,3.81)(2.10,3.81)
\psline(2.33,4.23)(2.33,4.23)
\psline(2.10,3.81)(2.33,4.23)
\psline(3.10,3.81)(3.10,3.81)
\pspolygon[fillstyle=solid,linewidth=1pt,fillcolor=lightyellow](3.63,2.96)(3.87,3.38)(2.87,3.38)
\psline(3.87,3.38)(2.87,3.38)
\psline(3.63,2.96)(3.63,2.96)
\psline(2.87,3.38)(3.63,2.96)
\psline(3.87,3.38)(3.87,3.38)
\psline(3.63,2.96)(3.87,3.38)
\psline(2.87,3.38)(2.87,3.38)
\pspolygon[fillstyle=solid,linewidth=1pt,fillcolor=lightyellow](2.87,3.38)(3.10,3.81)(3.87,3.38)
\psline(3.10,3.81)(3.87,3.38)
\psline(2.87,3.38)(2.87,3.38)
\psline(3.87,3.38)(2.87,3.38)
\psline(3.10,3.81)(3.10,3.81)
\psline(2.87,3.38)(3.10,3.81)
\psline(3.87,3.38)(3.87,3.38)
\pspolygon[fillstyle=solid,linewidth=1pt,fillcolor=lightyellow](4.40,2.54)(4.63,2.96)(3.63,2.96)
\psline(4.63,2.96)(3.63,2.96)
\psline(4.40,2.54)(4.40,2.54)
\psline(3.63,2.96)(4.40,2.54)
\psline(4.63,2.96)(4.63,2.96)
\psline(4.40,2.54)(4.63,2.96)
\psline(3.63,2.96)(3.63,2.96)
\pspolygon[fillstyle=solid,linewidth=1pt,fillcolor=lightyellow](3.63,2.96)(3.87,3.38)(4.63,2.96)
\psline(3.87,3.38)(4.63,2.96)
\psline(3.63,2.96)(3.63,2.96)
\psline(4.63,2.96)(3.63,2.96)
\psline(3.87,3.38)(3.87,3.38)
\psline(3.63,2.96)(3.87,3.38)
\psline(4.63,2.96)(4.63,2.96)
\pspolygon[fillstyle=solid,linewidth=1pt,fillcolor=lightyellow](5.17,2.12)(5.40,2.54)(4.40,2.54)
\psline(5.40,2.54)(4.40,2.54)
\psline(5.17,2.12)(5.17,2.12)
\psline(4.40,2.54)(5.17,2.12)
\psline(5.40,2.54)(5.40,2.54)
\psline(5.17,2.12)(5.40,2.54)
\psline(4.40,2.54)(4.40,2.54)
\pspolygon[fillstyle=solid,linewidth=1pt,fillcolor=lightyellow](4.40,2.54)(4.63,2.96)(5.40,2.54)
\psline(4.63,2.96)(5.40,2.54)
\psline(4.40,2.54)(4.40,2.54)
\psline(5.40,2.54)(4.40,2.54)
\psline(4.63,2.96)(4.63,2.96)
\psline(4.40,2.54)(4.63,2.96)
\psline(5.40,2.54)(5.40,2.54)
\pspolygon[fillstyle=solid,linewidth=1pt,fillcolor=lightyellow](5.93,1.69)(6.17,2.12)(5.17,2.12)
\psline(6.17,2.12)(5.17,2.12)
\psline(5.93,1.69)(5.93,1.69)
\psline(5.17,2.12)(5.93,1.69)
\psline(6.17,2.12)(6.17,2.12)
\psline(5.93,1.69)(6.17,2.12)
\psline(5.17,2.12)(5.17,2.12)
\pspolygon[fillstyle=solid,linewidth=1pt,fillcolor=lightyellow](5.17,2.12)(5.40,2.54)(6.17,2.12)
\psline(5.40,2.54)(6.17,2.12)
\psline(5.17,2.12)(5.17,2.12)
\psline(6.17,2.12)(5.17,2.12)
\psline(5.40,2.54)(5.40,2.54)
\psline(5.17,2.12)(5.40,2.54)
\psline(6.17,2.12)(6.17,2.12)
\psdot(2.33,4.23)\put(1.73,4.22){$Q^\prime$}
\psdot(6.17,2.12)\put(6.56,1.92){$R^\prime$}
\pspolygon[fillstyle=solid,linewidth=1pt,fillcolor=lightyellow](3.10,3.81)(3.33,4.23)(2.33,4.23)
\psline(3.33,4.23)(2.33,4.23)
\psline(3.10,3.81)(3.10,3.81)
\psline(2.33,4.23)(3.10,3.81)
\psline(3.33,4.23)(3.33,4.23)
\psline(3.10,3.81)(3.33,4.23)
\psline(2.33,4.23)(2.33,4.23)
\pspolygon[fillstyle=solid,linewidth=1pt,fillcolor=lightyellow](2.33,4.23)(2.57,4.65)(3.33,4.23)
\psline(2.57,4.65)(3.33,4.23)
\psline(2.33,4.23)(2.33,4.23)
\psline(3.33,4.23)(2.33,4.23)
\psline(2.57,4.65)(2.57,4.65)
\psline(2.33,4.23)(2.57,4.65)
\psline(3.33,4.23)(3.33,4.23)
\pspolygon[fillstyle=solid,linewidth=1pt,fillcolor=lightyellow](3.87,3.38)(4.10,3.81)(3.10,3.81)
\psline(4.10,3.81)(3.10,3.81)
\psline(3.87,3.38)(3.87,3.38)
\psline(3.10,3.81)(3.87,3.38)
\psline(4.10,3.81)(4.10,3.81)
\psline(3.87,3.38)(4.10,3.81)
\psline(3.10,3.81)(3.10,3.81)
\pspolygon[fillstyle=solid,linewidth=1pt,fillcolor=lightyellow](3.10,3.81)(3.33,4.23)(4.10,3.81)
\psline(3.33,4.23)(4.10,3.81)
\psline(3.10,3.81)(3.10,3.81)
\psline(4.10,3.81)(3.10,3.81)
\psline(3.33,4.23)(3.33,4.23)
\psline(3.10,3.81)(3.33,4.23)
\psline(4.10,3.81)(4.10,3.81)
\pspolygon[fillstyle=solid,linewidth=1pt,fillcolor=lightyellow](4.63,2.96)(4.87,3.38)(3.87,3.38)
\psline(4.87,3.38)(3.87,3.38)
\psline(4.63,2.96)(4.63,2.96)
\psline(3.87,3.38)(4.63,2.96)
\psline(4.87,3.38)(4.87,3.38)
\psline(4.63,2.96)(4.87,3.38)
\psline(3.87,3.38)(3.87,3.38)
\pspolygon[fillstyle=solid,linewidth=1pt,fillcolor=lightyellow](3.87,3.38)(4.10,3.81)(4.87,3.38)
\psline(4.10,3.81)(4.87,3.38)
\psline(3.87,3.38)(3.87,3.38)
\psline(4.87,3.38)(3.87,3.38)
\psline(4.10,3.81)(4.10,3.81)
\psline(3.87,3.38)(4.10,3.81)
\psline(4.87,3.38)(4.87,3.38)
\pspolygon[fillstyle=solid,linewidth=1pt,fillcolor=lightblue](4.63,2.96)(5.63,2.96)(5.40,2.54)
\psline(5.63,2.96)(5.40,2.54)
\psline(4.63,2.96)(4.63,2.96)
\psline(5.40,2.54)(4.63,2.96)
\psline(5.63,2.96)(5.63,2.96)
\psline(4.63,2.96)(5.63,2.96)
\psline(5.40,2.54)(5.40,2.54)
\put(5.22,2.63){$7$}
\pspolygon[fillstyle=solid,linewidth=1pt,fillcolor=lightgreen](4.63,2.96)(5.12,2.96)(5.12,3.84)
\psline(5.12,2.96)(5.12,3.84)
\psline(4.63,2.96)(4.63,2.96)
\psline(5.12,3.84)(4.63,2.96)
\psline(5.12,2.96)(5.12,2.96)
\psline(4.63,2.96)(5.12,2.96)
\psline(5.12,3.84)(5.12,3.84)
\put(4.86,3.01){$8$}
\psdot(2.33,4.23)\put(1.73,4.22){$Q$}
\psdot(2.57,4.65)\put(1.97,4.64){$P$}
\psdot(4.87,3.38)\put(4.42,3.70){$S_j$}
\psdot(5.63,2.96)\put(5.83,2.86){$S_{j+1}$}
\psdot(5.40,2.54)\put(5.59,2.44){$R_{j+1}$}
\psdot(5.60,10.15)\put(5.04,9.86){$A$}
\psdot(0.00,0.00)\put(0.00,-0.44){$B$}
\psdot(11.21,0.00)\put(11.21,-0.44){$C$}
\psdot(7.00,0.00)\put(7.00,-0.44){$U$}
\psline(4.87,3.38)(7.11,7.46)
\endpspicture
}

\def\FigureFourThousand{
\hskip 1in
\pspicture(2.0,2.8)
\psset{unit=0.8cm}
\newrgbcolor{lightblue}{0.8 0.8 1}
\newrgbcolor{pink}{1 0.8 0.8}
\newrgbcolor{lightgreen}{0.8 1 0.8}
\newrgbcolor{lightyellow}{1 1 0.8}
\pspolygon[fillstyle=solid,linewidth=1pt,fillcolor=lightblue](1.63,2.96)(0.00,0.00)(7.00,0.00)
\psline(0.00,0.00)(7.00,0.00)
\psline(0.23,0.42)(6.23,0.42)
\psline(0.47,0.85)(5.47,0.85)
\psline(0.70,1.27)(4.70,1.27)
\psline(0.93,1.69)(3.93,1.69)
\psline(1.17,2.12)(3.17,2.12)
\psline(1.40,2.54)(2.40,2.54)
\psline(1.63,2.96)(1.63,2.96)
\psline(7.00,0.00)(1.63,2.96)
\psline(6.00,0.00)(1.40,2.54)
\psline(5.00,0.00)(1.17,2.12)
\psline(4.00,0.00)(0.93,1.69)
\psline(3.00,0.00)(0.70,1.27)
\psline(2.00,0.00)(0.47,0.85)
\psline(1.00,0.00)(0.23,0.42)
\psline(0.00,0.00)(0.00,0.00)
\psline(1.63,2.96)(0.00,0.00)
\psline(2.40,2.54)(1.00,0.00)
\psline(3.17,2.12)(2.00,0.00)
\psline(3.93,1.69)(3.00,0.00)
\psline(4.70,1.27)(4.00,0.00)
\psline(5.47,0.85)(5.00,0.00)
\psline(6.23,0.42)(6.00,0.00)
\psline(7.00,0.00)(7.00,0.00)
\pspolygon[fillstyle=solid,linewidth=1pt,fillcolor=pink](0.93,1.69)(1.70,1.27)(1.47,0.85)
\psline(1.70,1.27)(1.47,0.85)
\psline(0.93,1.69)(0.93,1.69)
\psline(1.47,0.85)(0.93,1.69)
\psline(1.70,1.27)(1.70,1.27)
\psline(0.93,1.69)(1.70,1.27)
\psline(1.47,0.85)(1.47,0.85)
\pspolygon[fillstyle=solid,linewidth=1pt,fillcolor=pink](0.93,1.69)(1.47,0.85)(0.70,1.27)
\psline(1.47,0.85)(0.70,1.27)
\psline(0.93,1.69)(0.93,1.69)
\psline(0.70,1.27)(0.93,1.69)
\psline(1.47,0.85)(1.47,0.85)
\psline(0.93,1.69)(1.47,0.85)
\psline(0.70,1.27)(0.70,1.27)
\pspolygon[fillstyle=solid,linewidth=1pt,fillcolor=pink](1.93,1.69)(2.70,1.27)(2.47,0.85)
\psline(2.70,1.27)(2.47,0.85)
\psline(1.93,1.69)(1.93,1.69)
\psline(2.47,0.85)(1.93,1.69)
\psline(2.70,1.27)(2.70,1.27)
\psline(1.93,1.69)(2.70,1.27)
\psline(2.47,0.85)(2.47,0.85)
\pspolygon[fillstyle=solid,linewidth=1pt,fillcolor=pink](1.93,1.69)(2.47,0.85)(1.70,1.27)
\psline(2.47,0.85)(1.70,1.27)
\psline(1.93,1.69)(1.93,1.69)
\psline(1.70,1.27)(1.93,1.69)
\psline(2.47,0.85)(2.47,0.85)
\psline(1.93,1.69)(2.47,0.85)
\psline(1.70,1.27)(1.70,1.27)
\pspolygon[fillstyle=solid,linewidth=1pt,fillcolor=pink](2.93,1.69)(3.70,1.27)(3.47,0.85)
\psline(3.70,1.27)(3.47,0.85)
\psline(2.93,1.69)(2.93,1.69)
\psline(3.47,0.85)(2.93,1.69)
\psline(3.70,1.27)(3.70,1.27)
\psline(2.93,1.69)(3.70,1.27)
\psline(3.47,0.85)(3.47,0.85)
\pspolygon[fillstyle=solid,linewidth=1pt,fillcolor=pink](2.93,1.69)(3.47,0.85)(2.70,1.27)
\psline(3.47,0.85)(2.70,1.27)
\psline(2.93,1.69)(2.93,1.69)
\psline(2.70,1.27)(2.93,1.69)
\psline(3.47,0.85)(3.47,0.85)
\psline(2.93,1.69)(3.47,0.85)
\psline(2.70,1.27)(2.70,1.27)
\pspolygon[fillstyle=solid,linewidth=1pt,fillcolor=pink](3.47,0.85)(4.23,0.42)(4.00,0.00)
\psline(4.23,0.42)(4.00,0.00)
\psline(3.47,0.85)(3.47,0.85)
\psline(4.00,0.00)(3.47,0.85)
\psline(4.23,0.42)(4.23,0.42)
\psline(3.47,0.85)(4.23,0.42)
\psline(4.00,0.00)(4.00,0.00)
\pspolygon[fillstyle=solid,linewidth=1pt,fillcolor=pink](3.47,0.85)(4.00,0.00)(3.23,0.42)
\psline(4.00,0.00)(3.23,0.42)
\psline(3.47,0.85)(3.47,0.85)
\psline(3.23,0.42)(3.47,0.85)
\psline(4.00,0.00)(4.00,0.00)
\psline(3.47,0.85)(4.00,0.00)
\psline(3.23,0.42)(3.23,0.42)
\psline(5.60,10.15)(0.00,0.00)
\psline(0.00,0.00)(11.21,0.00)
\psline(5.60,10.15)(11.21,0.00)
\pspolygon[fillstyle=solid,linewidth=1pt,fillcolor=lightyellow](2.40,2.54)(2.63,2.96)(1.63,2.96)
\psline(2.63,2.96)(1.63,2.96)
\psline(2.40,2.54)(2.40,2.54)
\psline(1.63,2.96)(2.40,2.54)
\psline(2.63,2.96)(2.63,2.96)
\psline(2.40,2.54)(2.63,2.96)
\psline(1.63,2.96)(1.63,2.96)
\pspolygon[fillstyle=solid,linewidth=1pt,fillcolor=lightyellow](1.63,2.96)(1.87,3.38)(2.63,2.96)
\psline(1.87,3.38)(2.63,2.96)
\psline(1.63,2.96)(1.63,2.96)
\psline(2.63,2.96)(1.63,2.96)
\psline(1.87,3.38)(1.87,3.38)
\psline(1.63,2.96)(1.87,3.38)
\psline(2.63,2.96)(2.63,2.96)
\pspolygon[fillstyle=solid,linewidth=1pt,fillcolor=lightyellow](3.17,2.12)(3.40,2.54)(2.40,2.54)
\psline(3.40,2.54)(2.40,2.54)
\psline(3.17,2.12)(3.17,2.12)
\psline(2.40,2.54)(3.17,2.12)
\psline(3.40,2.54)(3.40,2.54)
\psline(3.17,2.12)(3.40,2.54)
\psline(2.40,2.54)(2.40,2.54)
\pspolygon[fillstyle=solid,linewidth=1pt,fillcolor=lightyellow](2.40,2.54)(2.63,2.96)(3.40,2.54)
\psline(2.63,2.96)(3.40,2.54)
\psline(2.40,2.54)(2.40,2.54)
\psline(3.40,2.54)(2.40,2.54)
\psline(2.63,2.96)(2.63,2.96)
\psline(2.40,2.54)(2.63,2.96)
\psline(3.40,2.54)(3.40,2.54)
\pspolygon[fillstyle=solid,linewidth=1pt,fillcolor=lightyellow](3.93,1.69)(4.17,2.12)(3.17,2.12)
\psline(4.17,2.12)(3.17,2.12)
\psline(3.93,1.69)(3.93,1.69)
\psline(3.17,2.12)(3.93,1.69)
\psline(4.17,2.12)(4.17,2.12)
\psline(3.93,1.69)(4.17,2.12)
\psline(3.17,2.12)(3.17,2.12)
\pspolygon[fillstyle=solid,linewidth=1pt,fillcolor=lightyellow](3.17,2.12)(3.40,2.54)(4.17,2.12)
\psline(3.40,2.54)(4.17,2.12)
\psline(3.17,2.12)(3.17,2.12)
\psline(4.17,2.12)(3.17,2.12)
\psline(3.40,2.54)(3.40,2.54)
\psline(3.17,2.12)(3.40,2.54)
\psline(4.17,2.12)(4.17,2.12)
\pspolygon[fillstyle=solid,linewidth=1pt,fillcolor=lightyellow](4.70,1.27)(4.93,1.69)(3.93,1.69)
\psline(4.93,1.69)(3.93,1.69)
\psline(4.70,1.27)(4.70,1.27)
\psline(3.93,1.69)(4.70,1.27)
\psline(4.93,1.69)(4.93,1.69)
\psline(4.70,1.27)(4.93,1.69)
\psline(3.93,1.69)(3.93,1.69)
\pspolygon[fillstyle=solid,linewidth=1pt,fillcolor=lightyellow](3.93,1.69)(4.17,2.12)(4.93,1.69)
\psline(4.17,2.12)(4.93,1.69)
\psline(3.93,1.69)(3.93,1.69)
\psline(4.93,1.69)(3.93,1.69)
\psline(4.17,2.12)(4.17,2.12)
\psline(3.93,1.69)(4.17,2.12)
\psline(4.93,1.69)(4.93,1.69)
\pspolygon[fillstyle=solid,linewidth=1pt,fillcolor=lightyellow](5.47,0.85)(5.70,1.27)(4.70,1.27)
\psline(5.70,1.27)(4.70,1.27)
\psline(5.47,0.85)(5.47,0.85)
\psline(4.70,1.27)(5.47,0.85)
\psline(5.70,1.27)(5.70,1.27)
\psline(5.47,0.85)(5.70,1.27)
\psline(4.70,1.27)(4.70,1.27)
\pspolygon[fillstyle=solid,linewidth=1pt,fillcolor=lightyellow](4.70,1.27)(4.93,1.69)(5.70,1.27)
\psline(4.93,1.69)(5.70,1.27)
\psline(4.70,1.27)(4.70,1.27)
\psline(5.70,1.27)(4.70,1.27)
\psline(4.93,1.69)(4.93,1.69)
\psline(4.70,1.27)(4.93,1.69)
\psline(5.70,1.27)(5.70,1.27)
\pspolygon[fillstyle=solid,linewidth=1pt,fillcolor=lightblue](5.47,0.85)(6.47,0.85)(6.23,0.42)
\psline(6.47,0.85)(6.23,0.42)
\psline(5.47,0.85)(5.47,0.85)
\psline(6.23,0.42)(5.47,0.85)
\psline(6.47,0.85)(6.47,0.85)
\psline(5.47,0.85)(6.47,0.85)
\psline(6.23,0.42)(6.23,0.42)
\put(6.06,0.51){$5$}
\pspolygon[fillstyle=solid,linewidth=1pt,fillcolor=lightgreen](5.47,0.85)(5.95,0.85)(5.95,1.72)
\psline(5.95,0.85)(5.95,1.72)
\psline(5.47,0.85)(5.47,0.85)
\psline(5.95,1.72)(5.47,0.85)
\psline(5.95,0.85)(5.95,0.85)
\psline(5.47,0.85)(5.95,0.85)
\psline(5.95,1.72)(5.95,1.72)
\put(5.69,0.90){$6$}
\psdot(1.63,2.96)\put(1.03,2.95){$Q$}
\psdot(1.87,3.38)\put(1.27,3.37){$P$}
\psdot(6.47,0.85)\put(6.66,0.74){$S_{k+1}$}
\psdot(6.23,0.42)\put(6.43,0.32){$R_{k+1}$}
\psline[linestyle=dashed](5.70,1.27)(8.12,5.65)
\pspolygon[fillstyle=solid,linewidth=1pt,fillcolor=lightyellow](2.63,2.96)(2.87,3.38)(1.87,3.38)
\psline(2.87,3.38)(1.87,3.38)
\psline(2.63,2.96)(2.63,2.96)
\psline(1.87,3.38)(2.63,2.96)
\psline(2.87,3.38)(2.87,3.38)
\psline(2.63,2.96)(2.87,3.38)
\psline(1.87,3.38)(1.87,3.38)
\pspolygon[fillstyle=solid,linewidth=1pt,fillcolor=lightyellow](1.87,3.38)(2.10,3.81)(2.87,3.38)
\psline(2.10,3.81)(2.87,3.38)
\psline(1.87,3.38)(1.87,3.38)
\psline(2.87,3.38)(1.87,3.38)
\psline(2.10,3.81)(2.10,3.81)
\psline(1.87,3.38)(2.10,3.81)
\psline(2.87,3.38)(2.87,3.38)
\pspolygon[fillstyle=solid,linewidth=1pt,fillcolor=lightyellow](3.40,2.54)(3.63,2.96)(2.63,2.96)
\psline(3.63,2.96)(2.63,2.96)
\psline(3.40,2.54)(3.40,2.54)
\psline(2.63,2.96)(3.40,2.54)
\psline(3.63,2.96)(3.63,2.96)
\psline(3.40,2.54)(3.63,2.96)
\psline(2.63,2.96)(2.63,2.96)
\pspolygon[fillstyle=solid,linewidth=1pt,fillcolor=lightyellow](2.63,2.96)(2.87,3.38)(3.63,2.96)
\psline(2.87,3.38)(3.63,2.96)
\psline(2.63,2.96)(2.63,2.96)
\psline(3.63,2.96)(2.63,2.96)
\psline(2.87,3.38)(2.87,3.38)
\psline(2.63,2.96)(2.87,3.38)
\psline(3.63,2.96)(3.63,2.96)
\pspolygon[fillstyle=solid,linewidth=1pt,fillcolor=lightyellow](4.17,2.12)(4.40,2.54)(3.40,2.54)
\psline(4.40,2.54)(3.40,2.54)
\psline(4.17,2.12)(4.17,2.12)
\psline(3.40,2.54)(4.17,2.12)
\psline(4.40,2.54)(4.40,2.54)
\psline(4.17,2.12)(4.40,2.54)
\psline(3.40,2.54)(3.40,2.54)
\pspolygon[fillstyle=solid,linewidth=1pt,fillcolor=lightyellow](3.40,2.54)(3.63,2.96)(4.40,2.54)
\psline(3.63,2.96)(4.40,2.54)
\psline(3.40,2.54)(3.40,2.54)
\psline(4.40,2.54)(3.40,2.54)
\psline(3.63,2.96)(3.63,2.96)
\psline(3.40,2.54)(3.63,2.96)
\psline(4.40,2.54)(4.40,2.54)
\pspolygon[fillstyle=solid,linewidth=1pt,fillcolor=lightyellow](4.93,1.69)(5.17,2.12)(4.17,2.12)
\psline(5.17,2.12)(4.17,2.12)
\psline(4.93,1.69)(4.93,1.69)
\psline(4.17,2.12)(4.93,1.69)
\psline(5.17,2.12)(5.17,2.12)
\psline(4.93,1.69)(5.17,2.12)
\psline(4.17,2.12)(4.17,2.12)
\pspolygon[fillstyle=solid,linewidth=1pt,fillcolor=lightyellow](4.17,2.12)(4.40,2.54)(5.17,2.12)
\psline(4.40,2.54)(5.17,2.12)
\psline(4.17,2.12)(4.17,2.12)
\psline(5.17,2.12)(4.17,2.12)
\psline(4.40,2.54)(4.40,2.54)
\psline(4.17,2.12)(4.40,2.54)
\psline(5.17,2.12)(5.17,2.12)
\pspolygon[fillstyle=solid,linewidth=1pt,fillcolor=lightyellow](5.70,1.27)(5.93,1.69)(4.93,1.69)
\psline(5.93,1.69)(4.93,1.69)
\psline(5.70,1.27)(5.70,1.27)
\psline(4.93,1.69)(5.70,1.27)
\psline(5.93,1.69)(5.93,1.69)
\psline(5.70,1.27)(5.93,1.69)
\psline(4.93,1.69)(4.93,1.69)
\pspolygon[fillstyle=solid,linewidth=1pt,fillcolor=lightyellow](4.93,1.69)(5.17,2.12)(5.93,1.69)
\psline(5.17,2.12)(5.93,1.69)
\psline(4.93,1.69)(4.93,1.69)
\psline(5.93,1.69)(4.93,1.69)
\psline(5.17,2.12)(5.17,2.12)
\psline(4.93,1.69)(5.17,2.12)
\psline(5.93,1.69)(5.93,1.69)
\pspolygon[fillstyle=solid,linewidth=1pt,fillcolor=lightyellow](2.87,3.38)(3.10,3.81)(2.10,3.81)
\psline(3.10,3.81)(2.10,3.81)
\psline(2.87,3.38)(2.87,3.38)
\psline(2.10,3.81)(2.87,3.38)
\psline(3.10,3.81)(3.10,3.81)
\psline(2.87,3.38)(3.10,3.81)
\psline(2.10,3.81)(2.10,3.81)
\pspolygon[fillstyle=solid,linewidth=1pt,fillcolor=lightyellow](2.10,3.81)(2.33,4.23)(3.10,3.81)
\psline(2.33,4.23)(3.10,3.81)
\psline(2.10,3.81)(2.10,3.81)
\psline(3.10,3.81)(2.10,3.81)
\psline(2.33,4.23)(2.33,4.23)
\psline(2.10,3.81)(2.33,4.23)
\psline(3.10,3.81)(3.10,3.81)
\pspolygon[fillstyle=solid,linewidth=1pt,fillcolor=lightyellow](3.63,2.96)(3.87,3.38)(2.87,3.38)
\psline(3.87,3.38)(2.87,3.38)
\psline(3.63,2.96)(3.63,2.96)
\psline(2.87,3.38)(3.63,2.96)
\psline(3.87,3.38)(3.87,3.38)
\psline(3.63,2.96)(3.87,3.38)
\psline(2.87,3.38)(2.87,3.38)
\pspolygon[fillstyle=solid,linewidth=1pt,fillcolor=lightyellow](2.87,3.38)(3.10,3.81)(3.87,3.38)
\psline(3.10,3.81)(3.87,3.38)
\psline(2.87,3.38)(2.87,3.38)
\psline(3.87,3.38)(2.87,3.38)
\psline(3.10,3.81)(3.10,3.81)
\psline(2.87,3.38)(3.10,3.81)
\psline(3.87,3.38)(3.87,3.38)
\pspolygon[fillstyle=solid,linewidth=1pt,fillcolor=lightyellow](4.40,2.54)(4.63,2.96)(3.63,2.96)
\psline(4.63,2.96)(3.63,2.96)
\psline(4.40,2.54)(4.40,2.54)
\psline(3.63,2.96)(4.40,2.54)
\psline(4.63,2.96)(4.63,2.96)
\psline(4.40,2.54)(4.63,2.96)
\psline(3.63,2.96)(3.63,2.96)
\pspolygon[fillstyle=solid,linewidth=1pt,fillcolor=lightyellow](3.63,2.96)(3.87,3.38)(4.63,2.96)
\psline(3.87,3.38)(4.63,2.96)
\psline(3.63,2.96)(3.63,2.96)
\psline(4.63,2.96)(3.63,2.96)
\psline(3.87,3.38)(3.87,3.38)
\psline(3.63,2.96)(3.87,3.38)
\psline(4.63,2.96)(4.63,2.96)
\pspolygon[fillstyle=solid,linewidth=1pt,fillcolor=lightyellow](5.17,2.12)(5.40,2.54)(4.40,2.54)
\psline(5.40,2.54)(4.40,2.54)
\psline(5.17,2.12)(5.17,2.12)
\psline(4.40,2.54)(5.17,2.12)
\psline(5.40,2.54)(5.40,2.54)
\psline(5.17,2.12)(5.40,2.54)
\psline(4.40,2.54)(4.40,2.54)
\pspolygon[fillstyle=solid,linewidth=1pt,fillcolor=lightyellow](4.40,2.54)(4.63,2.96)(5.40,2.54)
\psline(4.63,2.96)(5.40,2.54)
\psline(4.40,2.54)(4.40,2.54)
\psline(5.40,2.54)(4.40,2.54)
\psline(4.63,2.96)(4.63,2.96)
\psline(4.40,2.54)(4.63,2.96)
\psline(5.40,2.54)(5.40,2.54)
\pspolygon[fillstyle=solid,linewidth=1pt,fillcolor=lightyellow](5.93,1.69)(6.17,2.12)(5.17,2.12)
\psline(6.17,2.12)(5.17,2.12)
\psline(5.93,1.69)(5.93,1.69)
\psline(5.17,2.12)(5.93,1.69)
\psline(6.17,2.12)(6.17,2.12)
\psline(5.93,1.69)(6.17,2.12)
\psline(5.17,2.12)(5.17,2.12)
\pspolygon[fillstyle=solid,linewidth=1pt,fillcolor=lightyellow](5.17,2.12)(5.40,2.54)(6.17,2.12)
\psline(5.40,2.54)(6.17,2.12)
\psline(5.17,2.12)(5.17,2.12)
\psline(6.17,2.12)(5.17,2.12)
\psline(5.40,2.54)(5.40,2.54)
\psline(5.17,2.12)(5.40,2.54)
\psline(6.17,2.12)(6.17,2.12)
\psdot(2.33,4.23)\put(1.73,4.22){$Q^\prime$}
\psdot(6.17,2.12)\put(6.56,1.92){$R^\prime$}
\pspolygon[fillstyle=solid,linewidth=1pt,fillcolor=lightyellow](3.10,3.81)(3.33,4.23)(2.33,4.23)
\psline(3.33,4.23)(2.33,4.23)
\psline(3.10,3.81)(3.10,3.81)
\psline(2.33,4.23)(3.10,3.81)
\psline(3.33,4.23)(3.33,4.23)
\psline(3.10,3.81)(3.33,4.23)
\psline(2.33,4.23)(2.33,4.23)
\pspolygon[fillstyle=solid,linewidth=1pt,fillcolor=lightyellow](2.33,4.23)(2.57,4.65)(3.33,4.23)
\psline(2.57,4.65)(3.33,4.23)
\psline(2.33,4.23)(2.33,4.23)
\psline(3.33,4.23)(2.33,4.23)
\psline(2.57,4.65)(2.57,4.65)
\psline(2.33,4.23)(2.57,4.65)
\psline(3.33,4.23)(3.33,4.23)
\pspolygon[fillstyle=solid,linewidth=1pt,fillcolor=lightyellow](3.87,3.38)(4.10,3.81)(3.10,3.81)
\psline(4.10,3.81)(3.10,3.81)
\psline(3.87,3.38)(3.87,3.38)
\psline(3.10,3.81)(3.87,3.38)
\psline(4.10,3.81)(4.10,3.81)
\psline(3.87,3.38)(4.10,3.81)
\psline(3.10,3.81)(3.10,3.81)
\pspolygon[fillstyle=solid,linewidth=1pt,fillcolor=lightyellow](3.10,3.81)(3.33,4.23)(4.10,3.81)
\psline(3.33,4.23)(4.10,3.81)
\psline(3.10,3.81)(3.10,3.81)
\psline(4.10,3.81)(3.10,3.81)
\psline(3.33,4.23)(3.33,4.23)
\psline(3.10,3.81)(3.33,4.23)
\psline(4.10,3.81)(4.10,3.81)
\pspolygon[fillstyle=solid,linewidth=1pt,fillcolor=lightyellow](4.63,2.96)(4.87,3.38)(3.87,3.38)
\psline(4.87,3.38)(3.87,3.38)
\psline(4.63,2.96)(4.63,2.96)
\psline(3.87,3.38)(4.63,2.96)
\psline(4.87,3.38)(4.87,3.38)
\psline(4.63,2.96)(4.87,3.38)
\psline(3.87,3.38)(3.87,3.38)
\pspolygon[fillstyle=solid,linewidth=1pt,fillcolor=lightyellow](3.87,3.38)(4.10,3.81)(4.87,3.38)
\psline(4.10,3.81)(4.87,3.38)
\psline(3.87,3.38)(3.87,3.38)
\psline(4.87,3.38)(3.87,3.38)
\psline(4.10,3.81)(4.10,3.81)
\psline(3.87,3.38)(4.10,3.81)
\psline(4.87,3.38)(4.87,3.38)
\pspolygon[fillstyle=solid,linewidth=1pt,fillcolor=lightblue](4.63,2.96)(5.63,2.96)(5.40,2.54)
\psline(5.63,2.96)(5.40,2.54)
\psline(4.63,2.96)(4.63,2.96)
\psline(5.40,2.54)(4.63,2.96)
\psline(5.63,2.96)(5.63,2.96)
\psline(4.63,2.96)(5.63,2.96)
\psline(5.40,2.54)(5.40,2.54)
\put(5.22,2.63){$7$}
\pspolygon[fillstyle=solid,linewidth=1pt,fillcolor=lightgreen](4.63,2.96)(5.12,2.96)(5.12,3.84)
\psline(5.12,2.96)(5.12,3.84)
\psline(4.63,2.96)(4.63,2.96)
\psline(5.12,3.84)(4.63,2.96)
\psline(5.12,2.96)(5.12,2.96)
\psline(4.63,2.96)(5.12,2.96)
\psline(5.12,3.84)(5.12,3.84)
\put(4.86,3.01){$8$}
\psdot(2.33,4.23)\put(1.73,4.22){$Q$}
\psdot(2.57,4.65)\put(1.97,4.64){$P$}
\psdot(4.87,3.38)\put(4.42,3.70){$S_j$}
\psdot(5.63,2.96)\put(5.83,2.86){$S_{j+1}$}
\psdot(5.40,2.54)\put(5.59,2.44){$R_{j+1}$}
\pspolygon[fillstyle=solid,linewidth=1pt,fillcolor=lightyellow](3.33,4.23)(3.57,4.65)(2.57,4.65)
\psline(3.57,4.65)(2.57,4.65)
\psline(3.33,4.23)(3.33,4.23)
\psline(2.57,4.65)(3.33,4.23)
\psline(3.57,4.65)(3.57,4.65)
\psline(3.33,4.23)(3.57,4.65)
\psline(2.57,4.65)(2.57,4.65)
\pspolygon[fillstyle=solid,linewidth=1pt,fillcolor=lightyellow](2.57,4.65)(2.80,5.08)(3.57,4.65)
\psline(2.80,5.08)(3.57,4.65)
\psline(2.57,4.65)(2.57,4.65)
\psline(3.57,4.65)(2.57,4.65)
\psline(2.80,5.08)(2.80,5.08)
\psline(2.57,4.65)(2.80,5.08)
\psline(3.57,4.65)(3.57,4.65)
\pspolygon[fillstyle=solid,linewidth=1pt,fillcolor=lightyellow](4.10,3.81)(4.33,4.23)(3.33,4.23)
\psline(4.33,4.23)(3.33,4.23)
\psline(4.10,3.81)(4.10,3.81)
\psline(3.33,4.23)(4.10,3.81)
\psline(4.33,4.23)(4.33,4.23)
\psline(4.10,3.81)(4.33,4.23)
\psline(3.33,4.23)(3.33,4.23)
\pspolygon[fillstyle=solid,linewidth=1pt,fillcolor=lightyellow](3.33,4.23)(3.57,4.65)(4.33,4.23)
\psline(3.57,4.65)(4.33,4.23)
\psline(3.33,4.23)(3.33,4.23)
\psline(4.33,4.23)(3.33,4.23)
\psline(3.57,4.65)(3.57,4.65)
\psline(3.33,4.23)(3.57,4.65)
\psline(4.33,4.23)(4.33,4.23)
\pspolygon[fillstyle=solid,linewidth=1pt,fillcolor=lightyellow](4.87,3.38)(5.10,3.81)(4.10,3.81)
\psline(5.10,3.81)(4.10,3.81)
\psline(4.87,3.38)(4.87,3.38)
\psline(4.10,3.81)(4.87,3.38)
\psline(5.10,3.81)(5.10,3.81)
\psline(4.87,3.38)(5.10,3.81)
\psline(4.10,3.81)(4.10,3.81)
\pspolygon[fillstyle=solid,linewidth=1pt,fillcolor=lightyellow](4.10,3.81)(4.33,4.23)(5.10,3.81)
\psline(4.33,4.23)(5.10,3.81)
\psline(4.10,3.81)(4.10,3.81)
\psline(5.10,3.81)(4.10,3.81)
\psline(4.33,4.23)(4.33,4.23)
\psline(4.10,3.81)(4.33,4.23)
\psline(5.10,3.81)(5.10,3.81)
\pspolygon[fillstyle=solid,linewidth=1pt,fillcolor=lightyellow](3.57,4.65)(3.80,5.08)(2.80,5.08)
\psline(3.80,5.08)(2.80,5.08)
\psline(3.57,4.65)(3.57,4.65)
\psline(2.80,5.08)(3.57,4.65)
\psline(3.80,5.08)(3.80,5.08)
\psline(3.57,4.65)(3.80,5.08)
\psline(2.80,5.08)(2.80,5.08)
\pspolygon[fillstyle=solid,linewidth=1pt,fillcolor=lightyellow](2.80,5.08)(3.03,5.50)(3.80,5.08)
\psline(3.03,5.50)(3.80,5.08)
\psline(2.80,5.08)(2.80,5.08)
\psline(3.80,5.08)(2.80,5.08)
\psline(3.03,5.50)(3.03,5.50)
\psline(2.80,5.08)(3.03,5.50)
\psline(3.80,5.08)(3.80,5.08)
\pspolygon[fillstyle=solid,linewidth=1pt,fillcolor=lightyellow](4.33,4.23)(4.57,4.65)(3.57,4.65)
\psline(4.57,4.65)(3.57,4.65)
\psline(4.33,4.23)(4.33,4.23)
\psline(3.57,4.65)(4.33,4.23)
\psline(4.57,4.65)(4.57,4.65)
\psline(4.33,4.23)(4.57,4.65)
\psline(3.57,4.65)(3.57,4.65)
\pspolygon[fillstyle=solid,linewidth=1pt,fillcolor=lightyellow](3.57,4.65)(3.80,5.08)(4.57,4.65)
\psline(3.80,5.08)(4.57,4.65)
\psline(3.57,4.65)(3.57,4.65)
\psline(4.57,4.65)(3.57,4.65)
\psline(3.80,5.08)(3.80,5.08)
\psline(3.57,4.65)(3.80,5.08)
\psline(4.57,4.65)(4.57,4.65)
\pspolygon[fillstyle=solid,linewidth=1pt,fillcolor=lightyellow](5.10,3.81)(5.33,4.23)(4.33,4.23)
\psline(5.33,4.23)(4.33,4.23)
\psline(5.10,3.81)(5.10,3.81)
\psline(4.33,4.23)(5.10,3.81)
\psline(5.33,4.23)(5.33,4.23)
\psline(5.10,3.81)(5.33,4.23)
\psline(4.33,4.23)(4.33,4.23)
\pspolygon[fillstyle=solid,linewidth=1pt,fillcolor=lightyellow](4.33,4.23)(4.57,4.65)(5.33,4.23)
\psline(4.57,4.65)(5.33,4.23)
\psline(4.33,4.23)(4.33,4.23)
\psline(5.33,4.23)(4.33,4.23)
\psline(4.57,4.65)(4.57,4.65)
\psline(4.33,4.23)(4.57,4.65)
\psline(5.33,4.23)(5.33,4.23)
\pspolygon[fillstyle=solid,linewidth=1pt,fillcolor=lightyellow](3.80,5.08)(4.03,5.50)(3.03,5.50)
\psline(4.03,5.50)(3.03,5.50)
\psline(3.80,5.08)(3.80,5.08)
\psline(3.03,5.50)(3.80,5.08)
\psline(4.03,5.50)(4.03,5.50)
\psline(3.80,5.08)(4.03,5.50)
\psline(3.03,5.50)(3.03,5.50)
\pspolygon[fillstyle=solid,linewidth=1pt,fillcolor=lightyellow](3.03,5.50)(3.27,5.92)(4.03,5.50)
\psline(3.27,5.92)(4.03,5.50)
\psline(3.03,5.50)(3.03,5.50)
\psline(4.03,5.50)(3.03,5.50)
\psline(3.27,5.92)(3.27,5.92)
\psline(3.03,5.50)(3.27,5.92)
\psline(4.03,5.50)(4.03,5.50)
\pspolygon[fillstyle=solid,linewidth=1pt,fillcolor=lightyellow](4.57,4.65)(4.80,5.08)(3.80,5.08)
\psline(4.80,5.08)(3.80,5.08)
\psline(4.57,4.65)(4.57,4.65)
\psline(3.80,5.08)(4.57,4.65)
\psline(4.80,5.08)(4.80,5.08)
\psline(4.57,4.65)(4.80,5.08)
\psline(3.80,5.08)(3.80,5.08)
\pspolygon[fillstyle=solid,linewidth=1pt,fillcolor=lightyellow](3.80,5.08)(4.03,5.50)(4.80,5.08)
\psline(4.03,5.50)(4.80,5.08)
\psline(3.80,5.08)(3.80,5.08)
\psline(4.80,5.08)(3.80,5.08)
\psline(4.03,5.50)(4.03,5.50)
\psline(3.80,5.08)(4.03,5.50)
\psline(4.80,5.08)(4.80,5.08)
\pspolygon[fillstyle=solid,linewidth=1pt,fillcolor=lightyellow](5.33,4.23)(5.57,4.65)(4.57,4.65)
\psline(5.57,4.65)(4.57,4.65)
\psline(5.33,4.23)(5.33,4.23)
\psline(4.57,4.65)(5.33,4.23)
\psline(5.57,4.65)(5.57,4.65)
\psline(5.33,4.23)(5.57,4.65)
\psline(4.57,4.65)(4.57,4.65)
\pspolygon[fillstyle=solid,linewidth=1pt,fillcolor=lightyellow](4.57,4.65)(4.80,5.08)(5.57,4.65)
\psline(4.80,5.08)(5.57,4.65)
\psline(4.57,4.65)(4.57,4.65)
\psline(5.57,4.65)(4.57,4.65)
\psline(4.80,5.08)(4.80,5.08)
\psline(4.57,4.65)(4.80,5.08)
\psline(5.57,4.65)(5.57,4.65)
\psdot(3.27,5.92)\put(2.67,5.91){$E$}
\psdot(5.57,4.65)\put(5.76,4.55){$F$}
\psdot(3.69,6.69)\put(3.09,6.68){$$}
\pspolygon[fillstyle=solid,linewidth=1pt,fillcolor=lightgreen](3.69,5.69)(3.69,6.69)(3.27,5.92)
\psline(3.69,6.69)(3.27,5.92)
\psline(3.69,5.69)(3.69,5.69)
\psline(3.27,5.92)(3.69,5.69)
\psline(3.69,6.69)(3.69,6.69)
\psline(3.69,5.69)(3.69,6.69)
\psline(3.27,5.92)(3.27,5.92)
\psdot(5.60,10.15)\put(5.04,9.86){$A$}
\psdot(0.00,0.00)\put(0.00,-0.44){$B$}
\psdot(11.21,0.00)\put(11.21,-0.44){$C$}
\psdot(7.00,0.00)\put(7.00,-0.44){$U$}
\psline(5.57,4.65)(7.11,7.45)
\endpspicture
}

\title{Triangle Tiling I: the tile is similar\\ to $ABC$ or has a right angle}         
\author{Michael Beeson}        
\date{\today}          
\maketitle
\newtheorem{theorem}{Theorem}
\newtheorem{lemma}{Lemma}
\newtheorem{corollary}{Corollary}
\newtheorem{definition}{Definition}
\abstract{
An $N$-tiling of triangle $ABC$ by triangle $T$ is a way of writing $ABC$ as a union of $N$ triangles
congruent to $T$, overlapping only at their boundaries.   The triangle $T$ is the ``tile''.
  The tile  may or may not be similar to $ABC$.   This paper is the first of four papers, which 
together seek a complete characterization of the triples $(ABC,N,T)$ such that $ABC$ can be $N$-tiled
by $T$.   In this paper, we consider the case in which the tile is similar to $ABC$,
 the case in which the tile is a right triangle, and the case when $ABC$ is equilateral.
  We use (only) techniques from linear algebra and elementary field theory, as well as elementary geometry 
and trigonometry.  

Our results (in this paper) are as follows:  When the tile is similar to 
$ABC$, we always have ``quadratic tilings'' when $N$ is a square.  If the tile is similar to $ABC$
and is not a right triangle, then $N$ is a square.  If $N$ is a sum of two squares, $N = e^2 + f^2$, then 
a right triangle with legs $e$ and $f$  can be $N$-tiled by a tile similar to $ABC$; these tilings 
are called ``biquadratic''.  If the tile and $ABC$ are 30-60-90 triangles, 
then $N$ can also be three times a square.   If $T$ is similar to $ABC$, these are all the 
possible triples $(ABC, T, N)$.  

If the tile is a right triangle, of course it can 
tile a certain isosceles triangle when $N$ is twice a square, and in some cases when $N$ is six times a square.
  Equilateral triangles can be 3-tiled and 6-tiled
and hence they can also be $3n^2$  and $6n^2$ tiled for any $n$.  We also discovered a family of tilings when $N$ is 3 times a square, which we 
call the ``hexagonal tilings.''  These tilings exhaust all the possible triples $(ABC, T, N)$ in case $T$
is a right triangle or is similar to $ABC$.  Other cases are treated in \cite{beeson-nonexistence, beeson-triquadratics, beeson120}.  
}

\section{Examples of Tilings}      
We consider the problem of cutting a triangle into $N$ congruent triangles.   Figures 1
 through 4 show that,
at least for certain triangles,  this can be done with $N = 3$, 4, 5, 6, 9, and 16.   Such a configuration 
is called an $N$-tiling.
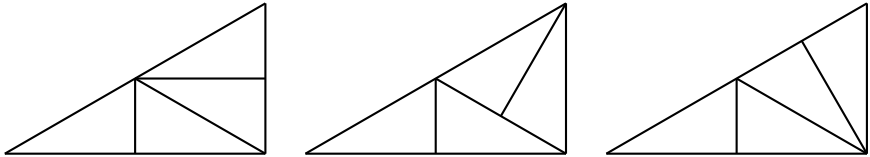
\begin{figure}[ht]
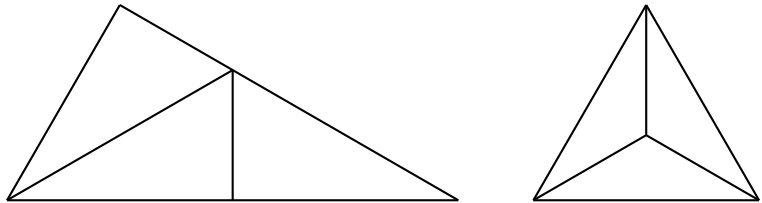
    
\caption{Two 3-tilings   }
\label{figure:3-tilings}
\hskip 2.5cm
\ThreeTilingA
\hskip 1cm
\ThreeTiling
\end{figure}

\begin{figure}[ht]    
\caption{A 4-tiling, a 9-tiling, and a 16-tiling}
\label{figure:nsquared}
\hskip 2cm
\FourTiling  
\hskip 1cm
\NineTiling
\hskip 1cm
\SixteenTiling
\end{figure}

\begin{figure}[ht]
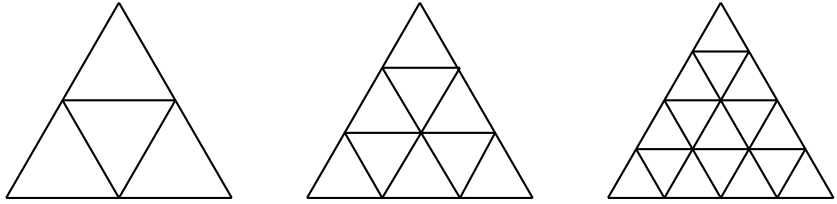
     
\caption{Three 4-tilings}
\label{figure:4-tilings}
\hskip 1.7cm
\FourTilingA
\FourTilingB
\FourTilingC
\end{figure}

\begin{figure}[ht]
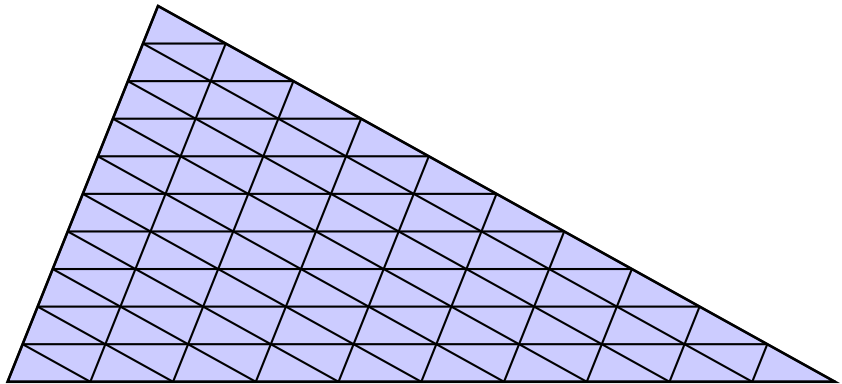

\caption{A quadratic tiling of an arbitrary triangle}
\label{figure:bigquadratic}
\FigureBigQuadratic
\end{figure}

The method illustrated for $N=4$ ,$9$, and $16$ clearly generalizes to any perfect square $N$. 
While the exhibited 3-tiling, 6-tiling, and 5-tiling clearly depend on the exactly angles of the triangle, 
{\em any} triangle can be decomposed into $n^2$ congruent triangles by  drawing $n-1$ lines, parallel to each edge and
 dividing the other two edges into $n$ equal parts.  Moreover, the large (tiled) triangle is similar to the 
small triangle (the ``tile'').  We call such a tiling a {\em quadratic tiling.} Fig.~\ref{figure:bigquadratic}
illustrates a quadratic tiling of an arbitrary triangle.

 It follows that if we have a tiling of a triangle $ABC$ into $N$ congruent triangles, and $m$ is any 
integer,  we can tile $ABC$ into $Nm^2$ triangles by subdividing the first tiling, replacing each of the $N$ triangles by $m^2$ smaller ones.
Hence the set of $N$ for which an $N$-tiling of some triangle exists is closed under multiplication by squares.

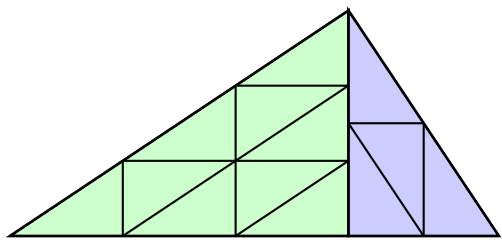
\begin{figure}[ht]     
\caption{Two 5-tilings}
\label{figure:5-tilings}
\hskip 2cm
\FiveTiling  
\hskip 1 cm
\FiveTilingB
\end{figure}

\begin{figure}[ht]
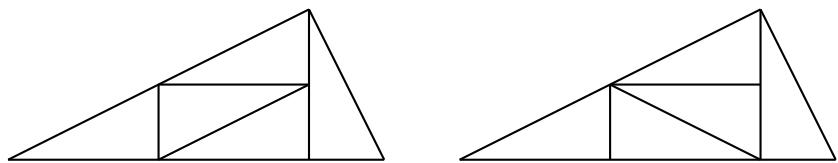
    
\caption{A 13-tiling}
\label{figure:13-tiling}
\hskip 2cm \ThirteenTiling
\end{figure}

Let $N$ be of the form $n^2 + m^2$.  Let triangle $T$ be a right triangle with perpendicular sides $n$ and $m$, say with $n \ge m$.
Let $ABD$ be a right triangle with base $AD$ of length $m^2$, the right angle at $D$ and altitude $mn$, so side $BD$ has length $mn$.  Then $ABD$ can be decomposed into $m$ triangles congruent to $T$, arranged with their short sides
(of length $m$) 
parallel to the base $AD$. Now, extend $AD$ to point $C$, located $n^2$ past $D$. Triangle $ADC$ can be tiled
with $n^2$ copies of $T$, arranged with their long sides parallel to the base.  The result is a tiling of 
triangle $ABC$ by $n^2 + m^2$ copies of $T$.  The first 5-tiling exhibited in Fig.~\ref{figure:5-tilings} is the simplest 
example, where $n=2$ and $m=1$.   The case $N = 13 = 3^2 + 2^2$ is illustrated in Fig.~\ref{figure:13-tiling}.  We call these tilings ``biquadratic.''
More generally, a {\em biquadratic tiling} of triangle $ABC$ is one in which $ABC$ has a right angle at $C$, and can be divided by an altitude 
from $C$ to $AB$ into two triangles, each similar to $ABC$, which can be tiled respectively by $n^2$ and $m^2$ copies of a triangle similar to $ABC$.  
The second 5-tiling in Fig.~\ref{figure:5-tilings}  shows that this can be sometimes be done more generally than by combining two quadratic tilings.  

\begin{figure}[ht]
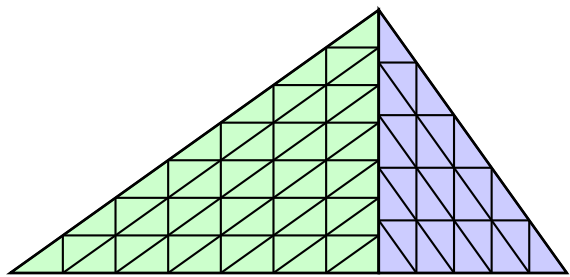

\caption{Biquadratic tiling with $N=74 = 5^2 + 7^2$}
\label{figure:biquadratic}
\FigureBiquadratic
\end{figure}

A larger biquadratic tiling, with  $n=5$ and $m=7$ and hence $N=74$, is shown in Fig.~\ref{figure:biquadratic}. 

If the original triangle $ABC$ is chosen to be isosceles,
then each of the $n^2$ triangles can be divided in half by an altitude;  hence any isosceles triangle can be decomposed into $2n^2$ congruent 
triangles.   If the original triangle is equilateral,  then it can be first decomposed into $n^2$ equilateral triangles, and then these 
triangles can be decomposed into 3 or 6 triangles each,  showing that any equilateral triangle can be decomposed into $3n^2$ or $6n^2$
congruent triangles.  These tilings are neither quadratic nor biquadratic.
For example we can 12-tile an equilateral triangle in two different ways,  starting with a 3-tiling and then subdividing each triangle 
into 4 triangles (``subdividing by 4''),  or starting with a 4-tiling and then subdividing by 3.

\begin{figure}[ht]
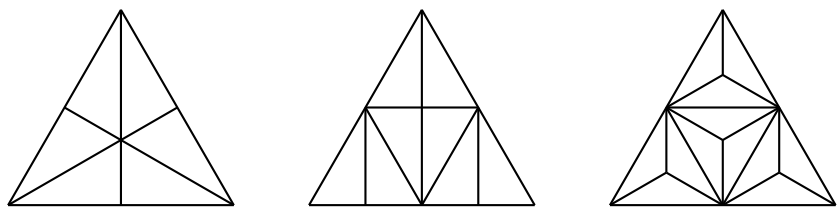
   
\caption{A 6-tiling, an 8-tiling, and a 12-tiling}
\label{figure:6-8-12-tilings}
\hskip 2cm
\SixTiling
\hskip 1cm
\EightTiling
\hskip 1cm
\TwelveTiling
\end{figure}

 Examples like these led us to the following definitions:

A tiling $E$ of triangle $ABC$ (with tile $T_2$) is a {\em subtiling} of another tiling $F$ of $ABC$ (with tile $T$), if $T$ can be tiled by the tile $T_2$ and the tiling $E$ is obtained by tiling each copy of $T$ in $F$ with triangle $T_2$.  It is not required that the {\em same}
tiling be used for each copy of $T$.  For example,  we could take $F$ to be one of the two five-tilings, and then tile each of the tiles 
in that tiling by one of its two five-tilings.  In this way we can obtain 64 different 25-tilings, none of them quadratic.   

A tiling of $ABC$ is called {\em composite} if it is a subtiling of some tiling.  It is called {\em prime} if it is 
not composite.  Note that a quadratic $N^2$-tiling is prime if and only if $N$ is a prime number.

The examples above do not exhaust all possible tilings, even when $N$ is a square.   For example, Fig.~\ref{figure:9-tiling} shows a 9-tiling that is not produced by 
those methods:

\begin{figure}[ht]
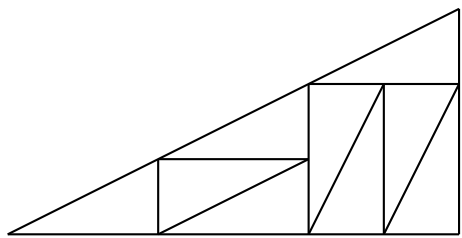
    
\caption{Another 9-tiling}
\label{figure:9-tiling}
\hskip 4cm
\NineTilingA
\end{figure}

There is another family of $N$-tilings, in which $N$ is of the form $3m^2$, and both the tile and the tiled triangle are 30-60-90 triangles.  We call these the ``triple-square'' tilings. 
The case $m=1$ is given in Fig.~\ref{figure:3-tilings}; the case $m=2$ makes $N=12$.  There are two ways to 12-tile a 30-60-90 triangle with 30-60-90 triangle.  
One is to first quadratically 4-tile it, and then subtile the four triangles with the 3-tiling of Figure 1.  This produces the first 
12-tiling in Fig.~\ref{figure:12-tilings}.  Somewhat surprisingly, there is another way to tile the same triangle with the same 12 tiles, also shown in Fig.~\ref{figure:12-tilings}; the second 
tiling is prime.
The next member of this family is $m=3$, which makes $N=27$.  Two 27-tilings are shown in Fig.~\ref{figure:27-tilings}; the first obtained by subtiling a quadratic tiling,
and the second one prime.   Similarly, there are two 48-tilings (not shown).

\begin{figure}[ht]    
\caption{Two 12-tilings   }
\label{figure:12-tilings}
\hskip 3cm
\TwelveTilingA
\hskip 2cm
\TwelveTilingB
\end{figure}

\begin{figure}[ht]
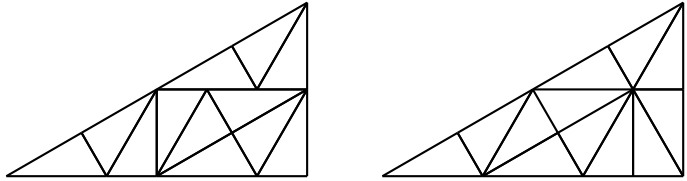
    
\caption{Two 27-tilings   }
\label{figure:27-tilings}
\hskip 3cm
\TwentySevenTilingA
\hskip 2cm
\TwentySevenTilingB
\end{figure}

\begin{figure}[ht]
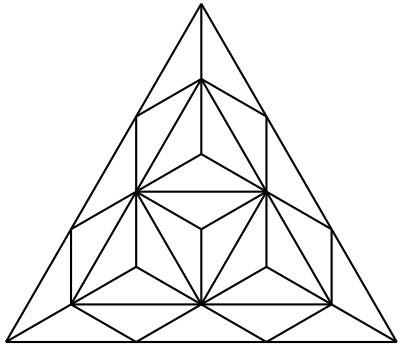
    
\caption{A prime 27-tiling}
\label{figure:prime27}
\hskip 5cm
\TwentySevenTiling
\end{figure}

\begin{figure}[ht]
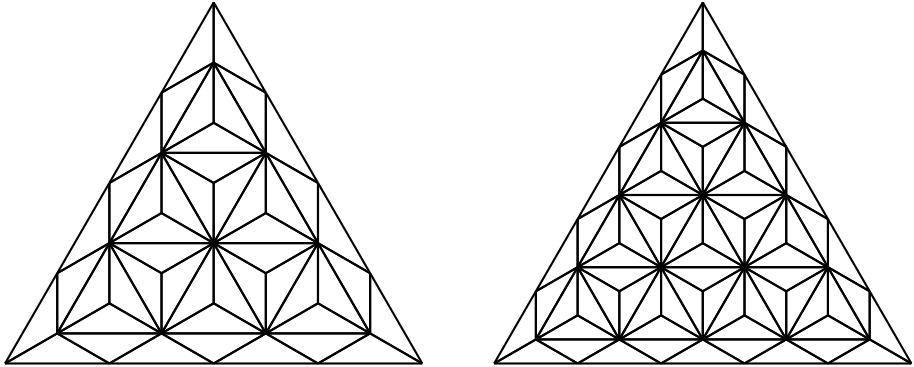
    
\caption{$3m^2$ (hexagonal) tilings for $m=4$ and $m=5$}
\label{figure:hexagonaltilings}
\hskip 1.5cm
\FortyEightTiling
\hskip 0.5cm
\OneHundredTwentyFiveTiling
\end{figure}

Until October 12, 2008,  we did not know any more complicated tilings than those illustrated above (and
there also none in \cite{soifer}).    Then we found the 
beautiful 27-tiling shown in Fig.~\ref{figure:prime27}.  This tiling is one of a family of $3k^2$ tilings (the case $k=3$).
The next case is a 48-tiling, made from six hexagons (each containing 6 tiles) bordered by 
4 tiles on each of 3 sides.  In general one can arrange $1+2+\ldots+k$ hexagons in bowling-pin fashion, and add $k+1$ tiles 
on each of three sides, for a total number of tiles of $6(1+2+\ldots+k)+3(k+1) = 3k(k+1) + 3(k+1) = 3(k+1)^2$. 
Fig.~\ref{figure:hexagonaltilings} shows more members of this family, which we call the ``hexagonal tilings.''%
\footnote{In January, 2012, I bought a puzzle at the exhibition at the AMS meeting, which contained the 
tiling in Fig.~\ref{figure:prime27} as part of a tiling of a larger hexagon.  The tiling is attributed 
to Major Percy Alexander MacMahon (1854-1929) \cite{macmahon}.}

In the third paper in this series, \cite{beeson-triquadratics}, we exhibit a new family of tilings called the 
``triquadratic tilings.''   The tiles and triangles $ABC$ involved do not fall into the cases
considered in this paper.   The main result of this paper is that,  when the tile is similar to 
$ABC$, or the tile is a right triangle, or $ABC$ is equilateral,  there are no more triples 
$(ABC,N,T)$ than we have mentioned.

\section{Previous work}
The examples given in Figures 1 through 6  are well-known.  They have been discussed,
in particular, in connection with ``rep-tiles'' \cite{golomb}.   A ``rep-tile'' is a set of points $X$ in the plane (not necessarily just a triangle) that
can be dissected into $N$ congruent sets, each of which is similar to $S$.
 An $N$-tiling in which the tiled triangle $ABC$ is similar to the 
triangle $T$ used as the tile is a special case of this situation.  That is the case, for example, for the $n^2$ family and the $n^2 + m^2$ family,  but not for the 3-tiling, 6-tiling, or the 12-tiling exhibited above.   Thus the concepts of an $N$-tiling and rep-tiles overlap, but neither subsumes the other.    The paper \cite{goldberg} also contains a diagram showing the $n^2$ family of tilings, but the problem considered there is different: one
is allowed to cut $N$ copies of the tile first, before assembling the pieces into a large figure,  but the large figure must be similar to the original
tile.  The two books \cite{boltyanski} and \cite{boltyanski2} have tantalizing titles, but deal with other problems. 

Only after completing the work in this paper did I encounter
 Soifer's book \cite{soifer},  when the second edition came out,  although the first edition had been out for 19 years.
 The book contains the observation that 
if the tile $T$ is similar to the tiled triangle then $\sqrt N$ is an eigenvalue of a certain matrix, so that observation
is, as it turns out, not new.  The book, however, does not contain any examples of tilings beyond the quadratic tilings,
 though it gives an indication 
that at least the biquadratic tilings were known, since it says that the 1989 Russian Mathematical Olympiad contained the problem to show that 
if $N$ is a sum of two squares then there is a triangle that can be $N$-tiled.   Soifer states (p. 48)  the open problem solved in this series of three papers, and says that Paul \Erdos\ offered a \$25 prize for the first solution.  He does not state where or when \Erdos\ mentioned these 
problems.  The problem statement is:  Find all positive integers $N$ such that at least one triangle can be cut into $N$ 
triangles congruent to each other.  This is Soifer's ``Problem 6.7.''

Soifer also states some related problems.  His ``Problem 6.5'' is:  For each triangle $ABC$, find all positive integers $N$ such that 
$T$ can be cut into $N$ triangles congruent to each other, and the number of distinct partitions of $T$ into $N$ congruent 
triangles.   Actually this is two problems, and we have made serious progress on them both 
in these papers.  Given a triangle $ABC$,
our results succinctly describe the pairs $(N, T)$ such that $ABC$ can be $N$-tiled by triangle $T$,
with the caveat that when the tile $T$ has a $120^\circ$ angle and is not isosceles,  
we have not yet (as of April, 2012) proved that there are no $N$-tilings of $ABC$ by $T$ (except, of 
course, if $ABC$ is similar to $T$).

Describing the pairs $(N,T)$ such that $T$ can $N$-tile $ABC$ is better
than just finding the possible $N$,  but not as good as completely classifying and counting the tilings.   Soifer says that his Problem 6.5 is ``open and very difficult.''

Soifer's ``Problem 6.6'' is also a \$25 \Erdos\ problem:  Find (and classify) all triangles that can only be cut into $n^2$
congruent triangles for any integer $n$.  We have solved this problem, except for the unsolved 
case when has a $120^\circ$ angle and is not isosceles.    The solution is presented in \cite{beeson-triquadratics}.  
  In fact, our main theorem, plus the conjectured absence of tilings 
by non-isosceles tiles with a $120^\circ$ angle,  would yield
 the stronger statement:  Given triangle $ABC$,  a necessary and sufficient 
  condition for $ABC$ to have $N$-tilings only when $N$ is a square is that $ABC$ be not isosceles, and not a 
  right triangle whose angles have rational tangents, and not of the form $N=2K^2-M^2$ with $K$ dividing $M^2$.   
We note that it is not the case that
  for triangles $ABC$ meeting these conditions, any tiling must be the quadratic tiling; there can still be 
integral relations between the sides, such as $3a=2b$, while the angles have irrational tangents.  See 
Lemma~\ref{lemma:nonquadratic} below, and the accompanying Fig.~\ref{figure:nonquadratic}, for an example.
  
Soifer claims, without publishing a proof,
 that if the sides and angles of $ABC$ are integrally independent, then $ABC$ 
admits only quadratic tilings.   It follows (as he points out) that the perfect squares are exactly the $N$ for 
which {\em every} triangle $ABC$ can be $N$-tiled by some triangle.

Dima Fan-Der-Flaas informed me that the problem of finding an $N$-tiling of some triangle when $N=1989$  was posed on the 
Russian  Mathematical Olympiad in 1989; it was solved by a few students,  who had to discover what we call here the ``biquadratic 
tilings'',  and realize that 1989 is a sum of two squares and the relevance of that fact.  I would like to thank Dima for his 
careful reading of parts of an early draft of this paper and \cite{beeson-nonexistence}.

\section{Definitions, notation, and some simple lemmas}
We give a mathematically precise definition of ``tiling'' and fix some terminology and notation.
Given a triangle $T$ and a larger triangle $ABC$,  a ``tiling'' of triangle $ABC$ by triangle $T$ is 
a set of triangles $T_1,\ldots,T_n$ congruent to $T$, whose interiors are disjoint, and the closure of whose union is triangle $ABC$. 

Let $a$, $b$, and $c$ be the sides of the tile $T$, and angles $\alpha$, $\beta$, and $\gamma$ be the angles 
opposite sides $a$, $b$, and $c$.
The letter ``$N$'' will always be used for the number of triangles used in the tiling.  An $N$-tiling of $ABC$ is a tiling that uses
$N$ copies of some triangle $T$.  
The meanings of $N$, $\alpha$, $\beta$, $\gamma$, $a$, $b$,$c$, $A$, $B$, and $C$ will be fixed throughout   these three papers.  
In this paper and the next, we  assume $\alpha \le \beta \le \gamma$;  in the third paper
we only assume $\alpha \le \gamma$ and $\beta \le \gamma$, but not $\alpha \le \beta$.

  An {\em interior  vertex}
in a tiling of $ABC$ is a vertex of one of $T_i$ that does not lie on the boundary of $ABC$.  A {\em boundary vertex}
is a vertex of one of the $T_i$ that lies on the boundary of $ABC$. A ``strict
vertex'' of the tiling is a vertex of one of the $T_i$ that does not lie on the interior of an edge of another $T_j$.
A ``strict tiling'' is one in which no $T_i$ has a vertex lying on the interior of an edge of another $T_j$, i.e. every vertex is strict.
For example, the biquadratic tilings (illustrated above for $N=5$, 13, and 74) are not strict, but all the other tilings shown above are strict.
 A ``non-strict vertex'' $V$ is one that lies on an edge of $T_j$, with 
$T_j$ on one side of the edge and (more than one) $T_i$ having vertex $V$ on the other side.

 Another way of describing a vertex is to say it is of ``type $\pi$'' or ``type $2\pi$'',
or ``has angle sum $\pi$'' or ``has angle sum $2\pi$''.   Strict interior vertices have angle sum $2\pi$, or are 
``of type $2\pi$'';   boundary 
vertices and non-strict interior vertices have angle sum $\pi$, or are ``of type $\pi$.''

By the law of sines we have 
$$ \frac a {\sin \alpha} = \frac b {\sin \beta} = \frac c {\sin \gamma} $$
Up to similarity then we may assume
\begin{eqnarray*}
a &=& \sin \alpha \\
b &=& \sin \beta \\
c &=& \sin \gamma 
\end{eqnarray*}

\begin{definition}\label{definition:maximalsegment}
A {\bf \em maximal segment} in a tiling is a line segment $S$ 
  contained in the union of the edges of the tiling, that cannot be extended to a longer line segment 
  so contained.
\end{definition}
 
Let $S$ be a maximal segment.  
Since there are triangles on each side of $S$, there are triangles on each side of $S$ at every point of $S$ (since
$S$ cannot extend beyond the boundary of $ABC$).  Hence the length of $S$ is a sum of lengths of sides of triangles $T_i$
in two different ways (though the summands may possibly be the same numbers in a different order).  Let us assume 
for the moment that the summands are not the same numbers.  Then 
it follows that some linear relation of the form 
$$ pa + qb + rc = 0$$
holds, with $p$, $q$, and $r$ integers not all zero (one of which must of course be negative),  and the sum of the absolute values
of $p$, $q$, and $r$ is less than or equal to $N$, since there are no more than $N$ triangles.  
This is called an ``edge relation.''   Of course since $a = \sin \alpha$, etc.,  we can express an 
edge relation in terms of the sines of the angles.  We can do that without mentioning $\gamma$: 
Since $\gamma = \pi - (\alpha + \beta)$ we have $\sin(\gamma) = \sin(\alpha + \beta)$, so 
$$ p \sin \alpha + q \sin \beta + r \sin(\alpha + \beta) = 0.$$

A {\em quadratic tiling} is one in which $N$ is a perfect square, say $N = m^2$, and the tiling is produced by drawing
$m-1$ equally spaced lines parallel to each side, dividing each edge into $m$ equal segments.   In such a tiling,
the tile $T$ is similar to the large triangle $ABC$.   An {\em angle relation} is an equation 
$$p \alpha + q\beta + r \gamma = 2\pi$$
where $p$, $q$, and $r$ are non-negative integers, not all equal.   (Since we always have $\alpha + \beta + \gamma = \pi$,
we do not count that equation or its multiples as an angle relation.)

A {\em split vertex} occurs when two copies of the tile in a triangle share one of the vertices of the large triangle.
Split vertices give rise to ``angle relations'', by which we mean an equation $p\alpha + q\beta + r \gamma = n\pi$,
with integers $p$, $q$, and $r$.   We always have the angle relation $\alpha + \beta + \gamma = \pi$, but 
the question of whether a tiling has or might have additional angle relations will be important in our work.

The following lemma is simple and fundamental:

\begin{lemma} \label{lemma:rightangle}
If, in a tiling, $P$ is a boundary vertex (or a non-strict interior vertex)
and only one interior edge emanates from $P$,
then both angles at $P$ are right angles and $\gamma = \pi/2$.   
\end{lemma} 
\noindent{\em Proof.}
 If the two angles at $P$ are different, then their sum is less than $\pi$, since the sum of all three angles 
is $\pi$.  Therefore the two angles are the same.  But $2\alpha \le \alpha + \beta < \pi$ and $2\beta \le \beta + \gamma < \pi$.
Therefore both angles are $\gamma$.  But then $2\gamma = \pi$, so 
$\gamma = \pi/2$.

The following result we call ``Euler's equation'',  because it is related to Euler's famous formula for the numbers of vertices, edges, and 
faces of a polygon,  although we have not derived it that way here.

\begin{lemma}
Let a tiling contain $N_b$ boundary vertices,  $N_s$ strict interior vertices, and $N_n$  non-strict interior vertices.  We then have 
\begin{eqnarray} \label{eq:euler}
 N-1  &=&  N_b + N_n  + 2 N_s
\end{eqnarray}
\end{lemma} 

\noindent {\em Proof}. 
We count vertices.  At each strict interior vertex, the sum of the angles of the tiles sharing that vertex is $2\pi$. 
At each non-strict interior vertex and at each boundary vertex (that is, vertex lying on the boundary of triangle $ABC$ but not 
equal to $A$, $B$, or $C$),  the angle sum is $\pi$.   The total angle sum of all $N$ copies of the tile is $N\pi$,
which must be accounted for by the $\pi$ at the vertices of $ABC$, plus the contributions at the other vertices.  

The following lemma identifies those relatively few rational multiples of $\pi$ that have rational tangents or 
whose sine and cosine satisfy a polynomial of low degree over $\Q$.

\begin{lemma}  \label{lemma:euler}
 Let $\zeta = e^{i\theta}$ be algebraic
of degree $d$ over $\Q$,  where
  $\theta$ is a rational multiple of $\pi$,  say $\theta = 2m \pi/n$, where $m$ and $n$ have no common factor.   
\smallskip

Then $d = \varphi(n)$, where $\varphi$ is the Euler totient function.
In particular if $d = 4$, which is the case when $\tan \theta$ is rational and $\sin \theta$ is not,
then $n$ is 5, 8, 10,  or 12; and if $d=8$ then $n$ is 15, 16, 20, 24, or 30.
\end{lemma}

\noindent{\em Remark}. For example, if $\theta = \pi/6$, we have $\sin \theta = 1 / 2$, which 
is of degree 1 over $\Q$.  Since $\cos \theta = \sqrt 3 / 2$, the number
$\zeta = e^{i\theta}$ is in $\Q(i,\sqrt{3})$, which is of degree 4 over $\Q$.   The number $\zeta$ is a 
12-th root of unity, i.e. $n$ in the theorem is 12 in this case;  so the minimal polynomial of $\zeta$ is
 of degree $\varphi(12) = 4$.   This example shows that the theorem is best possible.
\smallskip

\noindent{\em Remark}. The hypothesis that $\theta$ is a rational multiple of $\pi$ cannot be dropped.
For example, $x^4-2x^3+x^2-2x+1$ has two roots on the unit circle and two off the unit circle.
\smallskip

\noindent{\em Proof}.
Let $f$ be a polynomial with rational coefficients of degree $d$ satisfied by $\zeta$.
 Since $\zeta = e^{i 2m\pi/n}$,
$\zeta$ is an $n$-th root of unity, so its minimal polynomial has degree $d = \varphi(n)$,
where $\varphi$ is the Euler totient function.  Therefore $\varphi(n) \le d$.  If $\tan \theta$
is rational and $\sin \theta$ is not, then  $\sin \theta$ has degree 2 over $\Q$, so $\zeta$
has degree 2 over $Q(i)$, so $\zeta$ has degree 4 over $\Q$.  The stated values
of $n$ for the cases $d=4$ and $d=8$ follow from the well-known formula for $\varphi(n)$. That completes the 
proof of (ii) assuming (i).

\begin{corollary} \label{corollary:euler}
If $\sin \theta$ or $\cos \theta$ is rational, and $\theta < \pi$ is a rational multiple of $\pi$,
then $\theta$ is a multiple of $2\pi/n$ where $n$ is 5, 4, 8, 10, or 12.
\end{corollary}

\noindent{\em Proof}. Let $\zeta = \cos \theta + i \sin \theta = e^{i\theta}$.  Under the stated hypotheses,
the degree of $\Q(\zeta)$ over $\Q$ is 2 or 4.  Hence, by the lemma, $\theta$ is a multiple of $2\pi/n$,
where $n = 5$, 8, 10, or 12 (if the degree is 4) or $n = 3$ or $6$ (if the degree is 3).  But the cases
3 and 6 are superfluous, since then $\theta$ is already a multiple of $2\pi/12$.

\begin{lemma}\label{lemma:tworationalsquares} If $N$ is a sum of two squares of rational numbers, 
then $N$ is the sum of two squares of integers.
\end{lemma}

\noindent{\em Proof}.  See for example,  Proposition 5.4.9, p. 314 of \cite{cohen}, where the theorem 
is attributed to Fermat.   The proof given there is a one-paragraph appeal to the Hasse-Minkowski theorem,
which was not available to Fermat. A simpler proof was pointed out to me by Will Sawin (on MathOverflow),
which only uses the theorem that $N$ is a sum of two integer squares if and only no prime $p$ congruent 
to 3 mod 4 divides $N$ to (exactly) an odd power.  (Fermat knew that theorem.)  Here is the proof:
Suppose $N = p^2 + q^2$.  Then, clearing denominators, $w^2 N = u^2 + v^2$ for some integers $w$, $u$, 
and $v$.  Then any prime $p$ congruent to 3 mod 4 divides the right hand side to an even power; and 
hence $p$ also divides $N$ to an even power.  It follows that $N$ is a sum of two integer squares.  That completes
the proof of the lemma.

\begin{lemma}\label{lemma:sumsofsquares} Let $\mathcal S$ be the set of sums of two rational squares.
Then $\mathcal S$ is closed under multiplication and division.  
\end{lemma}

\noindent{\em Proof}. The members of $\mathcal S$ are just the norms $u^2 + v^2$ of complex numbers
$u + iv$ with rational real and imaginary parts.  Thus 
\begin{eqnarray*}
(u^2 + v^2)(p^2 + q^2) &=& \vert u + iv \vert^2  \vert p + i q \vert^2  \\
&=& \vert (u+iv)(p+iq) \vert^2 \\
&=& (up-vq)^2 + (vp + uq)^2 
\end{eqnarray*}
and 
\begin{eqnarray*}
\frac{u^2 + v^2}{p^2 + q^2} &=&\frac{ \vert u + iv \vert^2}{ \vert p + iq \vert^2} \\
&=& \bigg\vert \frac {u + iv}{p + iq} \bigg\vert^2 \\
&=& \bigg(\frac {up + vq} {p^2 + q^2}\bigg)^2 + \bigg( \frac{ vp -uq}{p^2 + q^2} \bigg)^2 
\end{eqnarray*}
That completes the proof of the lemma.

A product of two sums of integer squares is a sum of two integer squares, as is shown directly 
by the proof of the previous lemma, or follows from the previous two lemmas together.

\begin{lemma}\label{lemma:half-integer} Let $N$ be an integer.  If $N/2$ is a sum of two 
rational squares, then $N/2$ is an integer, i.e. $N$ is even.
\end{lemma}

\noindent{\em Proof}. Suppose $N/2 = p^2 + q^2$.  Then $2N = (2p)^2 + (2q)^2$ is also a sum of 
two rational squares.  Then by Lemma~\ref{lemma:tworationalsquares}, there are integers $u$ and $v$
such that $2N = u^2 + v^2$.  If $u$ and $v$ are both even, then the right hand side is divisible 
by 4, so $N$ is even.  If $u$ is even, then $v^2 = 2N-u^2$ is even, so $v$ is even too; so $N$
is even in that case.  Similarly if $v$ is even.  Finally,  if $u$ and $v$ are both odd, then 
$p = u+v$ and $q =u-v$ are even integers, and we have $p^2 + q^2 = 2(u^2 + v^2) = 4N$.  But 
each of $p^2$ and $q^2$ is divisible by 4, so $p^2 + q^2$ is divisible by 8.  Hence $4N$ 
is divisible by 8, which means $N$ is even.  That completes the proof of the lemma.

%

\section{Quadratic and non-quadratic tilings}

In this section we give a simple sufficient condition for a tiling to be quadratic.

\begin{lemma} \label{lemma:quadratic}  Suppose tile $T$ strictly tiles triangle $ABC$.  If the tile $T$ is similar to the 
triangle $ABC$, and there are no angle relations,  then the tiling is a quadratic tiling.
\end{lemma}

\noindent{\em Proof}. 
Note that since there are no angle relations, the three angles $\alpha$, $\beta$, 
and $\gamma$ are pairwise unequal:  for example, if $\alpha = \beta$, then the relation $\alpha + \beta + \gamma = \pi$
implies $2\alpha + \gamma = \pi$, which is an angle relation.  

Since $T$ is similar to triangle $ABC$,  and angle $A$ is the smallest angle of $ABC$, 
angle $A = \alpha$.  Then consider the copy $T_1$ of the tile that shares vertex $A$. Its two sides 
lie on the sides of triangle $ABC$.  We can relabel the vertices $B$ and $C$ if necessary so that 
the angle of $T_1$ at its vertex $P_1$ on side $AB$ is $\beta$, and its angle at its vertex $Q_1$ on side $AC$
is $\gamma$.  

There must be exactly three copies of the tile meeting at $P_1$, and the three angles 
at $P_1$ are (in some order) $\alpha$, $\beta$, and $\gamma$,  because any other vertex behavior gives rise to an 
angle relation.  Let the tiles meeting at $P_1$ be $T_1$, $T_2$, and $T_3$, numbered so that $T_2$ and $T_1$ share 
a side.   That shared side is $a$, since it is opposite angle $A$ in $T_1$.  Then $T_2$ does not have angle $\alpha$ 
at $P_1$, since the $\alpha$ vertex of $T_2$ has to be opposite side $P_1Q_1$.  $T_2$ does not have angle $\beta$ 
at $P_1$, since $T_1$ has angle $\beta$ there, and only one $\beta$ can occur at $P_1$.  Therefore $T_2$ has angle
$\gamma$ at $P_1$.  Therefore $T_3$ has angle $\alpha$ at $P_1$.  Since the tiling is strict, the angle of $T_3$
at its second vertex $P_2$ on side $AB$ must be $\beta$; otherwise the shared sides of $T_2$ and $T_3$ will have 
different length, since the length of that side of $T_2$ is $b$.  But now, we are in the same situation with $T_2$
as we originally were with $T_1$:  the two angles along side $AB$ are $\alpha$ and $\beta$ (in that order).  We 
can argue as before that the three triangles $T_3$, $T_4$, and $T_5$ meeting at $P_2$ have angles $\beta$,$\gamma$,
and $\alpha$ at $P_2$, in that order.  Continuing down side $AB$ in this fashion, we eventually reach a tile 
$T_{2m-1}$ that has $B$ for a vertex.  
Tile $T_{2m-1}$ has its $\beta$ angle at $B$; then the angle of $ABC$ at $B$ must be $\beta$ (not $\gamma$, 
which is {\em a priori} possible since we may have relabeled $B$ and $C$), since if the angle at $B$ is $\gamma$,
it splits, and $\gamma$ is a sum of $\beta$ plus some other angles, giving an angle relation, contrary to the
hypothesis.   
There will be $m$ copies of the tile sharing a side with $AB$; there will 
be $m-1$ vertices $P_,\ldots, P_{m-1}$ along $AB$, each shared by three triangles;  the number of tiles used is 
$2m-1$.  The third vertices of these triangles are points $Q_1,\ldots,Q_{m-1}$, lying on a line parallel to $AB$,
and the last point $Q_{m-1}$ lies on $BC$, since the angle at $B$ is $\beta$.
  The triangle $Q_1CQ_{m-1}$ is thus tiled by the restriction of the 
original tiling to that triangle.   This restricted tiling is still strict and has no angle relations.  By 
induction, we can assume that this restricted tiling is quadratic.  Since it has $m-1$ tiles along side $Q_1Q_{m-1}$,
we have $(m-1)^2 = N-(2m-1)$. Then $N = (m-1)^2 + 2m -1 = m^2$.  That completes the proof.
\smallskip

 {\em Remark}. 
The 5-tiling in Figure \ref{figure:5-tilings} has $T$ similar to $ABC$, but it has an angle relation $2\alpha + 2\beta = \pi$,
and it also has a non-strict vertex.  It is natural to ask if the hypotheses of the lemma can be weakened by dropping 
one or the other of the hypotheses.  Does there exist a strict non-quadratic tiling  in which $T$ is 
similar to $ABC$?  (Angle relations are OK.)   The following lemma shows that the hypothesis that there are 
no edge relations cannot be dropped. 

\begin{lemma}\label{lemma:nonquadratic} There exists a triangle $ABC$ and integer $N = m^2$ such that $ABC$ can 
be $N$-tiled by a tile similar to $ABC$, and there are no angle relations,
and the sines, cosines, and tangents of $\alpha$ and $\beta$ are irrational,  but the tiling is not quadratic.
(There is, however, an edge relation.)
\end{lemma}

\begin{figure}[ht]
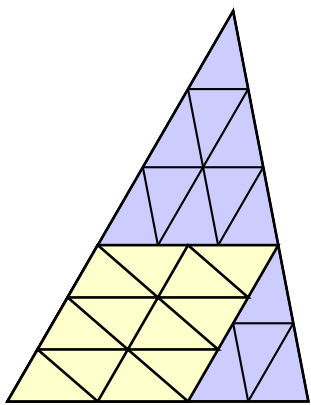

\caption{Here there is an edge relation $3a=2b$, but no angle relations.}
\label{figure:nonquadratic}
\FigureNonQuadratic
\end{figure}

\noindent{\em Proof}.  Take the tile to have $a=2$ and $b=3$.  Then $3a=2b$, which as illustrated in 
Fig.~\ref{figure:nonquadratic}  allows us to construct a nonquadratic tiling.  We are still free to 
choose the angle $\gamma$, or equivalently the side $c$ (as long as they are not too large).  Then $\alpha$,
$\beta$, and $c$ depend continuously and monotonically on $\gamma$, so there are only countably many 
values of $\gamma$ that satisfy an angle relation $p\alpha + q\beta + r\gamma = n\pi$, and only 
countably many values of $\gamma$ that make any of the trig functions of $\alpha$, $\beta$, or $\gamma$
rational.  Pick a value of $\gamma$ that is not one of these countably many values.   That  
completes the proof. 

Whether or not $\gamma$ avoids the countable set where there are angle relations, we can still draw 
a non-quadratic tiling in which all vertices are standard.  Fig.~\ref{figure:nonquadratic} illustrates
such a tiling for $\gamma = 60^\circ$. (Probably there are no angle relations in that case, but we could not 
prove it. If there are, then there are other values of $\gamma$ arbitrarily close to $60^\circ$ with no angle relations.)

The tiling in Fig.~\ref{figure:nonquadratic} has another interesting feature: it has ``one degree of freedom''.
That is,  with its base $AC$ fixed, point $B$ can be varied, but not arbitrarily (as with a quadratic tiling); 
the requirement that $3a=2b$ imposes a restriction on the possible variation of $B$.  The other tilings we have
seen (besides the quadratic tilings) have zero degrees of freedom.

\section{The d-matrix, and a related eigenvalue problem}
Let triangle $ABC$ be tiled by the tile $T$, whose sides are $a$, $b$, and $c$.   Let the sides of $ABC$ be 
$X$, $Y$, and $Z$.  We assume the triangle is labeled so that angles $A$, $B$, and $C$ are listed in non-decreasing order;
hence also $X \le Y \le Z$.  In case triangle $ABC$ is similar to the tile,  this implies that angle $A = \alpha$,
angle $B = \beta$, and angle $C = \gamma$.   

Each side $X$, $Y$, and $Z$ is a linear combination of $a$, $b$, and $c$, the coefficients specifying how many tiles 
share sides of length $a$, $b$, and $c$ with $X$, $Y$, or $Z$.  These nine numbers are the entries of the matrix 
$\dd$, such that 
$$ \vector X Y Z  = \dd\vector a b c $$
If the triangle $ABC$ is similar to the tile, then we have 
$$ \vector X Y Z = \sqrt N  \vector a b c$$
because each side of $ABC$ must be $\sqrt N$ times the corresponding side of the tile $T$, in order that the 
area of $ABC$ can be $N$ times the area of $T$.   Therefore
$$ \dd\vector a b c = \sqrt N \vector a b c.$$
That is, $\sqrt N$ is an eigenvalue of $\dd$, and $(a,b,c)$ is an eigenvector for that eigenvalue. 
If triangle $T$ is isosceles, then $\dd$ is not (yet) uniquely defined.  In that case we have either $a = b$ or 
$b = c$; our convention is to ignore $b$, so that when $T$ is isosceles, the middle column of the $\dd$ matrix is zero.
We will not make use of the $\dd$ matrix when $T$ is equilateral, but for completeness, we define the $\dd$ matrix 
in that case to have non-zero entries only in the first column.   
If $T$ is not isosceles, then the coefficients in the $\dd$ matrix are integers between 0 and $N-1$, inclusive, assuming $N > 2$:  Not all $N$ triangles
can share a side of triangle $ABC$, since if $N > 2$, there would be two adjacent vertices along that side at which 
only two triangles meet; but then by Lemma \ref{lemma:rightangle},  the copy of the tile between those vertices
would have two right angles.

For example, consider the 5-tiling shown in Figure \ref{figure:5-tilings}.  Here the shortest side of the large 
triangle consists of one $c$,  so the top row of the $\dd$ matrix is $0\  0 \ 1$.  The middle side of the large triangle 
consists of two $c$'s,  so the middle row of the $\dd$ matrix is $0\ 0\ 2$.  The longest side of the large triangles 
consists of one $a$ and two $b$'s, so the bottom row is $1\ 2\ 0$.   Thus the $\dd$ matrix for this example is 
$$ \matrix 0 0 1 0 0 2 1 2 0$$
and the eigenvalue equation is 
$$ \matrix 0 0 1 0 0 2 1 2 0  \vector a b c = \sqrt 5  \vector a b c  $$
In this example we have $\alpha = \pi/6$, $\beta = \pi/3$, and $\gamma = \pi/2$, so $a = \sin \pi/6 = 1/2$,
$b = \sin \pi/3 = \sqrt 3/2$, and $c = \sin \pi/2 = 1$.  One can check the eigenvalue equation numerically with these values.

Note that the $\dd$ matrix for a quadratic tiling is $\sqrt N$ times the identity.   At an 
early stage in this research, we conjectured that if $N$ is a perfect
square, say $m^2$,  and $\dd$ is $m$ times the identity, then the tiling is quadratic.  But now we know that 
this conjecture is not correct.  Observe  the tiling in Fig.~\ref{figure:nonquadratic}, which was constructed
starting with a quadratic tiling that, due to the edge relation $3a=2b$, contains a rhombus.  That rhombus 
can be tiled differently by turning the tiles of the quadratic tiling ``the other way.''  Had we started 
with a larger quadratic tiling, such a rhombus could be found entirely in the interior of the tiling, 
leaving the $\dd$ matrix the same as for a quadratic tiling.  The moral is that the $\dd$ matrix only tells 
us about the ``boundary tiling'', not about the interior of the tiling; and the larger $N$ is,  the more 
interior there is.

\section{Tilings with $T$ similar to $ABC$}
In this section, we assume triangle $ABC$ is $N$-tiled by triangle $T$ similar to $ABC$.  In case $N$ is a square,
we have the quadratic tiling of $ABC$; in this section we assume $N$ is not a square.
Let the sides of $T$ be $a$, $b$, and $c$,
in non-decreasing order; these are opposite the angles $\alpha$, $\beta$, and $\gamma$ of $T$.   
We start by disposing of a special case.

\begin{lemma} \label{lemma:equilateral}
Suppose $T$ and $ABC$ are both equilateral, and there is an $N$-tiling of $ABC$ by $T$.  Then
$N$ is a square and the tiling is a quadratic tiling.
\end{lemma}

\noindent{\em Proof}.  Since all the angles of $T$ and $ABC$ are equal, and all the sides of $T$ are equal, there is 
only one way to place tile $T_1$ at vertex $B$.  Along side $BC$ there must be a certain number $m$ of copies of $T$;
hence the side of $ABC$ is $mc$, where $c$ is the side of $X$.  We prove by induction on $m$ that such a tiling 
is a quadratic tiling using $m^2$ triangles.   There are $m$ tiles that share sides with $BC$.  
Call them $T_1,\ldots, T_m$.  This sawtooth-like configuration requires the placement of $m-1$ copies of $T$, one 
between each adjacent pair of triangles $T_1,\ldots,T_m$.  Now we have identified a total of $2m-1$ triangles that participate
in the original tiling, and the remaining triangles tile the smaller equilateral triangle formed by deleting the tiles 
identified so far from $ABC$.  The base of this triangle is smaller than the original base $BC$ by $c$, the side of $T$.
By the induction hypothesis, the tiling of this triangle is quadratic, using $(m-1)^2$ tiles.  Combining this with 
the row of $2m-1$ triangles along $BC$, we have a quadratic tiling with a total of $(m-1)^2 + 2m-1 = m^2$ tiles,
completing the inductive proof.

Next we review the computation of eigenvectors by cofactors.  To find an eigenvector  
of the $\dd$ matrix with eigenvalue $\sqrt N$, consider the matrix $X := \dd- \sqrt{N} I$.  An eigenvector 
can be found by picking any row, and then arranging the cofactors of the elements of that row as a (column) vector.
If these cofactors do not all vanish, then the result is an eigenvector.  (The reader may either verify this 
or just check directly that the particular eigenvalues produced this way below are indeed eigenvectors.)

Now we take up the general case of a tiling $T$ with $ABC$ similar to $T$, when $N$ is not a square.
\begin{lemma} \label{lemma:diagonalzeroes}  
Let triangle $ABC$ be $N$-tiled by tile $T$ similar to $ABC$,
and suppose $N$ is not a square.   Then the diagonal entries of the $\dd$ matrix are zero.
\end{lemma}

\noindent{\em Proof}.
  Since the area of $ABC$ is $N$ times
the area of $T$,  and $T$ is similar to $ABC$,  the sides of $ABC$ are $\sqrt N$ times $a$, $b$, and $c$.   
Then (as discussed in a previous section) we have the
eigenvalue equation 
$$ \dd\vector a b c = \sqrt{N} \vector a b c.$$
The characteristic polynomial $f(x)$ of the $\dd$ matrix, the determinant of $\dd- x I$, is a cubic polynomial with integer 
coefficients, yet has for a zero the number $\sqrt N$.  This is only possible if it factors into a quadratic factor and a 
linear factor.  Since $N$ is not a square, the quadratic factor must be a multiple of $\lambda^2-N$.  The coefficient of 
$x^3$ is $-1$,  and so for some $q$ we have
$$ f(x) = (x^2-N)(q-x) $$
 In general the coefficient of $x^2$ in the characteristic polynomial of any 3 by 3 matrix $\dd$ is the trace of $\dd$,
and the constant term is the determinant of $\dd$.    Hence $q$ is the trace of $\dd$ and $-Nq$ is the determinant of $\dd$.
Since the entries of $\dd$ are non-negative integers, the trace is non-negative, so $q \ge 0$.   

To avoid so many subscripts, we use separate letters for the entries in the $\dd$-matrix, writing it as 
$$ \dd= \matrix p d e g m f h \ell r $$
Since the similarity factor between $ABC$ and $T$ is $\sqrt N$, there cannot be more than $\sqrt N$ tiles with
$a$ sides along $X$, the short side of $ABC$.  That is, $p \le \sqrt N$.  More formally, $a \sqrt N = X = pa + db + ec \ge pa$,
so $p \le \sqrt N$.   Similarly $m \le \sqrt N$ and $r \le \sqrt N$.  Since $N$ is not rational, we have strict 
inequalities:  $p < \sqrt N$, $r < \sqrt N$, and $r < \sqrt N$.  It follows that $pm < N$, etc.

We also note that there is just one tile sharing vertex $A$, where $ABC$ has its $\alpha$ angle.  That tile 
must have its $b$ and $c$ sides along $AB$ and $AC$, or along $AC$ and $AB$, we don't know which.  Thus either 
$f$ and $\ell$ are nonzero, or $m$ and $r$ are nonzero.

Suppose, for proof by contradiction, that $q$, the trace of $\dd$, is not zero.  Then $q = p + m + r > 0$.
Since the three eigenvalues are distinct (because $q$ is rational and $\sqrt N$ is not), the eigenspace corresponding
to $\sqrt N$ is one-dimensional. The eigenvalue equation is
$$ (\dd- \lambda I) \vector u v w = 0$$
or showing the coefficients 
$$ \matrix { p - \sqrt N} d e g {m-\sqrt N} f h \ell {r - \sqrt N} \vector u v w = 0 $$
 We claim that there exists an eigenvector $(u,v,w)$ whose components lie in 
$\Q(\sqrt N)$.  To prove this we will use the cofactor method described above.
 
The resulting eigenvector is $(u,v,w)$,  provided all three components are nonzero,  where
\begin{eqnarray*}
u &=& \dettwo d e {m-\sqrt N} f \\
v &=& -\dettwo {p-\sqrt N} e g f \\
w &=& \dettwo {p-\sqrt N} d g {m-\sqrt N}
\end{eqnarray*}
Although we have not given a proof of the cofactor method's correctness, one can 
easily verify directly that the exhibited vector is indeed an eigenvector for $\sqrt N$; this 
also provides a check that no algebraic mistake has been made.  The fact that all three cofactors are nonzero
is really only needed to conclude directly that the eigenspace of $(u,v,w)$ is one-dimensional; but we 
know that directly in our case since the eigenvalues $\sqrt N$, $-\sqrt N$, and $q$ are distinct.  It therefore
suffices to check that one of the cofactors $u$, $v$, $w$ is nonzero; then the others must automatically be nonzero
because $(u,v,w)$ is a nonzero multiple of $(a,b,c)$.   But we give direct proofs that all three cofactors are nonzero
anyway, as it takes only one more paragraph.

We have $u = df-em + e\sqrt N$.  If $u=0$ then $e=0$ and hence $df = 0$.   If $v=0$ then similarly $f=0$ and $eg=0$.
Finally if $w=0$ then $p+m=0$ and hence $p=m=0$, so $N = dg$.  

Assume, for proof by contradiction, that $w=0$ . Then 
\begin{eqnarray*}
(m - \sqrt N)(p - \sqrt N) &=& dg \\
mp + N -\sqrt N(p + m) &=& dg
\end{eqnarray*}
Since $\sqrt N$ is  irrational this means $p+m=0$, and since $p$ and $m$ are nonnegative, that implies $p=0$
and $m=0$.  Hence $dg = N$.  But $d \le (a/b) \sqrt N$, with equality implying that $p=d=0$, 
and $g \le (b/a)\sqrt N$ with equality implying $m=f=0$.  Since $dg = N$, equality must hold in both
inequalities.  Hence $d = (a/b) \sqrt N$ and $g = (b/a) \sqrt N$
and $p=e=m=f=0$.  But we showed above that either $m$ and $r$ are both nonzero, or $f$ and $\ell$ are both 
nonzero.  That is now contradicted by $m=f=0$.   This contradiction shows that $w \neq 0$. 

Next we give the proof that $u \neq 0$; as remarked above, this is technically superfluous, but still it is 
interesting because the proof we give is not simply an abstract argument about projecting onto the one-dimensional eigenspace.
Assume, for proof by contradiction, that $u=0$.  Then $df-em + e \sqrt N$ = 0.  Since 
$\sqrt N$ is not rational,  $e=0$ and $df-em=0$.  Then $df=0$, so $d=0$ or $f=0$.  If $d=0$,
then since both $d$ and $e$ are zero,  side $X$ of triangle $ABC$ is composed of all $a$ sides 
and $X = (p - \sqrt N)a$.  But since $\sqrt N$ is the similarity factor between $T$ and $ABC$, we have 
$X = \sqrt N a$.  Hence $p -\sqrt N = \sqrt N$.  Hence $\sqrt N = p/2$,  so $N^2 = p^2/4$.  Hence $4N^2 = p^2$
and $p$ is even, so $N$ is a square, contradiction.   This contradiction proves $d \neq 0$.   Since 
$d=0$ or $f=0$, we have $f=0$.  Now assume, for proof by contradiction, that $g = 0$.  Then since $f=0$,
side $Y$ is composed entirely of $b$ sides of tiles, so $Y = mb$.  But $Y = b \sqrt N$ since $\sqrt N$ is 
the similarity factor between $T$ and $ABC$.  Hence $m = \sqrt N$, contradiction.  That proves $g \neq 0$.
As shown above, either $f$ and $\ell$ are both nonzero or $m$ and $r$ are both nonzero.  But $f=0$.  Hence
both $m$ and $r$ are nonzero.   Now 
\begin{eqnarray*}
X &=& pa + db \mbox{\qquad since $e=0$}\\
X &=& \sqrt N a  \\
(\sqrt N - p) a &=& db \\
\frac b a &=& \frac {\sqrt N - p} d \\
Y &=& \sqrt N b \\
Y &=& ga + m b \mbox{\qquad since $f=0$}\\
\frac b a &=& \frac g {\sqrt N-m} \\
\frac ba &=& \frac {\sqrt N -p} d \ = \ \frac g {\sqrt N-m}
\end{eqnarray*}
Cross-multiplying we have
\begin{eqnarray*}
dg &=& N-(m+p)\sqrt N + mp \\
m+p &=& 0   \mbox{\qquad as the coefficient of $\sqrt N$ must be zero}\\
m &=& p \ =  \ 0  \mbox{\qquad  as $m$ and $p$ are nonnegative}
\end{eqnarray*}
But $m$ was proved above to be nonzero.  This contradiction completes the proof that $u \neq 0$.
 
Now assume $v=0$.   Then 
$$ pf - eg +f \sqrt N = 0$$
Since $N$ is irrational we have $f=0$ and $pf = eg$, but since $f=0$ we have $eg=0$.
Hence either $e=0$ or $g=0$.  Assume, for proof by contradiction, that $g=0$.  Then the 
middle side of $ABC$ (corresponding to the middle row) is equal to $mb$ but also to $b\sqrt N$,
so $m = \sqrt N$, contradiction.  This contradiction proves $g \neq 0$.  Hence $e=0$.  Since either
$f$ and $\ell$ are both nonzero or $m$ and $r$ are both nonzero, and we have proved $f=0$, then
$m$ and $r$ are both nonzero.  Now that we have $e=0=f$, and $m \neq 0$, we reach a contradiction 
by the same computation as in the case $u=0$, shown in the series of displayed equations above.
Hence $v \neq 0$.  

 Thus none of the three cofactors is zero.  
That completes the proof that there is an eigenvector $(u,v,w)$  for the eigenvalue $\sqrt N$ 
with components in $\Q(\sqrt N)$.   
Since the eigenspace is one-dimensional, this eigenvector is a (not necessarily rational) multiple of $(a,b,c)$.

Recall that the third eigenvalue of the $\dd$ matrix is the trace $q = p + m + r$.  We can use the cofactor 
method to find an eigenvector for this eigenvalue as well, namely 
\begin{eqnarray*}
V &=& \bigg( \dettwo d e {m-q} f , - \dettwo {p-q} e g f, \dettwo {p-q} d g {m-q} \bigg) \\
&=& \vector {df-em + e(p+m+r)}
            {-pf + eg + f(p+m+r)}
            {pm-dg - (m+p)(m+p+r) + (m+p+r)^2} \\
&=& \vector {df+e(p+r)} {eg + f(m+r)} { pm-dg + r(m+p+r)}             
\end{eqnarray*} 
Technically, it is not an eigenvector until we prove that the components are not zero, but we do not need that
right now; it suffices that it satisfy the eigenvalue equation.  
The eigenvalue equation  $\dd V = (p+m+r) V$ is 
\begin{eqnarray}
 \hskip-0.16in \matrix p d e g m f h \ell r   \vector {df-e(p+r)} {eg + f(m+r)} { pm-dg + r(m+p+r)}  \hskip-0.12in &=& \hskip-0.1in (p + m + r) \vector {df-e(p+r)} {eg + f(m+r)} { pm-dg + r(m+p+r)} \label{eq:secondeigenvalue}
\end{eqnarray}       
The first component of this vector equation is 
\begin{eqnarray*}
p(df - e(p+r)) + d(eg+f(m+r)) + e(pm - dg + r(m+p+r)) &=& (p+m+r) (df-e(p+r))
\end{eqnarray*}
Multiplying out and cancelling like terms, and dividing by 2,  we find 
$$ epm + er(m+p+r) = 0.$$
We argue by cases, according to whether $e=0$ or not.  We first take up 
the case that $e \neq 0$.
Then $pm + r(m+p+r) = 0$.  Since these terms are  nonnegative, they are both zero.  Hence $pm = 0$ and  $r(m+p+r) = 0$.
Hence $r = 0$ or $m+p+r$ = 0.  In either case $r=0$. 
Writing out the third component of the eigenvalue equation, and setting $r=0$, we have
\begin{eqnarray*}
h(df-ep) + \ell (eg+fm) &=& (p+m)(pm-dg) \\
h(df-ep) + \ell (eg+fm) &=& -(p+m)dg \mbox{\qquad since $pm = 0$} \\
hdf - hep + \ell eg + \ell fm &=& -pdg - m dg \\
\end{eqnarray*}
 Now we write out the second component of the eigenvalue equation (\ref{eq:secondeigenvalue}), setting $r=0$:
\begin{eqnarray*}
g(df-ep) + m(eg + fm) + f(pm -dg) &=& (p+m)(eg + fm) \\
gdf - gep + m(eg + fm) + fpm - fdg &=& peg + pfm + m(eg + fm) \\
peg&=& 0 \\
pg &=& 0  \mbox{\qquad since $e \neq 0$ }
\end{eqnarray*}
Assume, for proof by contradiction, that $m \neq 0$.  Then since $mp = 0$ we have $p=0$.
The equation $N(p+m) = hdg + \ell pf + \ell eg$ becomes 
$$Nm = hdg + \ell eg.$$
The third component of the eigenvalue equation becomes, with $p=0$, 
\begin{eqnarray*}
hdf + \ell eg + \ell fm &=& - mdg 
\end{eqnarray*}
The left side is $\ge 0$ and the right side is $\le 0$.  Hence both sides are equal to zero.
Since $m \neq 0$ and $e \neq 0$, we have $dg = 0$ and $hdf = 0$ and $\ell g$ = 0 and $\ell f = 0$.  
We derived above (by observing that the $b$ and $c$ sides of the tile at vertex $A$ lie on the two adjacent 
sides of $ABC$) that either $f$ and $\ell$ are both nonzero or $m$ and $r$ are both nonzero. Since $r = 0$
we must have $f$ and $\ell$ both nonzero.  Hence $\ell f = 0$ is a contradiction.
That contradiction completes the proof that $m=0$.

Now assume, for proof by contradiction, that $p \neq 0$.  
 Then since $pg = 0$ we have $g=0$.    Then the equation $N(p+m) = hdg + \ell p f  + \ell eg$ becomes
$Np = \ell p f$.  Canceling $p$ we have $N = \ell f$.   But as proved above, $\ell f \le N$,
and equality holds if and only if $AC$ is composed only of $c$ sides of tiles and $AB$ is composed
only of $b$ sides.   Therefore we have $h=0$ as well as $g=m=r=0$.
Then since $h=0$ and $r=0$, the long side $AB$ of $ABC$ is composed entirely of $b$ sides of tiles.

If $T$ is isosceles, then by convention the middle column of the $\dd$-matrix is zero.  Since $\ell f = N$,
we now have $\ell \neq 0$, so the middle column is not zero, and $T$ is not isoceles.

Since side $AB$ is composed entirely of $b$ sides of tiles, there are equally spaced vertices $V_0 = A, V_1,\ldots, V_\ell$,
spaced $b$ apart, each one of which is one side of a tile $T_i$.  Tile $T_1$,  which has vertices at $A$ and $V_1$, has
its $\alpha$ angle at $V_0$.   All these tiles have their $\beta$ angles in the interior of $ABC$, and their $\alpha$
and $\gamma$ angles at the $V_i$.  If $\gamma > \pi/2$ then there is only one possible orientation for these tiles,
as two $\gamma$ angles will not fit at any $V_i$.  In that case the angle of the last tile at vertex $B$ must be 
$\gamma$,  contradiction, since the angle there cannot exceed $\beta$, and $\beta \neq \gamma$ since then $T$ would 
be isosceles.  Hence $\gamma \le \pi/2$.

In particular, the tile that shares vertex $B$ has its $b$ side along $AB$.   
Therefore the tile sharing vertex $B$ and part of side $AB$ 
has its $\alpha$ angle at $B$, and the angle $\beta$ at vertex $B$ splits into some number of $\alpha$ angles,
so for some number $J$, we have $\beta = J \alpha$.  Somewhere along $AB$ there must occur a vertex $V_k$ 
at which both the tile $T_k$ and the tile $T_{k+1}$ have angle $\gamma$.  There is not room at $V_k$ for 
a third tile, since $2\gamma + \alpha > \alpha + \beta + \gamma = \pi$.   Hence there are exactly those two 
tiles at $V_k$, and we have $\gamma = \pi/2$.   

Since $\gamma$ is a right angle, we must have $a^2 + b^2 = c^2$.   Since $(u,v,w)$ is a multiple of $(a,b,c)$ we
also have $u^2 + v^2 = w^2$.  We now compute these expressions from the formulas for $(u,v,w)$.  In view of $m=g=0$ we have
\begin{eqnarray*}
u  &=& ef - e\sqrt N \\
v &=& f\sqrt N - fp \\
w &=& N - p \sqrt N
\end{eqnarray*}
Squaring these equations we have 
\begin{eqnarray*}
u^2 &=& e^2(f^2 + N - 2f\sqrt N)\\
v^2 &=& f^2(N+p^2 - 2p\sqrt N) \\
w^2 &=& N^2 +p^2 N - 2p N \sqrt N
\end{eqnarray*}
Setting $u^2 + v^2 = w^2$ we find
\begin{eqnarray*}
e^2 f^2 + e^2 N  + f^2 N + f^2 p^2 - 2(e^2f- f^2p) \sqrt N &=& N^2 + p^2 N - 2pN\sqrt N \\
\end{eqnarray*}
Equating the coefficients of $\sqrt N$ and equating the rational parts, we have
\begin{eqnarray*}
e^2 f - f^2 p &=& pN \\
e^2f^2 + e^2 N + f^2 N + f^2 p^2 &=& N^2 + p^2 N
\end{eqnarray*}
Since $\gamma$ is a right angle, $\alpha + \beta = \pi/2$. Since $\beta = J\alpha$, we have 
$\alpha = \pi/(2(J+1))$, so $\alpha$ is a rational multiple of 2.   We have 
$\tan \alpha = b/a = v/u$,  which belongs to $\Q(\sqrt N)$.  We have 
$$ \cos \alpha = \frac u {u^2 + v^2}$$
which is also in $\Q(\sqrt N)$.  Similarly $\sin \alpha$ belongs to $\Q(\sqrt N)$.  Then $\zeta = d^{i\alpha}$
is of degree 4 over $\Q$, since $\Q(\zeta) = \Q(i,\sqrt N)$.  By Lemma \ref{lemma:euler}, $4(J+1)$ is 5, 8, 10, 
or 12.  Since $5$ and $10$ are not divisible by 4, we have $4(J+1) = 8$ or $10$.  But if $4(J+1) = 8$ then 
$J=1$, while we have $J \ge 2$ since $\beta = J\alpha$.  The only remaining possibility is $4(J+1) = 12$, 
which makes $J = 2$.   Then $\alpha = \pi/6$ and $2\beta = \alpha$, so $\beta = \pi/3$.   
Then $a = \sin \alpha = 1/2$, $b = \sqrt 3 / 2$ and $c = 1$.  But now $\overline{AC} = fc = f$,
and 
\begin{eqnarray*}
 \ell b &=&\overline{AB}\\
&=& \frac 2 {\sqrt 3} \overline{AC}\\
&=& \frac 2 {\sqrt 3} fc
\end{eqnarray*}
Now we put in $c = 1$ and $b = \sqrt 3 / 2$:
\begin{eqnarray*}
\ell \frac {\sqrt 3} 2 &=& \frac 2 {\sqrt 3} f \\
3 \ell &=& 4f
\end{eqnarray*} 
We have $N=\ell f = (4/3)f^2$, so $3N = 4f^2$, so $f$ is divisible by $3$, say $f = 3k$; then $N = 3(2k)^2$ 
is three times a square.
It remains to show that $p=0$;  in fact we claim $p=0$ and $e=0$, so side $AC$ is 
also composed entirely of $b$ sides of triangles.   We have 
\begin{eqnarray*}
\overline{BC} &=& \frac 1 2 \overline{AB} \\
&=& \frac 1 2 \ell b \\
&=& \frac {\ell \sqrt 3} 4 \\
&=& pa + db + ec 
\end{eqnarray*}
Now we put in the values $a = 1/2$, $b = \sqrt 3 /2$, and $c = 1$.
\begin{eqnarray*}
\frac {\ell \sqrt 3} 4 &=& \frac p 2 + d \frac {\sqrt 3} 2 + e
\end{eqnarray*}
This is an equation in $\Q(\sqrt 3)$.   Equating the rational parts we have $0 = p/2 + e$.  Since 
both $p/2$ and $e$ are nonnegative, we have $p=0$ and $e=0$, as claimed.   In particular $p=0$ so the 
diagonal elements are nonzero, which is the conclusion of the theorem; or we could say, in particular $e=0$,
contradicting the assumption $e\neq 0$ and completing the analysis of that case.

Therefore we may now assume $e = 0$.  Remember that $r=0$ was derived only under the assumption $e \neq 0$,
so the equation $r=0$ is no longer in force.    
The third component of the eigenvalue equation (\ref{eq:secondeigenvalue}) is (substituting $e=0$) 
\begin{eqnarray*}
hdf + \ell f(m+r) + r(pm-dg+r(m+p+r)) &=& (p+m+r)(pm + r(m+p+r)) 
 \end{eqnarray*}
Subtracting $r^2(m+p+r)$ from both sides we have 
\begin{eqnarray}
hdf + \ell fm + \ell fr + rpm -rdg &=& (p+m+r)pm + (p+m)r(m+p+r) \nonumber \\
hdf + \ell fm + \ell fr + rpm - rdg &=& (p+m+r)(pm + pr + mr) \label{eq:thirdcomponent8}
\end{eqnarray} 
To get rid of $h$ and $\ell$, we expand the determinant of the $\dd$ matrix by cofactors on the bottom row.
That determinant is $-Nq = -N(p+m+r)$, so we have (remembering $e=0$)
\begin{eqnarray*}
-N(p+m+r) &=& h\dettwo p d g m - \ell \dettwo p e g f + r \dettwo p d g m \\
&=& hdf + \ell pf + rpm - rdg
\end{eqnarray*}
Adding and subtracting $\ell pf$ to the left side of (\ref{eq:thirdcomponent8}) the expression for 
the determinant appears, and we have 
\begin{eqnarray*}
hdf + \ell pf +rpm - rdg + \ell fm + \ell fr - \ell pf &=& (p+m+r)(pm + pr + mr) \\
-N(p+m+r) + \ell fm + \ell fr - \ell pf &=& (p+m+r)(pm +pr +mr)
\end{eqnarray*}
Moving everything to the right side we have 
\begin{eqnarray}
0 &=& (p+m+r)(pm +pr + mr) + (N-\ell f)(m+r) + (N+\ell f) p  \label{eq:4}
\end{eqnarray}
Since $\ell$ is the number of $b$ sides of tiles on the long side $\sqrt N c$ of $ABC$, we have
$\ell b \le \sqrt N c$, or $\ell \le (c/b) \sqrt N$.  Since $f$ is the number of $c$ sides of tiles on $AC$,
whose length is $\sqrt N b$, we have $fc \le \sqrt N b$, or $f \le (b/c) \sqrt N$.  Hence 
\begin{eqnarray*}
\ell f &\le& \bigg( \frac c b \sqrt N \bigg)\bigg( \frac b c \sqrt N\bigg) \\
\ell f &\le& N
\end{eqnarray*}
Hence all the terms on the right of (\ref{eq:4}) are nonnegative.  Hence each of them is zero.
In particular $(N+\ell f)p = 0$; but $N + \ell f > 0$, so $p=0$.  Then then equation becomes
$$ (m+r)mr + (N-\ell f)(m+r) = 0$$
If $m+r = 0$ then $m=0=r$ and the lemma is proved.  Hence we may assume $mr = 0 $ and $N = \ell f$.
But if $N=\ell f$ then we must have equality in the two inequalities $\ell \le (c/b) \sqrt N$ and $f \le (b/c) \sqrt N$.
This implies that side $AC$ is composed only of $c$ sides of tiles and side $AB$ is composed only of $b$ sides of tiles, 
so $g=m=h=r=0$.  In particular $m=r=0$.  That completes the proof of the lemma.
\medskip

We pause to observe that the $\dd$ matrix for a biquadratic tiling, in case $N = m^2 + n^2$,  has the form 
\begin{eqnarray*}
\dd&=& \matrix 0 0 n 0 0 m n m 0
\end{eqnarray*}
which does satisfy the conditions above (as it must).  
The hypothesis that $N$ is not a square is necessary, as shown by the 9-tiling in Fig.~\ref{figure:9-tiling}.  Its $\dd$ matrix is 
$$ \matrix 1 1 0 2 2 0 0 0 3$$
and as predicted, the determinant is zero, but the trace is not zero, and  the characteristic polynomial is $-x(x-3)^2$.

Continuing with the general case of $N$ not a square,
some further conclusions can be drawn about the $\dd$ matrix. We have shown that
$p=m=r= 0$.  The determinant is then given by 
$$ {\rm det\ } \dd = dfh + eg\ell$$
Since the matrix entries are nonnegative,  that means that each of these two terms must contain a zero factor. 
In particular, at most four entries in the $\dd$ matrix are nonzero.

The negated coefficient of $\lambda$ in the characteristic equation is (since the diagonal elements are zero)
the sum of paired products of off-diagonal elements:
\begin{eqnarray} \label{eq:lambdacoef}
 N &=& dg + eh + fl
\end{eqnarray}
But at least one of these three terms will be zero, as shown above.

In view of Lemma~\ref{lemma:diagonalzeroes}, the $\dd$ matrix becomes
\begin{equation} \label{eq:dletters}
 \dd= \matrix 0 d e g 0 f h \ell 0 
\end{equation}
and the matrix equation 
$$ \dd\vector a b c =\sqrt N \vector a b c $$
becomes the three equations 
\begin{eqnarray*}
db + ec &=& \sqrt N a \\
ga + fc &=& \sqrt N b \\
ha + \ell b &=& \sqrt N c
\end{eqnarray*}

\begin{lemma} \label{lemma:penultimate}
 Suppose $ABC$ is $N$-tiled by tile $T$ similar to $ABC$, and $N$ is not a square.  Then 
$\gamma$ is a right angle.
\end{lemma}

\noindent{\em Proof}.  First we note that $T$ and $ABC$ are not equilateral, by Lemma \ref{lemma:equilateral}.  
Next we will prove that $T$ and $ABC$ are not isosceles with $\beta = \gamma$.
Assume, for proof by contradiction, that $\beta=\gamma$.
Then, by our definition of the $\dd$ matrix,
the middle column of the $\dd$ matrix is zero, i.e. $b$ is counted as $c$.   Then we have $d=\ell=0$ and 
\begin{equation}
 \dd= \matrix 0 0 e g 0 f h 0 0 
\end{equation}
That implies that the short side $BC$ of triangle $ABC$ has only $c$ sides of tiles on it, and the long side $AB$ has 
only $a$ sides of tiles on it.   At the vertex $A$, there can only be one tile, since the angle at $A$ is the smallest
angle $\alpha$ so there can be no vertex splitting.  This tile has one side of length $a$ opposite angle $A$ and another
along side $AB$.  Hence $a=b$.    Since $T$ is not equilateral, we must have $b < c$ and $\beta < \gamma$.  This contradicts
the assumption  $\beta=\gamma$, and thus completes the proof by contradiction that $\beta < \gamma$.

Since the $\dd$ matrix has zeroes on the diagonal, no $c$ sides of tiles occur along the longest side $AB$ of triangle $ABC$;
only $a$ and $b$ sides occur there.  
 Since there are no $c$ edges of tile on $AB$, every tile with an edge on $AB$ has a $\gamma$
angle on $AB$. No $\gamma$ angle occurs at the endpoints $A$ and $B$, since the angles there are 
$\alpha$ and $\beta$ respectively, and both are less than $\gamma$.
 By the pigeon-hole principle,
there is a vertex $V$ on $AB$ such that two tiles, say $T_1$ and $T_2$,  have their $\gamma$ angles
at $V$.   At that point we know $\gamma \le \pi/2$.  
Assume, for proof by contradiction, that $\gamma \neq \pi/2$. Then there is at least one 
additional copy  $T^\prime$ of the tile between $T_1$ and $T_2$, sharing vertex $V$.
If $T^\prime$ has 
its $\gamma$ angle at $V$ then there are exactly those three tiles meeting at $V$ (else $\gamma < \pi/3$) and we have 
$\gamma = \pi/3$, and hence $T$ is equilateral, which as noted above is impossible. 
 Hence none of the additional tiles $T^\prime$ meeting at $V$  have a $\gamma$ angle at $V$. 
None of the tiles $T^\prime$ can contribute a $\beta$ angle at $V$ either, since $2\gamma + \beta > \pi$.
Hence there is an angle relation $2 \gamma + p \alpha  = \pi$, where $p$ additional tiles contribute $\alpha$ each to the angle 
sum at $V$, and $p > 0$.   But $2 \gamma + p \alpha > \gamma + \beta + \alpha = \pi$, since $\gamma > \beta$.  This contradiction
completes the proof of the lemma.

\begin{lemma} \label{lemma:fnonzero}
 Suppose $ABC$ is $N$-tiled by tile $T$ similar to $ABC$, and $N$ is not a square.  Then $f = \dd_{21}$ is 
not zero.
\end{lemma}

\noindent{\em Proof}.  Suppose, for proof by contradiction, that $f=0$.  Then the middle row of the $\dd$ matrix 
is $(g,0,0)$, which means that all the tiles along side $AC$ of triangle $ABC$ share their $a$ sides with $AC$.
At vertex $A$, where $ABC$ has its smallest angle $\alpha$, there is exactly one tile $T_1$, with angle $\alpha$ at $a$.
Hence both the side of $T_1$ opposite that angle, and the side shared with $AC$, are equal to $a$.  Thus $T$ is isosceles.
In that case, by convention we have agreed to write the $\dd$ matrix with zeroes in the second column, so 
the $\dd$ matrix has the form
\begin{equation}
 \dd= \matrix 0 0 e g 0 0 h 0 0 
\end{equation}
Now the bottom row is $(h,0,0)$, which means that all the tiles along side $AB$ share their $a$ sides with $AB$.
In particular the tile at vertex $A$ has an $a$ side along $AB$. But we have already seen that its other two sides
are $a$.  Hence the tile is equilateral, contradicting Lemma \ref{lemma:equilateral}, since $N$ is not a square. That completes 
the proof.

There are six letters for coefficients in the $\dd$ matrix (since the three diagonal elements are zero), but for any specific tiling, at most four of those 
coefficients are nonzero.  We will analyze some special cases.   The case corresponding to the biquadratic 
tilings is $d=0$ and $g=0$.  We call that the ``biquadratic case''.  In the biquadratic case the $\dd$ matrix 
has the form
\begin{eqnarray} \label{eq:dsimple}
 \dd= \matrix 0 0 e 0 0 f h \ell 0 
\end{eqnarray}
Equation (\ref{eq:lambdacoef}) now becomes
\begin{eqnarray} \label{eq:lambdacoef2}
eh + \ell f &=& N
\end{eqnarray}
 
We compute the eigenvector in the biquadratic case, using the cofactor method described above.  Let 
$$ X = \matrix {-\sqrt N} 0 e 0 {-\sqrt N} f h \ell {-\sqrt N} $$
Taking the cofactors of the bottom row (notice the minus sign in the second component, which comes from the 
definition of ``cofactor'') we find the eigenvector 
\begin{eqnarray*}
 \bigg( \dettwo 0 e  {-\sqrt N} f,
               -\dettwo {-\sqrt N} e 0 f,
              \dettwo {-\sqrt N} 0  0 {- \sqrt N} \bigg)
  &=&  \vector { e \sqrt N}  { f \sqrt N}  N
\end{eqnarray*}
Note that $e \neq 0$ and $f \neq 0$, since the first two sides of $ABC$ are given by $ec$ and $fc$.
Hence the cofactors do not vanish.   

We claim that the bottom two rows of $\dd- \sqrt N I$, namely $(0,-\sqrt N,f)$ and $(h,\ell,-\sqrt N)$,  are linearly independent.
Indeed, suppose that for some constants $p$ and $q$ we have $p(0,-\sqrt N, f) + q(h,\ell,-\sqrt N) = 0$.  From the first 
component we see that $qh = 0$.    From the third component we see that $pf = q \sqrt N$.
If $q$ is not zero, then $\sqrt N = pf/q$, contradicting the irrationality of $\sqrt N$.  Hence $q = 0$. 
 Hence from the second component, $p \sqrt N = 0$. Hence $p=0$.
 This proves that the bottom two rows of $\dd- \sqrt N I$ are linearly independent.  Hence $\dd- \sqrt N I$ has rank 2;
hence the eigenspace associated with the eigenvalue $\sqrt N$ is one-dimensional.   It follows that the eigenvector computed above
is a multiple of $(a,b,c)$.  That is, for some constant $\mu$ we have
\begin{eqnarray} \label{eq:eigenvector-biquadratic}
\vector a b c &=& \mu \vector { e \sqrt N}  { f \sqrt N}  N
\end{eqnarray}
The constant $\mu$ is an arbitrary scale factor; changing $\mu$ just changes the size of the tile $T$ and the triangle $ABC$
by the same factor.  We are therefore free to choose $\mu$ to suit our convenience.  We choose to take $\mu = \sqrt N$; then we have 
\begin{eqnarray} \label{eq:abcvalues}
\vector a b c &=& \vector e f {\sqrt N}
\end{eqnarray}

\begin{lemma} \label{lemma:gammaright-sumofsquares}
 Let triangle $ABC$ be $N$-tiled by $T$, and suppose $N$ is not a square and $T$ is similar to $ABC$, and 
$d=g=0$ (the biquadratic case). Then $N$ is a sum of squares, specifically $N = e^2 + f^2$ where $e$ and 
$f$ are as above, and $\tan \alpha = e/f$.  In particular $\tan \alpha$ is rational.
\end{lemma}

\noindent{\em Proof}. 
By Lemma \ref{lemma:penultimate}, we have $\gamma = \pi/2$.  By the Pythagorean theorem and (\ref{eq:abcvalues}) 
we see that $\gamma = \pi/2$ if and only if $e^2 + f^2 = N$.  Since $\gamma  = \pi/2$, we have 
$\tan \alpha = a/b$, and by  (\ref{eq:eigenvector-biquadratic}), $\tan \alpha = e/f$. That completes the proof of the lemma.

\begin{lemma}  \label{lemma:vertex-splitting-occurs}
Suppose $ABC$ is $N$-tiled by tile $T$ similar to $ABC$, and $N$ is not a square, and $d=g=0$ (the biquadratic case).
 Then the right angle of $ABC$ is split by the tiling, and the tangents of the other angles of $ABC$ are rational.
\end{lemma}

\noindent{\em Proof}.  We suppose, as always, that the $\gamma$ angle of $ABC$ is at $C$, the $\beta$ angle at $B$, and
the $\alpha$ angle at $A$.
Since the $\dd$ matrix has the form given in (\ref{eq:dsimple}), all the tiles along side $BC$ share their $c$ sides with $BC$
(there are $e$ of them)  and all the tiles along side $AC$ share their $c$ sides with $AC$ (there are $f$ of them).  Suppose,
for proof by contradiction, that the 
vertex at $C$ is not split. Then a single tile shares vertex $C$,  so the tile has two $c$ sides, and hence is isosceles with $b=c$.
But by (\ref{eq:abcvalues}),  we have $c/b = \sqrt N/e$.  Hence if $b=c$ we have $N = e^2$, contradicting the hypothesis that 
$N$ is not a square.  Hence the vertex $C$ is split as claimed.   The tangents of the other two angles are $e/f$ 
and $f/e$, which are rational.  This completes the proof of the lemma.

We now turn to another important case, when $e=0$.   We call this the ``triple-square case'', because it will 
turn out that in this case $N$ must be three times a square.  The following lemma and its proof give a complete analysis 
of this case.

\begin{lemma} \label{lemma:triplesquare}
Suppose $ABC$ is not equilateral and is $N$-tiled by tile $T$ similar to $ABC$,
and $N$ is not a square, and $e=0$ (the triple-square case).
Then $\alpha = \pi/6$, $\beta = \pi/3$, and $N = 3d^2$  is three times a square.
\end{lemma}

\noindent{\em Remark}.  There do exist tilings for each $N$ of the form $3d^2$ that fall under the triple square case,
as we showed in Fig.~\ref{figure:27-tilings} and Fig.~\ref{figure:prime27}.
\smallskip

\noindent{\em Proof}.  Under the hypotheses of the lemma we have
$$ \dd= \matrix 0 d 0 g 0 f h \ell 0.$$
In this matrix, $d$ and $h + \ell$ are not zero, since they represent the number of tiles along $BC$ and $AB$, respectively.
By Lemma \ref{lemma:fnonzero} we have $f \neq 0$. We have
$$ X = \dd- \sqrt N  I = \matrix {-\sqrt N} d 0 g  {-\sqrt N} f h \ell  {-\sqrt N} $$
We will prove that the bottom two rows of $X$ are linearly independent.  If they are linearly 
dependent, then for some $p$ and $q$, we have 
\begin{eqnarray*}
0 &=&  p(g, - \sqrt N,f) + q(h,\ell,\sqrt N) \\
&=& pg + pf + q\ell + qh + \sqrt N (q-p)
\end{eqnarray*}
Since $N$ is not a square, the coefficient of $\sqrt N$ is zero, so $q=p$,  and 
\begin{eqnarray*}
0 &=& pg + pf + q\ell + qh \\
&=& pg + pf + p\ell + ph \\
&=& p(g + f + \ell + h)
\end{eqnarray*}
Since the entries of the $\dd$ matrix are non-negative, and $h+ \ell$ is strictly positive, we conclude $p = q = 0$.
That proves that the bottom two rows of $X$ are linearly independent, so $X$ has rank 2 and the eigenspace of $\sqrt N$ 
is one-dimensional. 
We then compute the eigenvector by the cofactor method.  Taking the cofactors of the bottom row, we 
find the eigenvector
\begin{eqnarray*}
 \bigg( \dettwo 0 e  {-\sqrt N} f,
               -\dettwo {-\sqrt N} 0 g f,
              \dettwo {-\sqrt N} d  g {- \sqrt N} \bigg)
  &=&  \vector { df}  { f \sqrt N}  {N-dg}
\end{eqnarray*}
Since $d \neq 0$ and $f \neq 0$, the first two components are not zero.  From the first row of the 
$\dd$ matrix we have $db \le \sqrt{N} a$, and from the second row, we have $ga < \sqrt{N} b$, 
with strict inequality because $f \neq 0$. 
Hence $dg < \sqrt{N}(a/b) \sqrt{N}(b/a) = N$.   Hence the third component is nonzero.

 Since $(a,b,c)$ is an eigenvector and 
the eigenspace is one dimensional, $(a,b,c)$ is a multiple of this computed eigenvector.  By scaling the triangle 
appropriately we can assume $(a,b,c)$ is actually equal to the computed eigenvector:
$$ \vector a b c = \vector  {df} {f\sqrt N} {N-dg}$$
It follows that $\sin \alpha = a/c = fd/(N-dg)$ is rational and $\tan \alpha = a/b = d/\sqrt N$ is of degree 2 over $\Q$.

According to the first row of the $\dd$ matrix, the tiles along $BC$ have only $b$ sides on $BC$.
Assume, for proof by contradiction, that vertex $B$ does not split.  Then there is a single tile $T_1$ at vertex $B$,
which therefore shares one side with $AB$ and one side with $BC$.  Triangle $T$ is not isosceles, since then by 
definition the $\dd$ matrix would have zeroes in the middle column.  Hence the unique $b$ side of $T_1$ must be 
opposite angle $B$; but it must also lie on $BC$, which is a contradiction.
  Hence vertex $B$ does split.  Therefore for some integer $P$ we have $\beta = P \alpha$.
Since by Lemma \ref{lemma:penultimate}, $\gamma = \pi/2$,we have 
\begin{eqnarray*}
\frac \pi 2 &=& \alpha + \beta \\
&=& \alpha + P \alpha \\
&=& (P+1)\alpha
\end{eqnarray*}
Therefore
$$ \alpha = \frac \pi {2(P+1)} = \frac {2\pi} {4(P+1)}.$$
By Lemma \ref{lemma:euler} we conclude that $4(P+1)$ is one of the numbers $n =3$, 4, 5, 8, 10, or 12 for which $\phi(n) = 4$.
Of these numbers, only 4, 1, and 12 are divisible by 4, which implies $P=2$, since the values $P=0$ and $P=1$ do not correspond
to vertex splitting.  Hence $P=2$ and we have $\beta = 2\alpha$, so $\alpha + \beta = \pi/2 = 3\alpha$, and $\alpha = \pi/6$.
 Hence $\tan \alpha = df/(f\sqrt N) = d/\sqrt N = 1\sqrt 3$.  Hence 
$N = 3d^2$.  
That completes the proof of the lemma.
 
Now we have dealt with the biquadratic case (when $d=g=0$) and the triple-square case (when $e = 0$).  It remains
to show that these are the only two possible cases, when $N$ is not a square and $T$ is similar to $ABC$.  Recall that
$dfh = 0$ and $eg\ell = 0$ since $\det \dd= 0$; that leaves only a few possibilities to consider.  We begin by 
showing that if $d =0$ then we are already in the biquadratic case.

\begin{lemma} \label{lemma:dandg}
Assume triangle $ABC$ is $N$-tiled by $T$, that $N$ is not a square, that $T$ is similar to $ABC$, and that 
$d$ and $g$ are two entries in the $\dd$ matrix of the tiling, in the notation used above (the ones that are zero
in the biquadratic case).  Then $d = 0$ implies $g = 0$, i.e. we are in the biquadratic case as soon as $d=0$.
\end{lemma}

\noindent{\em Proof}.  We have 
$$ X = \dd- \sqrt N I = \matrix {-\sqrt N} d e g {-\sqrt N} f h l {-\sqrt N} $$
By the cofactor method described above we compute the eigenvector
$$ \vector a b c = \vector { df + e \sqrt N} { f \sqrt N + e g} { N - dg} $$
Suppose, for proof by contradiction, that $d = 0$ and $g \neq 0$. 
Then 
$$ \vector a b c = \vector { e \sqrt N } { f \sqrt N + eg} N $$
By Lemma \ref{lemma:penultimate}, $\gamma$ is a right angle, so by the Pythagorean theorem, we have $a^2 + b^2 = c^2$.  That is,
\begin{eqnarray*}
 e^2 N + (f\sqrt N + eg)^2   &=&  N^2 \\
e^2 N + f^2 N + 2egf \sqrt N + e^2 g^2 &=&N^2 
\end{eqnarray*}
Since $N$ is not a square,  the coefficient of $\sqrt N$ is zero;  that is, $egf= 0$.  By hypothesis, $g \neq 0$, so $ef = 0$.
The first row of the $\dd$ matrix is $(0,d,e) = (0,0,e)$,  so $e \neq 0$ because there must be some triangles on the first side of $ABC$.
Therefore $f = 0$.   Then the $\dd$ matrix is 
$$ \matrix 0 0 e g 0 0 h \ell 0 $$
Hence  all the tiles on the middle side $AC$ of $ABC$ have their $a$ side on $AC$, and all the tiles on the hypotenuse $AB$ 
do not have their $c$ side on $AB$.  Consider the tile $T_1$ sharing vertex $A$ (there is only one, since $ABC$ has angle $\alpha$ there).
It has its $a$ side on $AC$ and does not have its $c$ side on $AB$.  Hence its $b$ side is on $AB$ and its $c$ side opposite angle $A$,
which is $\alpha$.  Hence $a=c$ and triangles $T$ and $ABC$ are equilateral, which is a contradiction since $\gamma = \pi/2$.
This contradiction shows that the assumption $d=0$ and $g \neq 0$ is untenable, which completes the proof of the lemma.

\begin{lemma}  \label{lemma:notbiquadraticimpliestriplesquare}
Assume triangle $ABC$ is $N$-tiled by $T$, that $N$ is not a square, that $T$ is similar to $ABC$,
and that $d \neq 0$.  Then $e=0$, i.e. we are in the triple-square case as soon as $d \neq 0$.
\end{lemma}

\noindent{\em Proof}.  We have as in the proof of the previous lemma
$$ \dd= \matrix 0 d e g 0 f h \ell 0$$
$$ X = \dd- \sqrt N I = \matrix {-\sqrt N} d e g {-\sqrt N} f h \ell {-\sqrt N} $$
$$ \vector a b c = \vector {df + e \sqrt N} {eg + f \sqrt N} {N-dg} $$
By Lemma \ref{lemma:penultimate} and the Pythagorean theorem we have
\begin{eqnarray*}
0 &=& c^2 - a^2 - b^2 \\
&=& (N-dg)^2 - (df + e \sqrt N)^2 - (eg + f\sqrt N)^2  \\
&=& -2(def + egf) \sqrt N + \mbox{\ rational} 
\end{eqnarray*}
Since $N$ is not a square and the entries of the $\dd$ matrix are nonnegative integers, we have $def = 0$ and $egf = 0$.
Since $f \neq 0$ by Lemma \ref{lemma:fnonzero}, and $d \neq 0$ by hypothesis, we have $e = 0$.  
That completes the proof of the lemma.

The following theorem completely answers the question, ``for which $N$ does there exist an $N$-tiling
in which the tile is similar to the tiled triangle?''
\begin{theorem} \label{theorem:Similar}
 Suppose $ABC$ is $N$-tiled by tile $T$ similar to $ABC$.  Then either $N$ is a square, or a sum of two squares,
or three times a square.  
\end{theorem}

\noindent{\em Proof}.  Suppose $N$ is not a square.  Then by Lemma \ref{lemma:penultimate}, $\gamma$ is a right angle.
Now consider the $\dd$ matrix.  By Lemma \ref{lemma:diagonalzeroes} the diagonal entries are zero, so 
as stated in (\ref{eq:dletters}) the $\dd$ matrix has the form
$$ \dd= \matrix 0 d e g 0 f h \ell 0 $$
By Lemma \ref{lemma:dandg},  if $d=0$ then also $g=0$, i.e. we are in the 
``biquadratic case''.  Then by Lemma \ref{lemma:gammaright-sumofsquares},  $N$ is a sum of squares.
If $e=0$ then by Lemma \ref{lemma:triplesquare}, $N$ is three times a square and $T$ is a 30-60-90 
triangle.  Finally, Lemma \ref{lemma:notbiquadraticimpliestriplesquare} shows that the cases $d=0$ 
and $e=0$ are exhaustive.   That completes the proof of the theorem.

As for characterizing the triples $(ABC, N, T)$,  we have  
the following  results:
\begin{theorem} \label{theorem:Similar2} 
 Suppose $ABC$ is $N$-tiled by tile $T$ similar to $ABC$.
If $N$ is not a square, then $T$ and $ABC$ are right triangles.  Then either 
\smallskip

(i) $N$ is three times
a square and $T$ is a 30-60-90 triangle, or 

(ii)  $N$ is a sum of squares $e^2 + f^2$,  the right angle of $ABC$ is split 
by the tiling, and the acute angles of $ABC$ have rational tangents $e/f$ and $f/e$,
\smallskip

\noindent
  and these two alternatives are 
mutually exclusive. 
\end{theorem}

\noindent{\em Proof}.  First we prove that $N$ cannot be both a sum of squares and three times a square,
since the equation $x^2 + y^2 = 3z^2$ has no integer solutions.   To see that, we can assume without loss of 
generality that $x$, $y$, and $z$ are not all even.  Note that squares are 
always congruent to 0 or 1 mod 4, so the left side is 0, 1, or 2 mod 4.  Then $z^2$ must be congruent to 
0 mod 4, since if not, the right side is congruent to 3 mod 4.  Hence $z$ is even.  But $x$ and $y$ must also be even to make the left side congruent to 0 mod 4, contradiction.  Hence the equation has no solutions.  Thus the alternatives in the 
theorem are mutually exclusive, as claimed.  The rest of the theorem follows from {lemma:gammaright-sumofsquares}, 
Lemma \ref{lemma:penultimate} and Lemma \ref{lemma:triplesquare}.

Note that the 9-tiling in Fig.~\ref{figure:9-tiling} shows that not every $m^2$-tiling is a quadratic tiling, so we have not 
even classified all the $m^2$-tilings.   Briefly we conjectured that a tiling in which the $\dd$-matrix 
is $m$ times the identity should be a quadratic tiling,  but that is not true.  One can extend the 
9-tiling in Fig.~\ref{figure:9-tiling} by adding more triangles to the right and below,  producing a 25-tiling in which the 
$\dd$-matrix is 5 times the identity.

\section{Tilings of an isoceles triangle by a right triangle}
\begin{figure}[ht]
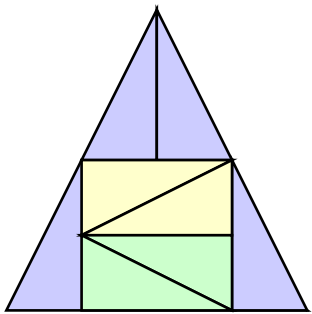

\caption{A tiling that does not include the altitude of $ABC$}
\label{figure:isoscelescounterexample}
\FigureIsoscelesCounterexample
\end{figure}

Let us first review the known examples of tilings of an isosceles triangle $ABC$    by a right triangle $T$.
 There is always the ``double quadratic'' tiling, in which one divides $ABC$ into two halves
by the altitude and then quadratically tiles each half. One might conjecture that {\em any} tiling contains the altitude, 
i.e. is a subtiling of the 2-tiling defined by the altitude.  But this is easily refuted:  consider the second tiling in Fig. 6.  In
that case $ABC$ is equilateral; take the base $BC$ as one of the slanting sides in Fig. 6.  This is not a very 
satisfying counterexample since the triangle can be rotated so that, with a different vertex and base, the tiling 
does contain the altitude.   Fig.~\ref{figure:isoscelescounterexample} gives a counterexample in which
$ABC$ is not equilateral.

Whenever there is an $N$-tiling of the right triangle $ABM$, there is a $2N$-tiling of
the isosoceles triangle $ABC$.  Using the biquadratic tilings (see Fig.~\ref{figure:5-tilings}
and Fig.~\ref{figure:13-tiling} and  triple-square tilings (see Fig~\ref{figure:12-tilings} and Fig.~\ref{figure:27-tilings}),
we can produce $2N$-tilings when $N$ is a sum of squares or three times a sum of squares.  We call these tilings ``double biquadratic''
and ``hexquadratic''.  
  For example,  one has two 10-tilings and two 26-tilings, obtained by reflecting Figs. 4 and 5 about either of the 
sides of the triangles shown in those figures;  and one has 24-tilings and 54-tilings obtained from Figs. 8 and 9.
Note that in the latter two cases, $ABC$ is equilateral.   

\begin{figure}[ht]
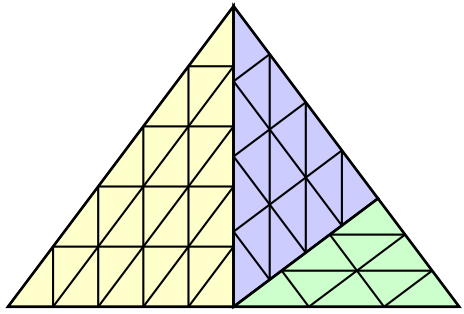

\caption{A tiling related to a Pythagorean triple $a^2 + b^2 = c^2$.}
\label{figure:Pythagorean}
\FigurePythagorean
\end{figure}

In the case when the sides of the tile $T$ form a Pythogorean triple $n^2 + m^2 + k^2 = N/2$, then we can tile one 
half of  $ABC$ with a quadratic tiling and the other half with a biquadratic tiling.  The smallest example is when the tile
has sides 3, 4, and 5, and $N = 50$.  See Fig.~\ref{figure:Pythagorean}.
 One half is 25-tiled quadratically, and the other half is divided into two smaller
right triangles which are 9-tiled and 16-tiled quadratically.  This shows that the tiling of $ABC$ does not have to be 
symmetric about the altitude.  

One looks for other ways to tile the two halves of $ABC$ differently.  We cannot use two different biquadratic tilings,
even if $N/2$ can be written as a sum of squares in two different ways, because the tiles would have to be different shapes.
We cannot use a biquadratic tiling and a triple-square tiling, since $m^2 + n^2 = 3k^2$ has no solutions.  In 
addition to the quadratic-biquadratic
tilings corresponding to Pythagorean triples, we could, for certain isosceles $ABC$, tile the two halves with 
tilings differing by having the tiles in a certain internal square ``turned the other way.'' For example,
one half could be 9-tiled by a quadratic tiling, and the other half 9-tiled as in Fig.~\ref{figure:9-tiling}.

 In this section, our normal convention that $A$ is the smallest angle of $ABC$ is 
temporarily suspended.  Instead, $A$ is the vertex of the isosceles triangle and $BC$ is the base.  Similarly,
we also suspend the convention that $\alpha < \beta$.  The following lemma shows that the base angles of $ABC$
are either $\alpha$ or $\beta$;  we will make the convention for this section that $\beta$ is the base angle.

The area of triangle $ABC$ must be $N$ times the area of the tile $T$.   By hypothesis, the triangle $T$ is similar to half the 
triangle $ABC$;   the similarity factor is the square root of the ratio of the half the area of $ABC$ to the area of $T$,
namely $\sqrt{\frac N 2}$.   The sides of $T$ are $a = \sin \alpha$, $b = \cos \alpha$, and $c = 1$. 
Let $\overline{AB}$ be the length of $AB$.  Since $AB$ is opposite the right angle of one of the halves of $ABC$ that is similar to $T$, we have
\begin{eqnarray*}
\overline{AB}  &=& \sqN \sin \gamma = \sqN \mbox{\qquad since $\gamma = \pi/2$} \\
\frac {\overline{BC}} 2 &=& \sqN \sin \alpha 
\end{eqnarray*}
On the other hand, $\overline{AB}$ and $\overline{BC}$ must be integer linear combinations of $a$, $b$, and $c$.  So we have, for 
some non-negative integers $p$, $d$, and $e$, that $\overline{AB} = pa + db + ec$.  Since $c=1$ and $\overline{AB} = \sqN$, we have 
\begin{equation}\label{eq:pqrsqrtN}
p \sin \alpha  + d \cos \alpha  + e = pa + db + e = \sqN
\end{equation}
The notation in this equation will be used throughout this section.  In this section we also assume that 
$X = \overline{AB} = Y = \overline{AC}$ are the equal sides, and $Z = \overline{BC}$ is the base of $ABC$,
and the $\dd$ matrix satisfies 
$$ \dd\vector abc = \vector X Y Z$$
This differs from the convention in other sections since we do not know the relative sizes of the angles at 
the vertices $A$ and $B$.   With this convention,  $(p,d,e)$ is the first row of the $\dd$ matrix. 

\begin{lemma}  \label{lemma:degree2}
Let $ABC$ be an isosceles triangle with base $BC$,  $N$-tiled by a right triangle $T$ similar to half of $ABC$.
Suppose $N/2$ is not a rational square.     Then $a$ belongs to $\Q(b)$,  and $b$ belongs to $\Q(a)$,
and both $a$ and $b$ belong to $\Q\big(\sqN\big)$. 
\end{lemma}
 
\noindent{\em Proof}.
Squaring both sides of (\ref{eq:pqrsqrtN}) and simplifying as a polynomial in $a$ we have 
\begin{eqnarray*}
(pa + db + e)^2 &=& \frac N 2 \\
p^2a^2 + 2pa(db+e) + (db+e)^2 - \frac N 2 &=& 0
\end{eqnarray*}
Since $a = \sin \alpha$ and $b = \cos\alpha$ we have $a^2 + b^2 -1 = 0$. Replacing $a^2$ by $1-b^2$ we find
$$ 2ap(db+e) + (db+e)^2 + p^2(1-b^2) -\frac N 2= 0$$
which shows that $a$ belongs to $\Q(b)$.  Similarly we find that $b$ belongs to $\Q(a)$, so 
$\Q(a) = \Q(b) = \Q(a,b)$. 
Starting again from (\ref{eq:pqrsqrtN}) we have
\begin{eqnarray*}
db &=& \sqN  - pa - e 
\end{eqnarray*}
For typographical simplicity we write 
$$\lambda := \sqN.$$ 
This equation shows that if $d=0$, $a$ belongs to $\Q(\lambda))$, and since $\Q(b) = \Q(a)$, the 
conclusion of the theorem holds.   Therefore we may assume without loss of generality that $d \neq 0$.
Continuing, we have
\begin{eqnarray*}
d^2 b^2 &=& (\lambda - pa - e)^2 \\
d^2(1-a^2) &=& (\lambda - pa - e)^2  \\
&=& p^2a^2 - 2ap(\lambda-e) + (\lambda-e)^2
\end{eqnarray*}
Writing it as a polynomial in $a$ we have
\begin{equation} \label{eq:galois}
0 = a^2(p^2+d^2) - 2ap (\lambda-e)  + (\lambda-e )^2 -d^2,  
\end{equation}
which shows that $a$ has degree 2 over $\Q (\lambda)$ or is in $\Q(\lambda)$.
Since $d \neq 0$,  the quadratic term in (\ref{eq:galois})
does not vanish.  

Solving (\ref{eq:galois}) by the quadratic formula we have
\begin{eqnarray}
  a &=& \frac {p (\lambda-e)}{p^2+d^2} \pm \frac {\sqrt{ (\lambda-e)^2p^2+(p^2+d^2)(d^2- (\lambda-e )^2)}}{p^2+d^2}\nonumber \\
  && \frac {p (\lambda-e )}{p^2+d^2} \pm \frac {\sqrt{ p^2d^2 + d^4 -d^2  (\lambda-e )^2)}}{p^2+d^2} \nonumber \\
  &=& \frac {p (\lambda-e )}{p^2+d^2} \pm \frac {d\sqrt{ p^2+ d^2 -  (\lambda-e )^2)}}{p^2+d^2} 
  \label{eq:quadraticsolutionfora}
\end{eqnarray}
Define 
$$ \mu := \sqrt{ p^2+ d^2 - (\lambda-e )^2)}$$
For proof by contradiction, assume that $a$ (and $\mu$) do not belong to $\Q(\lambda)$.
Then $1$, $\lambda$, $\mu$, and $\lambda\mu$ constitute a basis for $\Q(\lambda)$ over $\Q$,
as shown by the equation for $a$ above.  Let $\sigma$ be the automorphism of $\Q(\sqN)$ that 
takes $\lambda$ to $-\lambda$. 
We extend $\sigma$   to be defined on $\Q\big(a,\sqN\big)$ and fix $\mu$.  Therefore $\mu = \mu\sigma$.
We have
\begin{equation}
\mu \sigma = \sqrt{ p^2+ d^2 -(-\lambda-e)^2)} \label{eq:musigma}
\end{equation}
Then
\begin{eqnarray*}
\mu &=& \mu\sigma  \\
\mu^2 &=& (\mu\sigma)^2 \\
 p^2+ d^2 -  (\lambda-e )^2)  &=&  p^2+ d^2 -  (-\lambda-e )^2)
\end{eqnarray*}
Subtracting the right hand side  from both sides we have 
\begin{eqnarray}
0 &=&  p^2+ d^2 -  (\lambda-e)^2) -(p^2+ d^2 -  (-\lambda-e )^2))  \nonumber  \\
0 &=& -  (\lambda-e )^2) +  (-\lambda -e )^2))\nonumber \\
 0 &=&  4e\sqN  \nonumber \\
e &=& 0    
\end{eqnarray} 
But with $e=0$ the expressions for $\mu$ and $\mu \sigma$ simplifies considerably:
\begin{eqnarray*}
\mu &=&  d \sqrt{p^2 + d^2 - \lambda^2}\\
&=& d \sqrt{p^2 + d^2 - \frac N 2}  \mbox{\qquad since $\lambda^2 = N/2$} \\
 \mu\sigma =&=& \sqrt{p^2 + d^2 + \frac N 2} \mbox{\qquad by (\ref{eq:musigma})}
\end{eqnarray*}
Since $d \neq 0$ we have $\mu \sigma > \mu$.  But this contradicts $\mu= \mu\sigma$.
This contradiction completes the proof by contradiction that $a$ belongs to $\Q(\lambda)$,
and hence the proof of the lemma.
  
\begin{lemma}\label{lemma:isoscelessplitting}
 Let $ABC$ be an isosceles triangle with base $BC$, tiled by triangle $T$ similar to half of $ABC$, in which 
angle $\alpha$ is not a rational multiple of $\pi$. 
Then the base angles do not split, and the vertex angle $A$ is shared by exactly two tiles with the same angle at $A$,
 i.e.the vertex angle at $A$ splits into 
exactly two equal angles.
\end{lemma}

\noindent{\em Proof}.
   We have
$\alpha + \beta = \pi/2$.  We also have an equation arising from the vertex splitting.  Let 
$P$, $Q$, and $R$ be the total number of $\alpha$, $\beta$,
and $\gamma$ angles of tiles at the vertices of $ABC$.  Then because the angles at these vertices must add to $\pi$,
we have $P\alpha + Q\beta + R\gamma = \pi$, and since $\gamma = \pi/2$ this becomes
$$P\alpha + Q \beta = \pi\big(1-\frac R 2\big).$$
If $Q \neq P$ then we can put $\beta = \pi/2 - \alpha $ and solve for $\alpha$:
$$ \frac \alpha \pi = \frac 1 2  \frac {Q + R - 2} {Q - P}$$
making $\alpha$ a rational multiple of $\pi$.  But by hypothesis, $\alpha$ is not a rational multiple of $\pi$.
 Therefore we can assume $P=Q$.  We claim $P \ge 2$.   Certainly $P$ cannot be 1, since if   one angle 
 $\alpha$ appears at a vertex of $ABC$, the rest of the angle cannot be filled without using another $\alpha$.
 If $P=0$ then $Q\beta = \pi(1-R/2)$, so $\beta$ is a rational multiple of $\pi$, and hence $\alpha = \pi/2-\beta$
 is also a rational multiple of $\pi$, contradicting the hypothesis.  Hence $P \ge 2$.
  Therefore $P \alpha + Q \beta \ge 2 \alpha + 2 \beta = \pi$; but
$P \alpha + Q \beta = \pi(1-\frac R 2)$, which is strictly less than $\pi$ unles $R=0$.
  Hence $R=0$ and $P=Q=2$.
 
If neither base angle splits, then the base angles are $\beta$ (by convention--as mentioned, in this section
we do not assume $\alpha \le \beta$), and the vertex angle splits into two $\alpha$ angles.  
  The only other possibility is that one of the base angles of $ABC$  
 splits into two $\alpha$ angles.  In that case the vertex angle is $\beta$ and $ABC$ is 
equilateral.  Hence $\alpha = \pi/6$, which is a rational multiple of $\pi$, contrary to hypothesis.   
  Therefore   it is the vertex angle of $ABC$ that splits.  That completes the proof of the lemma.
  \medskip

Fix any vertex $V$ of the tiling, and let $n$, $m$, and $\ell$ count the number of $\alpha$, 
$\beta$, and $\gamma$ angles at $V$, and let $k\pi$ be the angle sum at $V$, so $k=1$ at a non-strict vertex or a boundary vertex,
and $k=2$ at a strict interior vertex.  Then we have 
$$n \alpha + m \beta = \pi\big(k-\frac \ell 2\big)$$
and if $n \neq m$, we can solve this equation together with $\alpha + \beta = \pi/2$, obtaining
$$\frac \alpha \pi = \frac 1 2  \frac {m + \ell - 2k} {m-n}$$
Since $\alpha/\pi$ is not rational, we must have $n=m$ at each  vertex.  This means that at each 
boundary or non-strict vertex, there are three possibilities:  one each of $\alpha$, $\beta$, and $\gamma = \pi/2$,
or two right angles, or two each of $\alpha$ and $\beta$.  

\begin{lemma} \label{lemma:rationalcase}
Suppose $ABC$ is isosceles and tiled by triangle $T$ similar to half of $ABC$, and 
assume $\alpha$ is a rational multiple of $\pi$.  Then $N$ is even and either
\smallskip

(i) $N/2$ is a square, or
\smallskip

(ii) $N/2$ is a twice a square (that is, $N$ is a square) and $\alpha = \pi/4$, or 
\smallskip

(iii) $N/2$ is three times a square and $\alpha = \pi/6$.
\end{lemma}

\noindent{\em Proof}.
Suppose that $\alpha$ is a rational multiple of $\pi$.  By Lemma \ref{lemma:degree2}, $e^{i\alpha}$ has degree 2 or 4 over $\Q$.
We can therefore apply
 Lemma \ref{lemma:euler} to conclude that $\alpha = 2\pi/n$, where $n = 5$, 8, 10, or 12.

In case $n=8$ we have $\alpha = \pi/4$; hence the left hand side of (\ref{eq:pqrsqrtN}) belongs to $\Q(\sqrt 2)$;  hence $\sqrt{N/2}$ belongs 
to $\Q(\sqrt 2)$.  Then $\sqrt{N/2}$ has the form $u + v \sqrt 2$ with $u$ and $v$ rational. Squaring both sides 
we have $N/2 = u^2 + 2v^2 + 2uv \sqrt 2$. Hence $uv=0$.  In case $v=0$ then $N/2$ is a square. In case $u=0$ then 
$N/2$ is twice a square.

In case $n=12$, $\alpha = \pi/6$, so $\cos \alpha = \sqrt{3}/2$ and $\sin \alpha = 1/2$; hence the left hand side belongs to $\Q(\sqrt 3)$; hence
$\sqrt{N/2}$ belongs to $\Q(\sqrt 3)$.  Then $\sqrt{N/2}$ has the form $u + v \sqrt 3$ with $u$ and $v$ rational. Squaring both sides 
we have $N/2 = u^2 + 3v^2 + 2uv \sqrt 3$.    Hence $uv = 0$.  Hence either $u=0$ or $v=0$. In case $u=0$ then 
$N/2$ is three times a square (which is possible, for example by bisecting each tile in Fig.~\ref{figure:prime27}, producing a 54-tiling); in case $v=0$ then $N/2$ is a square.

In case $n=10$ we have $\alpha = \pi/5$.  Then  $\cos \alpha = (1/4)( 1 + \sqrt 5)$, and 
$$\sin \alpha =  \frac 1 2 \sqrt{ \frac 1 2( 5 - \sqrt 5)}$$
But by Lemma \ref{lemma:degree2}, $\sin \alpha$ must belong to $\Q(\cos \alpha)$;  hence $\sqrt{(5 - \sqrt 5)/2}$  belongs to $\Q(\cos \alpha) = \Q(\sqrt 5)$.
That is, for some rational numbers $u$ and $v$, we have $\sqrt{(5-\sqrt 5)/2} = u + v \sqrt 5$.  A bit of algebra, not reproduced here,
shows that this is impossible, so the case $n=10$ cannot actually arise.

In case $n=5$ we have $\alpha = 2\pi /5$ and 
\begin{eqnarray*}
\sin \alpha &=&  \frac 1 2 \sqrt{ \frac 1 2( 5 + \sqrt 5)} \\
\cos \alpha &=& \frac 1 4 (-1 + \sqrt 5)
\end{eqnarray*}
and in this case also $\sin \alpha$ does not belong to $\Q(\cos \alpha)$, so by Lemma \ref{lemma:degree2}, this case cannot actually arise.
That completes the proof of the lemma.

\begin{lemma} \label{lemma:isosceles-irrational} Suppose the isosceles triangle $ABC$ is $N$-tiled by a right triangle  similar to half of $ABC$.
 Suppose that $\alpha$ is not a rational multiple of $\pi$, and $N/2$ is not an integer square or a sum of 
 two integer squares.  Then 
 \smallskip
 
 \noindent
 (i) in the $\dd$ matrix we have $d \neq 0$ and $p \neq 0$ and $e=0$, i.e. there 
 are no $c$ edges on the two equal sides $AB$ and $AC$, and there are some $a$ and some $b$ edges there, 
 and 
 \smallskip
 
 \noindent
 (ii) the second row of the $\dd$ matrix is identical to the first, and 
 \smallskip
 
 \noindent
 (iii) there are no ``edge relations'', i.e. no
 relations $ua + vb + wc = 0$ with rational $u$, $v$, and $w$.   
\end{lemma}

\noindent{\em Proof}. First, we prove that a half-integer $N/2$ is a square of an integer if and only if it is 
a rational square.  Suppose $N/2 = (P/Q)^2$, with $P$ and $Q$ relatively prime. 
  Then $N Q^2 = 2P^2$, so $N$ is even, since 2 divides the right side to an odd power, and 2 divides $Q^2$ to 
  an even power.  Then $N/2$ is an integer and $(N/2) Q^2 = P^2$.  Since $P$ and $Q$
  are relatively prime, we must have $Q=1$ and $N/2 = P^2$.
  
  Next, we note that by Lemma~\ref{lemma:half-integer} and Lemma~\ref{lemma:tworationalsquares},
   a half-integer $N/2$ is a sum of two integer squares if and only if it is a sum of two 
  rational squares.   Hence it does not matter whether the hypothesis of the lemma mentions sums of 
  two rational squares, or sums of two integer squares.  

 For notational simplicity, in this proof we continue to use
$$\lambda := \sqN.$$ 
Since half of $ABC$ is similar to the tile $T$, with scale factor $\lambda$ we have 
$$ \dd\vector a b 1 = \lambda  \vector 1 1 {2a}$$
The first row of this equation says that 
$$ \lambda = pa + db + 1,$$
from which it follows that if $a$ and $b$ are both rational, $\lambda^2 = N/2$ is a rational square,
and hence an integer square.  Therefore 
 not both $a$ and $b$ are rational. 

Suppose $a$ is rational. Then by Lemma~\ref{lemma:degree2}, since $N/2$ is not a square,  
 $b$ belongs to $\Q(a)$, so $b$ is 
also rational, contradiction, since we have shown   not both $a$ and $b$ are rational.
Therefore $a$ is not rational.  Similarly, if $b$ is rational, then $a$ belongs to $\Q(b)$,
so $a$ is also rational, contradiction.  Hence $b$ is not rational.

Recall that the $\dd$ matrix is 
$$ \dd= \matrix p d e g m f h \ell r.$$ 
We first consider the case $d=0$.  
Assume $d=0$.   Since $pa + db + e = \lambda$,  with $d=0$ we
have $pa + e = \sqN$.  If $p=0$ then $\lambda = e$, so $N/2$ is a square, and we are finished.
Hence we may assume $p \neq 0$.  If $e=0$ then $ap = \lambda$, so $N/2$ is a square, and
we are finished. 
  
By Lemma~\ref{lemma:degree2}, $b$ belongs to $\Q(\lambda)$, so for some integers $n$ and $k$ we have
\begin{eqnarray*}
b &=& \sqrt{1-a^2} \\
  &=& n + k \lambda \\
\sqrt{1-a^2}  &=& n + k \lambda
\end{eqnarray*}
Squaring both sides we have 
\begin{eqnarray*}
1-a^2 &=& n^2 + k^2 + 2kn\lambda
\end{eqnarray*}
 Since $pa + e = \lambda$, we have 
\begin{eqnarray*}
a &=& -\frac e p + \frac{\lambda}{p}
\end{eqnarray*}
Putting in this expression for $a$ we have
\begin{eqnarray*}
1-\bigg(-\frac e p + \frac{\lambda}{p}\bigg)^2 &=& n^2 + k^2 + 2kn\lambda\\
1 - \frac{e^2}{p^2}  - \frac N {2p^2} + \frac{2e}{p^2} \lambda&=& n^2 + k^2 + 2kn\lambda
\end{eqnarray*}
If $N/2$ is a square, we are finished, so we may assume $N/2$ is not a square, and 
then we can equate the rational parts:
\begin{eqnarray*}
1 - \frac{e^2}{p^2}  - \frac N {2p^2} &=& n^2 + k^2
\end{eqnarray*}
Since $b$ is not zero, and $b= n + k \lambda$, we cannot have both $n=0$ and $k=0$.  Hence the right hand side
is a positive integer,  so it is at least 1.  The left hand side, however, is less than 1,
since $N$ is positive (even if $e=0$).   This is a contradiction.   We have shown 
that if $d=0$, then $N/2$ is either a square or a sum of two squares; but that 
contradicts our assumptions.   Therefore $d \neq 0$.  

Our next aim is to prove $e=0$. 
Since $d \neq 0$, the expression  (\ref{eq:quadraticsolutionfora}) is valid.  That is, 
$$ a = \frac{ p(\lambda -e) \pm d\sqrt{p^2 + d^2 - (\lambda -e)^2} } {p^2 + d^2}. $$
Since $a$ belongs to $\Q(\lambda)$ (by Lemma~\ref{lemma:degree2}), the expression 
under the square root is a square in $\Q(\lambda)$.  It therefore suffices to prove that 
if $p^2 + d^2 + (\lambda -e)^2$ is a square in $\Q(\lambda)$, then $e=0$.

Let $\sigma$ be the automorphism of $\Q(\lambda)$ determined by $\lambda \sigma = -\lambda$.
Define
$$ \xi := \sqrt{p^2 +d^2 - (\lambda -e)^2} $$
Then 
\begin{eqnarray*}
\xi\sigma &=&\pm \sqrt{p^2 + d^2 - (\lambda + e)^2} \\
Norm(\xi) &=& \xi(\xi \sigma)  \\
\pm (Norm(\xi))^2 &=& (  p^2 +d^2 - (\lambda -e)^2)(p^2 +d^2 - (\lambda +e)^2)\\
&=& (p^2 + d^2)^2 + (\lambda - e)^2 (\lambda + e)^2 + (p^2 + d^2)((\lambda+e)^2 - (\lambda - e)^2) \\
&=&  (p^2 + d^2)^2  + (\lambda^2-e^2)^2  + 4e(p^2 + d^2 ) \lambda \\
&=&  (p^2 + d^2)^2  + \big(\frac N 2 -e^2\big)^2  + 4e(p^2 + d^2 ) \lambda  \mbox{\qquad since $\lambda^2 = N/2$} \\
\end{eqnarray*}
The left hand side of this equation is rational.  Hence, unless $N/2$ is a square (in which case
we are finished),  the coefficient of $\lambda$ on the right is zero.  Hence 
$e(p^2 + d^2) = 0$.  But  since $d \neq 0$, this implies $e=0$ as desired.  

Next we claim $p \neq 0$ and $r \neq 0$.
 Consider the tile, say Tile 1, at vertex $B$.   There is only one, since the angle 
there is $\beta$, and if it splits, then $\beta$ is a multiple of $\alpha$, and hence $\alpha$ is 
a rational multiple of $\pi$, contrary to hypothesis.  Since  Tile 1 has its $b$ side opposite its 
$\beta$ angle, its $a$ or $c$ side must be on $AB$.  But since $e=0$, it cannot be the $c$ side.  Hence 
the $a$ edge of Tile 1 is on $AB$.  But since $p$ is the number of $a$ edges on $AB$, we have $p \neq 0$.
Then the $c$ side of Tile 1 is on $BC$, and since  $r$ is the number of $c$ edges on $BC$, we have $r \neq 0$.
This completes the proof of claim (i) of the lemma.

 Now that we have proved $e=0$, (\ref{eq:quadraticsolutionfora})
 simplifies considerably:
 \begin{eqnarray}
a &=& \frac {p \lambda}{p^2 + d^2} \pm \frac {d \sqrt{p^2+d^2 - N/2}}{p^2+d^2} \label{eq:formulaforawhene=0}
\end{eqnarray}

We next derive a similar equation for $b$.  
\begin{eqnarray*}
pa + db  &=& \lambda  \mbox{\qquad (since $e=0$,  $e$ does not appear)} \\
pa &=& \lambda- db \\
p^2a^2 &=& (\lambda - db )^2 \\
p^2(1-b^2) &=&   (\lambda-db )^2 \\
&=& d^2b^2 - 2db\lambda+ \frac N 2   \\
0 &=& b^2(p^2+d^2) - 2bd\lambda+ \frac N 2 - p^2 
\end{eqnarray*}
Solving by the quadratic formula and simplifying,  we find the desired formula for $b$:
\begin{eqnarray}
b &=& \frac {d\lambda}{p^2+d^2} \pm \frac {p \sqrt{p^2+d^2 - N/2}}{p^2+d^2} \label{eq:formulaforbwhene=0}
\end{eqnarray}

We now define
\begin{equation} \label{eq:mudefn}
\mu := \sqrt{p^2 + d^2 - N/2}
\end{equation}
We will show that either $N/2$ is a sum of two squares, or $\mu$ is rational.  Suppose  
 that $\mu$ is not rational.  Since $a$ belongs to $\Q(\sqN)$, 
so does $\mu$. But the only rationals that have irrational square roots in $\Q(\sqN)$ are rational 
squares times $\sqN$.  Hence for some rational $\xi$ we have $\mu = \xi \sqN$.
Then
\begin{eqnarray*}
\mu^2 &=& \xi^2 \frac N 2 \\
p^2 + d^2 - \frac N 2 &=& \xi^2 \frac N 2 \\
p^2 + d^2 &=& (\xi^2 + 1) \frac N 2 \\
\frac N 2 &=& \frac {p^2 + d^2}{\xi^2 + 1} 
\end{eqnarray*}
This is a quotient of sums of two rational squares.  By Lemma~\ref{lemma:sumsofsquares}, it is 
a sum of two rational squares, and since it is a half-integer, it is a sum of two integer squares, 
by Lemma~\ref{lemma:half-integer}.  Hence, as claimed, either $N/2$ is a sum of squares,
or $\mu$ is rational.  But we have earlier assumed that $N/2$ is not a sum of two integer squares;
hence $\mu$ is rational.
 We also have 
$\mu \neq 0$, since if $\mu =0$ we have $N/2 = p^2 + d^2$.

In terms of $\lambda$ and $\mu$, our formulas for $a$ and $b$ are 
\begin{eqnarray*}
a &=& \frac {p \lambda \pm d\mu}{p^2 + d^2} \\
b &=& \frac {d \lambda \pm p\mu}{p^2 + d^2} 
\end{eqnarray*}
We now claim that one must take opposite signs for the $\pm$ in these two formulas.
To prove that, we calculate as follows:
\begin{eqnarray*}
1 &=& a^2 + b^2 \\
&=& \big( \frac{ p \lambda \mp d \mu} {p^2 + d^2} \big)^2 + \big( \frac{ d \lambda \pm p \mu} {p^2 + d^2} \big)^2  \\
&=& \frac {(p^2 + d^2)( \lambda^2 +\mu^2)}{(p^2 + d^2)^2}  + \mbox{\ \ cross terms} 
\end{eqnarray*}
Here the ``cross terms'' are either zero, if $a$ and $b$ have opposite signs on the coefficients of $\mu$, 
or they are $\pm 2dp\lambda \mu /(p^2 + d^2)$, which is not zero unless $\mu = 0$, since $p$ and $d$ 
are not zero.   Now $\lambda^2 = N/2$ and $\mu^2 = p^2 + d^2 - N/2$, so $\lambda^2 + \mu^2 = p^2 + d^2$
and we get 1 for $a^2 + b^2$  {\em without} the cross terms.  
Hence the cross terms {\em must} cancel out.  We will see below that one must take the plus sign for $a$
and the minus sign for $b$.    The fact that the signs must be opposite can also be proved by 
starting from $pa + db = \lambda$ instead of from $a^2 + b^2 = 1$. 
 
  Now we turn to the third row of the $\dd$ matrix, which tells us 
$$ ha + \ell b + r = 2\lambda a.$$
We will use this to determine the signs in the formulas for $a$ and $b$, which so far have an 
ambiguous $\pm$ sign, except that we know the signs must be opposite.  Substituting in the formulas for 
$a$ and $b$ we have 
\begin{eqnarray*}
2 \lambda \big(\frac{p \lambda \pm d\mu} {p^2 + d^2}\big) &=& 
h \big(\frac{p \lambda \pm d\mu} {p^2 + d^2}\big) + \ell big(\frac{d \lambda \mp p\mu} {p^2 + d^2}\big) + r 
\end{eqnarray*}
Clearing denominators we have
\begin{eqnarray}
2 \lambda (p \lambda \pm d\mu) &=& h (p \lambda \pm d \mu) + \ell (d \lambda \mp p \mu) + r(p^2 + d^2) 
\label{eq:2869}
\end{eqnarray}
Equating the coefficients of $\lambda$,   and noting that $\mu$ and $\lambda^2$ are real, we have 
\begin{eqnarray*}
\pm 2d\mu &=&hp +\ell d 
\end{eqnarray*}
Since the right-hand side is nonnegative, and since $d$ and $\mu$ are not zero, we must take the positive 
sign on the left, which means we must take the positive sign in the formula for $a$, and hence the 
negative sign in the formula for $b$.  That is, 
\begin{eqnarray}
a &=& \frac {p\lambda + d\mu}{p^2 + d^2} \label{eq:aformula}\\
b &=& \frac {d \lambda - p\mu}{p^2 + d^2} \label{eq:bformula}
\end{eqnarray}

We next claim that the second row of the $\dd$ matrix is identical to the first, i.e. the 
numbers of edges of each length are the same on $AB$ as on $AC$.  From the second row 
of the $\dd$ matrix equation we have $ga + mb + fc = \lambda$.  This has the same 
right hand side as the equation $pa + db + ec  = \lambda$ from the first row.  We can 
change $(p,d,e)$ to $(g,m,f)$ in all the above calculations, and we find instead of $e=0$
that $f=0$, and instead of $d \neq 0$, that $g\neq 0$, and finally we find formulas for $a$ 
and $b$:
\begin{eqnarray*}
\nu &:=& \sqrt{g^2 + m^2 - N/2} \\
a &=& \frac {g\lambda + m\nu}{g^2 + m^2} \\
b &=& \frac {m \lambda - g\nu}{g^2 + m^2}
\end{eqnarray*}
Equating the coefficients of $\lambda$ in the two formulas for $a$ we have
$$ \frac d {d^2 + p^2} = \frac g {g^2 + m^2}$$
This implies that the points $(p,d)$ and $(g,m)$ lie on the same line through the origin,
and hence for some real $t$, $m = td$ and $g =tp$.  Then let $x$ be the coefficient of $\lambda$ in $a$.
We have 
\begin{eqnarray*}
x &=&  \frac g {g^2 + m^2}\\
  &=& \frac {td} {(td)^2 + (tp)^2 } \\
  &=& \frac 1 t   \frac d { d^2 + p^2 } \\
  &=& \frac  x t 
\end{eqnarray*}
Thus $x = x/t$, which implies $t=1$.  Hence 
\begin{eqnarray}
g &=& p  \label{eq:g=p} \\
m &=& d \label{eq:m=d}
\end{eqnarray}
 That proves part  (ii) of the lemma.
 
We now turn our attention to part (iii) of the lemma. 
Suppose, for proof by contradiction, that $ua + vb + w = 0$ 
 with rational $u$, $v$ and $w$.   
 Then 
 \begin{eqnarray*}
 0 &=& ua + vb + w \\
 &=& u \frac{p\lambda + d\mu}{p^2 + d^2} + v \frac{d\lambda -p\mu}{p^2 + d^2} + w \\
 0 &=& u (p\lambda + d\mu) + v(d\lambda -p \mu) + w(p^2 + d^2) \\
 &=& \lambda(up + vd) + \mu(ud -vp) 
\end{eqnarray*}
Since $p > 0$ and $d > 0$,  we must have $u = v = 0$ to make the coefficient of $\lambda$ zero.
But then $0 = ua + vb + w = w$, so $w$ is also zero.   We note in passing the consequence that 
$a/b$ is irrational, since if $a/b = v$ than $a-bv = 0$.
That completes the proof of the lemma.

\begin{theorem} \label{theorem:isosceles-irrational} Suppose the isosceles triangle $ABC$ is $N$-tiled by a right triangle  similar to half of $ABC$.
 Suppose that $\alpha$ is not a rational multiple of $\pi$.
  Then  $N/2$ is a square (of an integer) or a sum of two integer squares.    
\end{theorem}

\noindent{\em Proof}.   Suppose, for proof by contradiction, that there is such a tiling and $N/2$ is 
neither a square nor a sum of two squares.  Then, by Lemma~\ref{lemma:isosceles-irrational}, there are 
no edge relations.   That means that each maximal segment in the tiling has equal numbers of $a$ edges
on each side,  equal numbers of $b$ edges on each side, and equal numbers of $c$ edges on each side, 
for otherwise an edge relation exists.

Consider two points $Q$ on $BA$ and $U$ on $BC$ such that $BQU$ is similar to the tile, 
with a right angle at $Q$ and an $\alpha$ angle at $U$.  We say that $BQU$ is ``nicely tiled''
if the tiling of $ABC$, restricted to $BQU$,  is obtained from a quadratic tiling of $BQU$ 
by replacing some pairs of adjacent tiles that form a rectangle by the same rectangle with 
the other diagonal.  
Choose
segment $QU$   as far to the northeast as possible with $BQU$ nicely tiled.  There is 
some such segment $QU$ because the tile at $B$ has its $c$ side on $BC$, not on $AB$, since 
$e=0$ by Lemma~\ref{lemma:isosceles-irrational};  and since $d \neq 0$ (also by Lemma~\ref{lemma:isosceles-irrational}), there are some $a$ edges on $AB$, so $Q$ lies on the interior of $AB$.
Point $U$ lies on the interior of $BC$, since the tile at $C$ has its $\beta$ angle at $C$, but 
such a tile will not be part of the lattice tiling of $BQU$, since those tiles have their 
$\alpha$ angle to the east.
Let $P$ be  the next vertex on $QA$ above $Q$, and let $W$ be the next vertex on $BC$ east of $U$. 
Then all the tiles below $QU$ with an edge on $QU$ have their $c$ edge on $QU$.  Since there 
are no linear integral relations between $a$, $b$, and $c$, all the tiles above $QU$ with an edge
on $QU$ also have their $c$ edges on $QU$.   Suppose  that $PQ$ has length $a$.

Let Tile 1 be the tile with a vertex at $Q$ and its southwest edge on $QU$; let $R$ be its
southeast vertex on $QU$.  Then there are two cases:  the third vertex of Tile 1 is either $P$
or the point $S$ such that $SR$ is perpendicular to $QR$ and $SR$ has length $a$.   In case the 
third vertex of Tile 1 is $S$, then Tile 1 has its $\alpha$ angle at $Q$, and hence the 
remaining angle $PQR$ is $\beta$, which must be filled by a single tile, Tile 2.  Since $PQ$
has length $a$, Tile 2 must have its $c$ edge along $QS$ and hence the rectangle $PSRQ$ is 
nicely tiled.

\begin{figure}[ht]
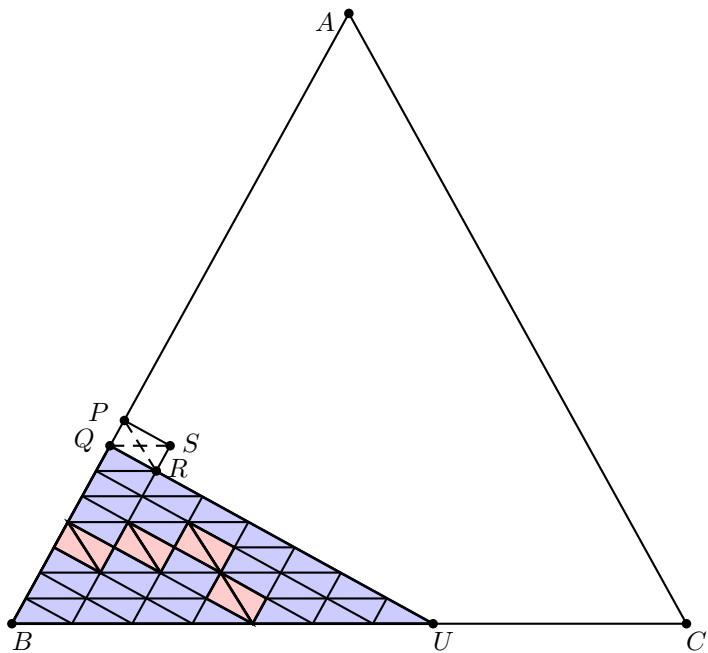

\caption{$PQRS$ must be nicely tiled, assuming $\overline{PQ} = a$.}
\label{figure:3024}
\FigureThreeZeroTwoFour
\end{figure}

In case  the third vertex of Tile 1 is $P$, then $PR$ has length $c$, and $PR$ is a maximal segment, so the 
tile northeast of $PR$, say Tile 2, shares its $c$ edge with $PR$.  Suppose, for proof 
by contradiction, that Tile 2 is not $PRS$.  Then it has its $\beta$ angle at $P$,
and its $\alpha$ angle at $R$.  Tile 1 also has its $\alpha$ angle at $R$, so the 
angle between Tile 2 and $RU$ is $2\beta$.  That is partly filled by the angle at 
$R$ of Tile 3, the tile above $RU$ with a $b$ edge on $RU$.  Tile 3 has its $b$ edge
on $RU$, so its angle at $R$ is either $\alpha$ or $\gamma$; but since it must fit into $2\beta$,
it cannot be $\gamma$.  It must therefore be $\alpha$.  But that  leaves
an unfilled angle of $\beta - \alpha$, which cannot be filled by any number of $\alpha$ angles,
since $\alpha$ is not a rational multiple of $\pi$.  That contradiction shows that Tile 2 
is in fact $PRS$.  Hence in this case also, rectangle $PSRQ$ is nicely tiled.   See 
Fig.~\ref{figure:3024}.

  Now let $W$ be the intersection point of $UC$ and 
the line containing $PS$.   Let the vertices on $QU$ be 
$Q, R, R_2, R_3, \ldots$, spaced apart by $b$, and let the points $P, S_1, S_2, \ldots$ on $PW$
also be spaced apart by $b$, so $S=S_1$.  Let $k$ be the largest integer such that $PS_kR_kQ$
is nicely tiled.  If $R_k$ lies on $BC$, then (since $QW$ was chosen as far to the northeast as possible),
the tile with vertex at $R_k$ and an edge on $R_kC$ has its $\beta$ angle at $R_k$, and its $c$ 
edge against the $a$ edge $R_kS_k$.   

Otherwise point $R_{k+1}$ exists.  Let Tile 5 be the tile above $QU$ with its $b$ edge on $Q_kU$ and a vertex at $Q_k$.
The third vertex of Tile 5 is either at $S_k$ or at $S_{k+1}$.  First consider the case
when the third vertex of Tile 5 is at $S_{k+1}$.  Then there is exactly one more tile, Tile 6, with 
a vertex at $R_k$, and Tile 6 has its $\beta$ angle at $R_k$.  If Tile 6 has its $a$ side on 
$R_kS_k$, then the rectangle $PQR_{k+1}S_{k+1}$ is nicely tiled, contradicting the choice of $k$. 
Hence Tile 6 has its $c$ side along $R_kS_k$.  Then there is a maximal segment through $R_kS_k$
with only $a$ edges on the northwest below $S_k$ and at least one $c$ edge on the southeast.  
See Fig.~\ref{figure:3048}.

\begin{figure}[ht]
\caption{Tile 5 has a vertex at $S_{k+1}$, so along $R_kS_k$ there are an $a$ edge and a $c$ edge.}
\label{figure:3048}
\FigureThreeZeroFourEight
\end{figure}

Next consider the case when the third vertex of Tile 5 is at $S_{k}$.  Let Tile 6 be 
the tile with a vertex at $R_{k+1}$ and sharing an edge with Tile 5.  We claim that Tile 6
has its $\beta$ angle at $R_{k+1}$.  If $R_{k+1}$ lies on $BC$ (so $R_{k+1}  = U$), 
then angle $S_kUC$ is $2\beta$,  so it can only be filled by two tiles with their $\beta$ 
angles at $R_k$, one of which is Tile 6.   On the other hand, if $R_{k+1}$ does not lie on $BC$, 
then $R_{k+2}$ exists, and Tile 7 above $R_{k+1}R_{k+2}$ has its $b$ edge on $R_{k+1}R_{k+2}$.
Then if Tile 7 has its right angle at $R_{k+1}$, that leaves a $\beta$ angle for Tile 6, 
and if Tile 7 has its $\alpha$ angle at $R_{k+1}$, that leaves $2\beta$ unfilled at $R_{k+1}$,
which can only be filled by two tiles with their $\beta$ angles at $R_{k+1}$, one of which is 
Tile 6.  Hence in any case, Tile 6 has its $\beta$ angle at $R_{k+1}$.  It does not have its
$c$ edge on $R_{k+1}S_k$, as that would make $PQR_{k+1}S_{k+1}$ nicely tiled, contradicting
the choice of $k$.  Hence Tile 6 has its $a$ edge on $R_{k+1}S_k$.   Then there is a maximal
segment through $R_{k+1}S_k$,  with only $c$ edges on the southwest side below $S_k$, and 
at least one $a$ edge on the northeast side.  See Fig.~\ref{figure:3068}.

\begin{figure}[ht]
\caption{Tile 5 has a vertex at $S_k$, so along $S_kR_{k+1}$ there are an $a$ edge and a $c$ edge.}
\label{figure:3068}
\FigureThreeZeroSixEight
\end{figure}

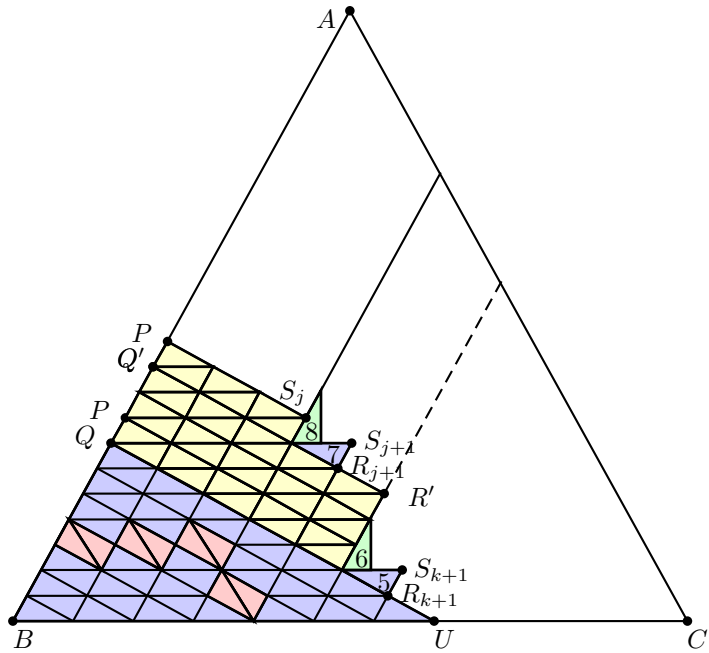
\begin{figure}[ht]
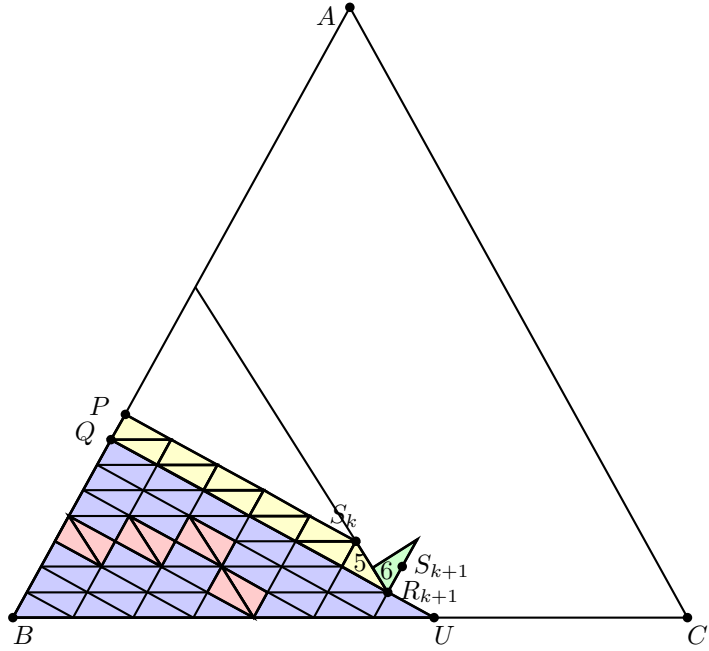

\caption{The second stage of the construction; the strip under $Q^\prime R^\prime$ is nicely tiled.}
\label{figure:3578}
\FigureThreeFiveSevenEight
\end{figure}

   Thus we are in the same situation, whether $R_k$ lies on $BC$ or not, and regardless of the orientation 
of Tile 5:  to the southeast of a nicely-tiled strip $PQS_kR_k$
there is a line, either containing $R_kS_k$ or $R_{k+1}S_k$,  with a $c$ edge on one side and 
an $a$ edge on the other.  Let $F$ be the point $R_k$ or $R_{k+1}$ where this line starts.
Let $T$ be the first point above $S_k$ on this line, such that $T$ is 
a vertex both of a tile on the left and a tile on the right of $FT$.  Then there 
must be equal numbers of $a$ edges on each side of $FT$, and equal numbers of $c$ edges on 
each side, because otherwise there would be a nontrivial integer linear relation between $a$, $b$, and $c$.
  Let $Q^\prime $ and $R^\prime$ be 
points on $QA$ and $FT$ respectively such that $Q^\prime R^\prime$ is parallel to $QR$ and 
rectangle $Q^\prime QR_kR^\prime$ is nicely tiled, and $Q^\prime R^\prime$ is as far northeast as 
possible.  There must be such points $Q^\prime$ and $R^\prime$, with $R^\prime$ on the interior of 
segment $FT$, since otherwise $FT$ has all $a$ edges on its left side (if $FT$ is parallel to $AB$),
or $FT$ has all $c$ edges on its left side (otherwise),  either of which 
contradicts Lemma~\ref{lemma:isosceles-irrational}. 
 Then there are $b$ edges all along $Q^\prime R^\prime$.   Now we repeat the argument that
we made above, with $QU$ replaced by $Q^\prime R^\prime$,  and we find another line $R_j S_j$ with $j < k$
such that there is one row of nicely tiled rectangles above $Q^\prime R^\prime$ between $AB$ and line 
$R_jS_j$.  Let $R^{\prime\prime}$ be the intersection point of $Q^\prime R^\prime$ and $R_jS_j$.
Then there is a tile, say Tile 7,  with its $b$ edge on $Q^\prime R^\prime$ southwest of $R^{\prime\prime}$.
As before there are two possibilities for the orientation of Tile 7.  One of them is shown in 
Fig.~\ref{figure:3578}, namely when Tile 7 has its $\alpha$ angle at $R_j$.  Then Tile 8 has its 
$\beta$ angle at $R_j$ and its $c$ edge against the $a$ edge of the tile northwest of $R_jS_j$, 
so $R_jS_j$ must extend northeast until the numbers of $a$ edges on the left and right are equal and 
the number of $c$ edges on the left and right are equal and the number of $b$ edges on the left and right are equal.
If Tile 7 has the other orientation, then the same is true of the line $R_{j+1}S_j$.   

\FloatBarrier

\begin{figure}[ht]
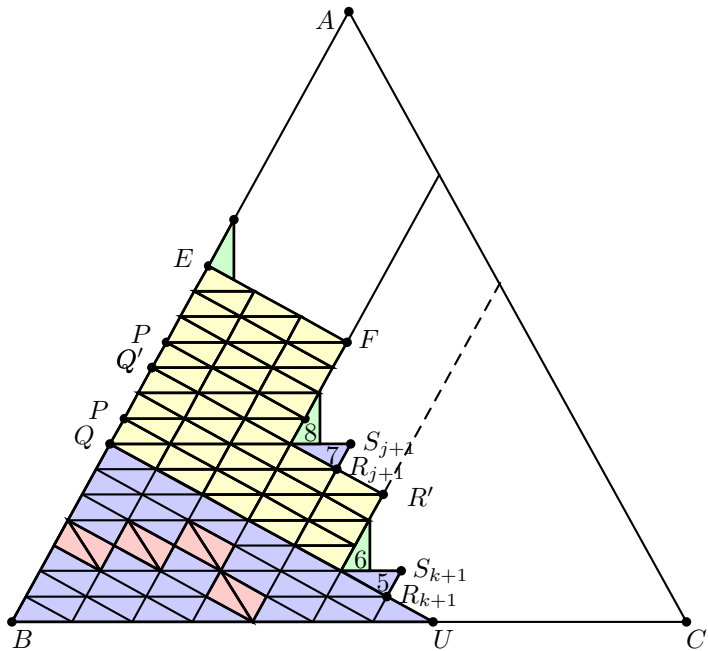

\caption{The final contradiction when we encounter a $b$ edge above $E$ on $AB$.}
\label{figure:4000}
\FigureFourThousand
\end{figure}

 We continue in this fashion, defining narrower strips of nicely-tiled rectangles bounded
on the east by lines with $a$ edges on the left and $c$ edges on the right.  Eventually one of two 
things happens.  Either
we run out of $a$ edges on $AB$, or we reach a point where 
one of the boundary lines $R_{j+1}S_j$ that run towards the northwest intersects $AB$.  For example,
in Fig.~\ref{figure:3068},  the nicely-tiled area could fill up the region under the line $R_{j+1}S_j$ 
before we run out of $a$ edges on $BC$.   We will show that both possibilities are impossible.  The 
second one is immediately contradictory:  if $R_{j+1}S_j$ reaches $AB$ at a point $G$ 
with all the area underneath it 
nicely tiled, then there are only $c$ edges underneath $R_jG$,  but above it there is at least one $a$ 
edge,  giving rise to an edge relation, which contradicts Lemma~\ref{lemma:isosceles-irrational}.
Therefore this possibility does not happen.  We cannot reach the vertex $A$ without running out of $a$ 
edges on $AB$, since 
 by Lemma~\ref{lemma:isosceles-irrational}, $d \neq 0$, so there is at least one $b$ edge on $AB$.
   Let $E$ be the highest vertex on $AB$ such that below $E$,
there are only $a$ edges on $AB$.  Then let Tile 9 be the tile with an edge on $AB$ north of $E$
and a vertex at $E$.  Since Tile 9 does not have an $a$ edge on $AB$, and by Lemma~\ref{lemma:isosceles-irrational}, there are no $c$ edges
on $AB$, Tile 9 must have its $b$ edge on $AB$.   Let $EF$ be the line perpendicular to $AB$ at $E$,
where $EF$ is a maximal segment.  Then $F$ lies on the southwest boundary $L$ of the last strip of 
nicely-tiled rectangles.  That boundary $L$ does not terminate at $F$ since there are only $a$ edges
on it south  of $F$. ($L$ may run southwest to northeast, or southeast to northwest; as we have seen, 
at each stage of the construction there are two cases, so we just say ``south'' here.) 
 Thus  every tile southwest of $EF$ with an edge on $EF$ has its $b$ edge on $EF$,
and there is a tile northeast of $EF$ with a vertex at $F$.  Hence all the tiles north of $EF$ with 
an edge on $EF$ have their $b$ edges on $EF$.
 If Tile 9 has its right angle at $E$ then it has 
its $a$ edge on  $EF$, contradiction.  This situation is shown in Fig.~\ref{figure:4000}.

 Hence Tile 9 has its $\alpha$ angle at $E$, leaving an angle of 
$\beta$ unfilled at $E$.   This $\beta$ angle must be filled by Tile 10,  which then cannot have its 
$b$ edge on $EF$.  That is a contradiction (not illustrated in a separate figure, but imagine the last
tile at $E$ in Fig.~\ref{figure:4000} placed with its $\alpha$ angle at $E$ instead of its right angle as shown).
   That contradiction completes the proof of the theorem.
 \medskip

 \begin{theorem} \label{theorem:isosceles}
Suppose triangle $ABC$ is isosceles with $AB = AC$.  Let $M$ be the midpoint of base $BC$,  and suppose $ABC$ is $N$-tiled by the right 
triangle $T$, and suppose $T$ is similar to triangle $ABM$  (half of triangle $ABC$).  Then $N$ is even and 
one of the following cases holds:
\smallskip

(i) $N/2$ is a square. 
\smallskip

(ii) $N$ is a square, and $\alpha = \pi/4$.
\smallskip

(iii)  $ABC$ is equilateral and $N$ is six times a square (so $\alpha = \pi/6$).
\smallskip

(iv)  $ABC$ has base angles $\pi/6$ and $N$ is six times a square.
\smallskip

(v) $N/2$ is a sum of two squares, and $\alpha$ is not a rational multiple of $\pi$,  and $a$ and $b$ are rational.
\end{theorem}

\noindent{\em Remark}.  All the cases mentioned in the theorem do actually occur.
\medskip

 \noindent{\em Proof}.
By Lemma \ref{lemma:isoscelessplitting}, the vertex angles at $B$ and $C$ do not split.  As mentioned 
above, we make the convention that $\beta$ is the base angle;  then by the lemma, $2\alpha$ is the vertex angle.
For this section only, we do not assume $\alpha \le \beta$. 
If $\alpha$ is a rational multiple of $\pi$, then by Lemma \ref{lemma:rationalcase}, one of the first four 
conclusions of the theorem holds.  We can therefore assume that $\alpha$ is not a rational multiple of $\pi$. 
Then  by Theorem~\ref{theorem:isosceles-irrational}, $N/2$ is a square or a sum
 of two squares, so conclusion (i) or conclusion (v) holds.  
That completes the proof of the theorem.

\section{Tilings by a non-isosceles right triangle $T$}
Many of the tilings exhibited in the introduction have the 
tile $T$ similar to the tiled triangle $ABC$.  Indeed that is the case for the quadratic tilings and the 
biquadratic tilings, but not the case for the $3$-tiling of the equilateral triangle, and various composite
tilings involving that tiling as a subtiling (such as the exhibited 12-tilings, and a 6-tiling obtained from the 
equilateral 3-tiling).  There is also, of course, the 2-tiling of an equilateral triangle, which can in turn 
be used with a quadratic tiling of an equilateral triangle to produce various 8-tilings. The 27-tiling
shown in Fig. 10 is an example of a prime tiling with $N>3$ in which the tile is not similar to $ABC$.
    In this section, we deal with the special case in which the tile is a right triangle.
We begin with a special case of this special case.

\begin{lemma} \label{lemma:piover8not} Suppose $T$ is a right triangle with $\alpha = \pi/8$, 
and suppose $ABC$ is a right triangle and $ABC$ is $N$-tiled by $T$ for some $N$. Then  $ABC$ is similar to $T$.
\end{lemma}

\noindent{\em Proof}. Let $\zeta =  e^{i \pi/8}$. Then the degree of $\Q(\zeta)$ over $\Q$ is $\varphi(16) = 8$, by Lemma~\ref{lemma:euler}.
Let $\sigma = \sigma_9$ be the automorphism taking $\zeta$ to $\zeta^9 = -\zeta$.  Then $\sigma$ fixes $i$, 
fixes  $\sin J\alpha $ for $J$ even and changes the sign of $\sin J\alpha$ for $J$ odd, since 
\begin{eqnarray*}
(2i \sin J \alpha) \sigma &=&  (\zeta^J - \zeta^{-J})\sigma \\
&=& (\zeta \sigma)^J - (\zeta \sigma)^{-J} \\
&=& (-\zeta)^J - (-\zeta)^{-J} \\
&=& (-1)^J 2i \sin J \alpha
\end{eqnarray*}
Hence $\sigma$ changes the signs of $a = \sin \alpha$ and $b = \sin \beta = \sin 3\alpha$.  Since $T$ is 
a right angle we have $c = 1$. 
Let $X = pa + qb + r$ and $Y = ma + nb + \ell$, where $p$, $q$, $r$, $m$, $n$, and $\ell$ give the numbers
of tiles with sides $a$, $b$, and $c$ along $X$ and $Y$.   

The area of the tile $\A_T$ is given by $ab/2$, since $\gamma$ is a right angle.
 Since triangle $ABC$ is also a right triangle its area $\A_{ABC}$ is given by $XY/2$.  Since 
 there is a tiling, we have the area equation $N\A_T = \A_{ABC}$,  which now becomes
 $Nab = XY$.  Writing out $X$ and $Y$ in terms of the edges that compose them, we have 
 $$ N ab = (pa + db + e)(ga + mb+ f) $$
Applying $\sigma$ we have 
$$ N ab = (-pa -db + e)(-ga -mb + f)$$
Subtracting the two equations  we have
\begin{eqnarray*}
0 &=& (pa + db + e)(ga + mb + f) - (-pa -db + e)(-ga-mb+f) \\
  &=&  2e(ga+mb) + 2f(pa+db)  \mbox{\qquad since the other terms cancel} 
\end{eqnarray*}
Dividing by 2 we have 
\begin{eqnarray*}
0 &=& gea + emb + fpa + fdb
\end{eqnarray*}
Since all these quantities are nonnegative, and $a$ and $b$ are not zero,
 we have $ge$, $em$, $fp$ and $fd $ each equal to zero.
Assume, for proof by contradiction, that $f \neq 0$. Then $p=d=0$.  Hence $e \neq 0$, 
since $X = pa + db + ec$ and $X \neq 0$.   Hence $m=g=0$.
Hence both $X$ and $Y$ are composed entirely of $c$ sides, and hence are rational, since $c=1$.
Hence $ab = XY/N$ is 
also rational.  But it can be expressed as 
\begin{eqnarray*}
ab &=& \sin \alpha \sin 3\alpha \\
&=& (\zeta - \zeta^{-1})(\zeta^3-\zeta^{-3}) \\
&=& \zeta^4 - \zeta^2 - \zeta^{-2} + \zeta^{-4} 
\end{eqnarray*}
Multiplying by $\zeta^4$ we have 
\begin{eqnarray*}
\zeta^4 ab &=& \zeta z^8 - \zeta^6 + \zeta^2 + 1 \\
&=& -\zeta^6 + \zeta^2  \mbox{\qquad since $\zeta^8 = -1$} \\
0 &=& -\zeta^6 -ab\zeta^4 + \zeta^2
\end{eqnarray*}
Since $ab$ is rational, this is a polynomial with rational coefficients satisfied by $\zeta$,
and it is not identically zero; but its degree is less than the degree of the minimal polynomial
of $\zeta$, which we have seen is 8.  
That contradiction completes the proof of the lemma.
\medskip

\begin{lemma} \label{lemma:special2}
Suppose that $T$ is a right triangle with $\alpha = \pi/8$, and that $ABC$ has
angle $C = 5\pi/8$ and angle $A = \alpha$.  
Let $N$ be a positive integer.  Then there is no 
$N$-tiling of $ABC$ by $T$.  
\end{lemma}

\noindent{\em Proof}.  We use the same cyclotomic field and automorphism $\sigma$ 
as in the previous proof.  The area equation can be written as
$$ Nabc = XY \sin 5\alpha $$
since angle $ABC =5\alpha$. Then $\sigma$ changes the signs of $a$ and $b$,
since these are respectively the sines of $\alpha$ and  $3\alpha$. 
Since $5$ is odd, $\sigma$ changes the sign of $\sin 5\alpha$.  On the other 
hand $\sigma$ fixes $c = \sin 4\alpha$.  Hence, when we apply $\sigma$ to 
the area equation, we get
$$ Nabc = -(X\sigma)(Y\sigma) \sin 5\alpha.$$
Now it will not help to subtract the equations, but we can divide them instead.  
We obtain
\begin{eqnarray*}
-1 &=& \frac { (X \sigma)(Y\sigma)} {XY}
\end{eqnarray*}
Multiplying by the denominator and adding $XY$ to both sides, we have
\begin{eqnarray*}
-XY &=& (X\sigma)(Y\sigma) \\
0 &=& (X \sigma)(Y \sigma) + XY
\end{eqnarray*}
Writing $X$ and $Y$ out in terms of the $\dd$ matrix, and remembering that $\sigma$
changes the signs of $a$ and $b$ but not $c$, we have 
\begin{eqnarray*}
0 &=& (-pa - db + e)(-ga-mb + f) + (pa +db + e)(ga+mb+f) \\
&=& 2(pa + db)^2 + 2(ef)^2   \mbox{\qquad since the other terms cancel out} \\
0 &=& (pa+db)^2 + (ef)^2
\end{eqnarray*}
Since all the letters denote nonnegative quantities, we have $ef = 0$ and 
$pa+db = 0$; since $a > 0$ and $b > 0$ we have $p=d=0$.   Since $X \neq 0$,
$p=d=0$ implies $e \neq 0$; then $ef=0$ implies 
\begin{equation}
f=0. \label{eq:4710}
\end{equation}
The area equation can also be written as 
$$ Nabc = YZ \sin \alpha$$
since angle $A$ of triangle $ABC$ is $\alpha$.  Since $\sin \alpha = a$ we have 
$$ Nbc = YZ.$$
Applying $\sigma$ we have 
\begin{eqnarray*}
-Nbc  &=& Y\sigma Z\sigma \\
\end{eqnarray*}
Adding that to the previous equation we have 
\begin{eqnarray*}
0 &=& (Y \sigma)(Z\sigma) + YZ \\
\end{eqnarray*}
Writing this out in terms of the $\dd$ matrix we have 
\begin{eqnarray*}
0 &=& (-ga-mb+fc)(-ha-\ell b + rc) + (ga+mb+fc)(ha+\ell b + rc) \\
&=& 2(ga+mb)^2 + 2(fr)^2 \\
0 &=& (ga+mb)^2 + (fr)^2
\end{eqnarray*}
Hence   $g=m = 0$.  But by (\ref{eq:4710}) we have $f=0$.
Hence $Y = ga + mb + fc  = 0$.  But that is a contradiction, 
since $Y$ is the length of side $AC$ of triangle $ABC$.
That completes the proof of the lemma.
\medskip

\noindent{\em Remark}.  The proofs given above of the previous two lemmas
are not the only proofs.  In both the case of a right angle and the case when 
angle $C = 5\alpha$,  the area of $ABC$ is given by $\A_{ABC} = X^2/2$.  That 
is obvious for the right triangle, and not very difficult in the other case.  
It is possible to show that this formula leads to a contradiction.  We sketch 
the idea rather than give the complete proof.   One writes $X$ in terms of the 
$\dd$ matrix as $X = pe + db + e$ (since $c=1$, $c$ does not appear here).
Since $b = \sin 3\alpha$, we have $b = 3a-4a^3$ by simple trigonometry.  
Substituting that in the expression for $X$ we have 
$$ X = pa + d(3a-4a^3) + e.$$
Since the area of the tile is $ab/2$,  the area equation 
$$ 2\A_{ABC} = 2N\A_T$$
becomes $X^2-Nab = 0$.  That is,
\begin{equation}
 (pa + d(3a-4a^3) + e)^2 -Na(3a-4a^3) = 0. \label{eq:4759}
\end{equation}
The left side is a polynomial in $a$ of degree 6.  One shows that this polynomial 
cannot be zero.  To do that one has to first show  
\begin{equation}
 a^4 - a^2 + \frac 1 8 = 0, \label{eq:4764}
 \end{equation}
   which follows from the formula 
$$ \sin\big( \frac \pi 8 \big) = \sqrt { \frac 1 2 - \frac 1 4 \sqrt 2}.$$
Applying this formula one can reduce a polynomial in $a$ of 
degree 6   to  a polynomial of degree 4.   One cannot entirely avoid 
cyclotomic fields, however, as one must prove that $a$ really has degree 4 over $\Q$.
Given that, the polynomial (\ref{eq:4759}), after reduction to degree 4, must 
be a multiple of (\ref{eq:4764}), and that soon leads to a contradiction.
The  calculations are somewhat lengthier than the proofs
with automorphisms given above.    

\medskip


\begin{lemma}\label{lemma:pioverten}  Suppose the tile $T$ is a right triangle with $\alpha = \pi/10$, and 
suppose that triangle $ABC$  is $N$-tiled by $T$, and angle $C$ is at least a right angle.
  Then either $T$ is similar to $ABC$, or $ABC$ is isoceles
and $T$ is similar to half of $ABC$.
\end{lemma}

\noindent{\em Proof}.   The possible shapes of $ABC$ are as follows:
\smallskip

(i) $C$ is a right angle and angle $A = \alpha$ and angle $B = 4\alpha$.  In that case, $T$ is similar 
to $ABC$, so we are done.

(ii) Angle $C$ is $\pi/2 + 3\alpha$ and angles $A$ and $B$ are both $\alpha$.  In that case $ABC$ is isosceles
and $T$ is similar to half of $ABC$, so we are done.

(iii) Angle $C$ is $\pi/2 + \alpha = 3\pi/5 $ and angles $A$ and $B$ are both $2\alpha$.

(iv) $C$ is a right angle and angle $A = 2\alpha$ and angle $B = 3\alpha$. 

(v) Angle $C$ is $\pi/2 + \alpha = 3\pi/5$, and angle $A = \alpha$ and angle $B=3\alpha$.

(vi) Angle $C$ is $\pi/2 + 2\alpha = 7\pi/10$, and angle $A$ is $\alpha$ and angle $B=2\alpha$.
\smallskip

  Your favorite computer algebra system 
 will tell you that 
  \begin{eqnarray*}
 a \ = \ \sin \frac \pi {10} &=& \frac 1 4 (\sqrt 5 - 1) \\
 b \ =\  \cos \frac \pi {10} &=& \frac 1 2 \sqrt{\frac 1 2 (5 + \sqrt 5)}
 \end{eqnarray*}
 Therefore $\Q(a)$ has degree 2 over $\Q$ and $a^2$ can be expressed linearly in $a$:
 \begin{eqnarray*}
 (4a+1)^2 &=& 5 \\
 16a^2 + 8a + 1 &=& 5 \\
 4a^2 + 2a-1 &=& 0 
 \end{eqnarray*}
 and the result is
 \begin{equation}\label{eq:asqlinear}
 a^2 = \frac 1 4 - \frac 1 2 a
 \end{equation}
 Similarly $\Q(b)$ has degree 4 over $\Q$, so $b^4$ can be expressed in terms of lower powers of $b$:
 \begin{eqnarray*}
 4b^2 &=&\frac 1 2(5 + \sqrt 5) \\
 8b^2-5 &=& \sqrt 5 \\
 (8b^2-5)^2 &=& 5 \\
 64b^4 - 80b^2 + 25 &=& 5 \\
 16b^4 - 20b^2 + 5 &=& 0 
 \end{eqnarray*}
 and the result is
 \begin{equation}\label{eq:bfourth}
 b^4 = \frac 5 4 b^2 - \frac 5 {16}
 \end{equation}
 We can also see that $a$ belongs to $\Q(b)$:  a few lines above we showed $\sqrt 5 = 8b^2 -5$,
 and hence
 \begin{eqnarray*}
 a &=& \frac 1 4 (\sqrt 5 - 1) \\
   &=& \frac 1 4 (8b^2-5 - 1)
\end{eqnarray*}
with the result
\begin{equation}\label{eq:aquadraticinb}   
 a = 2b^2 - \frac 3 2
\end{equation}
Since $a$ belongs to $\Q(b)$, a basis for $\Q(a,b)$ over $\Q$ is $1$, $b$, $b^2$, and $b^3$.  We will 
need to express other quantities in that basis.  For example
\begin{eqnarray*}
ab^2 &=& (2b^2 - \frac 3 2) b^2 \\
&=& 2b^4 - \frac 3 2 b^2 \\
&=& 2\big( \frac 5 4 b^2 - \frac 5 {16}\big) - \frac 3 2 b^2 
\end{eqnarray*}
with the result
\begin{equation} \label{eq:absquared}
ab^2 =  b^2 - \frac 5 8
\end{equation}
We also have 
\begin{eqnarray*}
a^3 &=& a^2 a \\
&=& (\frac 1 4 - \frac 1 2 a) a  \mbox{\qquad by (\ref{eq:asqlinear})} \\
&=& \frac 1 4 a - \frac 1 2 a^2 \\
&=& \frac 1 4 a  - \frac 1 2\big(\frac 1 4 - \frac 1 2 a\big) \\
&=& \frac 1 2 a - \frac 1 8 \\
&=&\frac 1 2 \big(2b^2 - \frac 3 2\big) - \frac 1 8 
\end{eqnarray*}
with the result
\begin{equation} \label{eq:acubeinb}
a^3 = b^2 - \frac 7 8 
\end{equation}
We will also need an expression for $\sin 3\alpha$.  We work that out as follows:
\begin{eqnarray}
\sin 3\alpha &=& 3 \sin \alpha - 4\sin^2\alpha \nonumber \\
&=& 3a - 4a^3 \nonumber\\
&=& 3(2b^2 - \frac 3 2) - 4(b^2 - \frac 7 8 ) \mbox{\qquad by (\ref{eq:aquadraticinb}) and (\ref{eq:acubeinb})} \nonumber \\
\sin 3 \alpha &=& 2b^2 - 1 \label{eq:3326}
\end{eqnarray}
With these expressions in hand, we consider side $AC$ of triangle $ABC$, whose length we call $Y$. 
Let $g$, $m$, and $f$ (respectively) be the numbers of $a$ sides, $b$ sides, and $c$ sides of tiles 
that lie along $AC$ in the given tiling.   Since $c = \sin \gamma = 1$, we have 
$$ Y = ga + mb + f$$
We next work out $Y^2$ in the basis $1$, $b$, $b^2$, $b^3$.
\begin{eqnarray*}
Y^2 &=&(ga+mb+f)^2 \\
&=& \big(g\big(2b^2- \frac 3 2\big) + mb + f\big)^2 \mbox{\qquad by (\ref{eq:aquadraticinb})} \\
&=& \big(2gb^2 + mb  + \big(f - \frac {3g} 2\big)\big)^2 \\
&=& 4g^2 b^4 + 4gmb^3 + \big[ m^2 + 4g\big( f-\frac {3g} 2\big)\big] b^2 + 2m\big( f-\frac {3g}2\big) b + \big(f-\frac {3g}2 \big)^2 \\
&=& 4g^2 b^4 + 4gmb^3 + \big[m^2 -6g^2 + 4gf\big] b^2 + 2m\big( f-\frac {3g}2\big) b + (f^2 -3gf + \frac 9 4 g^2)
\end{eqnarray*}
Substituting in for $b^4$ from (\ref{eq:bfourth}) we have
\begin{eqnarray*}
Y^2 &=& 4g^2\big(\frac 5 4 b^2 - \frac 5 {16}\big) + 4gmb^3+ \big[m^2 -6g^2 + 4gf\big] b^2 + 2m\big( f-\frac {3g}2\big) b + (f^2 -3gf + \frac 9 4 g^2) 
\end{eqnarray*}
and writing it as a polynomial in $b$ we find
\begin{equation}\label{eq:Xsqinb}
Y^2 =  4gmb^3 + (m^2 -g^2 + 4gf)b^2 + 2m\big( f-\frac {3g}2\big) b + (f^2 -3gf +  g^2)
\end{equation}
We use the letters $A$ and $C$ to stand for the angles of $ABC$ at vertices $A$ and $C$, as well as for the 
vertices themselves.  Let $X$ be the length of side $BC$, opposite angle $A$.   Then by the law of sines we have
$$ \frac Y {\sin B} = \frac X { \sin A} $$
Therefore 
$$ X = Y \frac {\sin A} { \sin B} $$
The area of triangle $ABC$ is given, according to the usual cross product formula, as
$$ \A_{ABC} = \frac 1 2 XY \sin C$$
Substituting for $X$ we have 
$$ \A_{ABC} = \frac 1 2 Y^2 \frac {\sin A \sin C} {\sin B} $$
Let $\A_T = ab/2$ be the area of the tile $T$. Our fundamental equation is
$$N \A_T = \A_{ABC}$$
That is  
$$ Nab = Y^2 \frac {\sin A \sin C} {\sin B} $$
(Since $\A_T = ab/2$, the factors of $1/2$ cancel and disappear.)
Now we take up the argument by cases according as the possible shapes of $ABC$, as enumerated at the beginning 
of the proof.   First we take up case (iii), where angles $A$ and $B$ are each $2\alpha$.   Then $\sin B$ and 
$\sin A$ are equal, and cancel out, so we have 
$$ N ab = Y^2 \sin C$$
Since in this case angle $C = \pi/2 + \alpha$, we have $\sin C = \cos \alpha = b$, so we can cancel $b$ on 
the left with $\sin C$ on the right, obtaining
$$ N a = Y^2.$$
Putting in on the right the expression for $Y^2$ found in (\ref{eq:Xsqinb}), and using (\ref{eq:aquadraticinb}) on the left,
 we have 
$$ N\big(2b^2 - \frac 3 2\big) =  4gmb^3 + (m^2 -g^2 + 4gf)b^2 + 2m\big( f-\frac {3g}2\big) b + (f^2 -3gf + g^2)
$$
Since both sides are expressed in the basis $\{1,b,b^2,b^3\}$,  the coefficients of like powers of $b$ are equal.
Equating the coefficients of $b^3$ we find $gm=0$.
If $g=0$ then from the coefficients of $b^2$ we have $2N = m^2$.  Then from the coeficients of $b$ we have $mf = 0$,
and since $m \neq 0$ we have $f=0$.  Then the constant coefficient on the right is zero, but on the left it is $-3N/2$.
Hence $g \neq 0$.  But since $gm=0$ that implies $m=0$.  Then our equation becomes
$$ N\big(2b^2 - \frac 3 2\big) =  (4gf -g^2) b^2 +   (f^2 -3gf + g^2)
$$
Equating coefficients of $b^2$ we have 
\begin{eqnarray}
2N &=&  4gf-g^2 \label{eq:3375}
\end{eqnarray}
Equating the constant terms we have 
\begin{eqnarray*}
-\frac {3N}2 &=& \ f^2-3gf + g^2\\
3N&=& 2(3 gf - g^2 - f^2)  \\
&=& 2(4 gf - g^2 -gf -f^2) 
\end{eqnarray*}
 Using (\ref{eq:3375}) on the right hand side we find
\begin{eqnarray*}
3N &=& 2(2N-gf-f^2) \\
&=& 4N-2(gf+f^2) \\
N &=& 2(gf+f^2)\\
2N &=& 4gf + 4f^2
\end{eqnarray*}
But also $2N = 4gf - f^2$, by (\ref{eq:3375}).  Subtracting this from the last equation we have $0 = 4f^2 + g^2$.
Hence $g=f=0$.  But this contradicts $m=0$, since $g+m+f$ is the number of tiles along side $AC$.
That contradiction completes the proof in case (iii), when $ABC$ is isosceles with base angles $2\alpha$.
\smallskip

Next we take up case (iv), when angle $C$ is a right angle, angle $A$ is $2\alpha$, and angle $B$ is $3\alpha$.  Then we have
\begin{eqnarray*}
2N\A_T &=& Y^2 \frac {\sin A \sin C} {\sin B} \\
N ab &=& Y^2 \frac {\sin 2\alpha}{\sin 3\alpha}
\end{eqnarray*}
We have $\sin 2\alpha = 2 \sin \alpha \cos \alpha = 2ab$, so the $ab$ on the left cancels with $\sin 2\alpha$ on the right:
$$ N = \frac {2Y^2} {\sin 3 \alpha}$$
We have $\sin 3 \alpha = 3\sin \alpha \cos^2 \alpha - \sin^2 \alpha = 3ab^2 - a^3$, so 
\begin{eqnarray*}
N(3ab^2-a^3) &=& 2Y^2 
\end{eqnarray*}
We must express the left side in terms of the basis.  We have already done the work; substituting the values for $ab^2$ and $a^3$ found in (\ref{eq:absquared}) and (\ref{eq:acubeinb}) we have
\begin{eqnarray*}
 N(3(b^2 - \frac 5 8) - (b^2 - \frac 7 8)) &=&  2Y^2 \\
 N(2b^2-1) &=& 2Y^2\\
 N(b^2 - \frac 1 2) &=& Y^2
\end{eqnarray*}
and substituting the value for $Y^2$ found in (\ref{eq:Xsqinb}) we have
\begin{eqnarray*}
N(b^2-\frac 1 2) &=&  4gm b^3 +  (m^2-g^2+4gf) b^2 + 2m\big(f-\frac{3g} 2\big) b +(f^2 - 3gf + g^2)
\end{eqnarray*}
The right hand side is the same as in the previous case, but the left hand side is slightly different.  As before,
the coefficient of $b^3$ must be zero on the right, so $gm=0$, and as before, if $g=0$ then we compare coefficients
of $b^2$.  This time we find $N = m^2$.  Comparing the coefficients of $b$ we find $mf=0$, but since $m \neq 0$ 
we have $f=0$, and hence the constant term on the right is zero, while on the left it is $-N/2$.  Therefore $g \neq 0$;
but since $gm=0$, we have $m = 0$.   Then from equating the coefficients of $b^2$ we have 
\begin{equation} \label{eq:N38}
N = 4gf - g^2
\end{equation}
and by equating the constant terms we find
\begin{eqnarray*}
-\frac N 2 &=& f^2 -3gf + g^2\\
\frac N 2 &=& 3gf - g^2 -f^2 \\
&=& 4 gf - g^2 -gf -f^2 
\end{eqnarray*}
 Using (\ref{eq:N38}) on the right hand side we find
\begin{eqnarray*}
\frac N 2 &=& N-gf-f^2 
\end{eqnarray*}
Solving for $N$ we have 
\begin{equation} \label{eq:N39}
 N = 2gf + 2f^2  
\end{equation}
 Now we have two expressions equal to $N$, namely (\ref{eq:N38}) and (\ref{eq:N39}).
 Equating them, we have 
 \begin{eqnarray*}
 2gf + 2f^2 &=& 4gf - g^2 \\
 2f^2 &=& 2gf - g^2 \\
 f^2 + (f^2-2gf + g^2) &=& 0 \\
 f^2 + (f-g)^2 &=& 0
 \end{eqnarray*}
 It follows that $f$ and $f-g$ are both zero.  Hence $f$ and $g$ are both zero, which contradicts $m=0$ since $g+m+f$ 
 is the number of tiles along side $AC$.  That completes the proof in case (iv).
 \smallskip
 
 Next we take up case (v), in which angle $C$ is $\pi/2 + \alpha$, angle $A$ is $\alpha$, and angle $B$ is $3\alpha$.  Then we have
 \begin{eqnarray*}
Nab &=& Y^2 \frac {\sin A \sin C} {\sin B} \\
  &=& Y^2 \frac {\sin \alpha \cos \alpha}{\sin 3\alpha} \\
     &=& Y^2 \frac {ab} {\sin 3 \alpha}
\end{eqnarray*}
and canceling $ab$ we have
\begin{eqnarray*}
N &=& \frac {Y^2}{\sin 3 \alpha}
\end{eqnarray*}
This equation differs from case (iv) only by a factor of 2 on the right side, but it is not exactly the same,
so we must carry out the calculations.  Expressing $\sin 3\alpha$ in terms of $b$  by means of (\ref{eq:3326}),
we find  
\begin{eqnarray*}
N(2b^2-1) &=& Y^2
\end{eqnarray*}
Substituting for $Y^2$ as before, we find 
\begin{eqnarray*}
N(2b^2-1) &=&  4gm b^3 + (m^2-g^2 + 4gf)b^2 + 2m\big(f-\frac{3g} 2\big) b + (f^2-3gf + g^2)
\end{eqnarray*}
Comparing the coefficients of $b^3$ we again find $gm=0$.  Assume, for proof by contradiction, that $m\neq 0$.
Then from the coefficients of $b$,  if $g=0$ then
$f=0$ and vice-versa.    If $f=g=0$, then the constant term is zero on the right but $-N$ on the left.  Hence 
$g \neq 0$.  But since $gm = 0$
that contradicts the assumption $m \neq 0$.  That contradiction shows that  $m = 0$.

Equating the coefficients of $b^2$ we have 
\begin{equation} \label{eq:N40}
2N = 4gf-g^2
\end{equation}
 and from the constant terms we have 
\begin{eqnarray*}
N &=& 3gf-g^2-f^2 \\
&=& 4gf - g^2 - gf - f^2 \\
&=& 2N - gf - f^2   \mbox{ \qquad by (\ref{eq:N40})}
\end{eqnarray*}
Solving for $N$ we have
\begin{eqnarray*}
N &=& gf + f^2 \\
2N &=& 2gf + 2f^2
\end{eqnarray*}
and equating this  expression for $2N$ with the one in (\ref{eq:N40}), we have 
\begin{eqnarray*}
 4gf - g^2 &=& 2gf + 2f^2 \\
 0 &=& g^2 -2gf + 2f^2  \\
   &=& (g-f)^2 + f^2
\end{eqnarray*}
Hence $f=0$ and $g=0$; and since we already have $m=0$, this is a contradiction.  That completes case (v).

The last case is case (vi), in which angle $C$ is $\pi/2 + 2\alpha = 7\pi/10$, and angle $A$ is $\alpha$ and angle $B=2\alpha$.
Then we have 
\begin{eqnarray*}
Nab &=& Y^2 \frac {\sin A \sin C} {\sin B} \\
  &=& Y^2 \frac {\sin \alpha \cos 2\alpha}{\sin 2\alpha} \\
  &=& Y^2 \frac { \sin \alpha (\cos^2 \alpha - \sin^2 \alpha)}{2\sin \alpha \cos\alpha}\\
     &=& Y^2 \frac {a(b^2-a^2)} {2ab} 
\end{eqnarray*}
Canceling $a$ and multiplying by the denominator, we have 
\begin{eqnarray*}
2Nab^2 &=& Y^2(b^2-a^2)
\end{eqnarray*}
Putting in the value for $a^2$ from (\ref{eq:asqlinear}) on the right, and the value of $ab^2$ from (\ref{eq:absquared}) on the left, 
 we find 
\begin{eqnarray*}
2N\big( b^2 - \frac 5 8\big) &=& Y^2\big(b^2- \big(\frac 1 4 - \frac 1 2 a\big)\big) \\
2Nb^2 - \frac 5 4 N &=& Y^2\big(b^2 - \frac 1 4 + \frac 1 2 \big(2b^2 - \frac 3 2\big)\big) \mbox{\qquad by (\ref{eq:aquadraticinb})} \\
&=& Y^2(2b^2 -1)
\end{eqnarray*}
Substituting for $Y^2$ from (\ref{eq:Xsqinb}) we have
\begin{eqnarray*}
2Nb^2 - \frac 5 4 N &=& (2b^2-1)\bigg\{
4gm b^3 + (m^2-g^2+4gf) b^2 + 2m\big(f-\frac{3g} 2\big) b + (f^2 -3fg + g^2)\bigg\} \\
&=& 8gm b^5 + (m^2-g^2 + 4gf)2b^4 + 4m\big(f-\frac{3g} 2\big) b^3 - 4gmb^3  \\
&& +  ( -m^2 +g^2 -4gf  +2(f^2-3fg + g^2)) b^2  \\
&&-2m\big(f-\frac{3g} 2\big) b -  (f^2 - 3fg + g^2) 
\end{eqnarray*}
Multiplying  (\ref{eq:bfourth}) by $b$ we have 
$$
 b^5 = \frac 5 4 b^3 - \frac 5 {16}b
$$
We use this to eliminate the $b^5$ term, and also we use (\ref{eq:bfourth}) to eliminate the $b^4$ term:
\begin{eqnarray*}
2Nb^2 - \frac 5 4   N
&=& 8gm \big( \frac 5 4 b^3 - \frac 5 {16}b \big) + (m^2-g^2 + 4fg)2\big(\frac 5 4 b^2 - \frac 5 {16}\big) \\
&& + 4m\big(f-\frac{3g} 2\big) b^3 - 4gmb^3  \\
&& +  (2f^2 -10gf + 3g^2 -m^2  ) b^2  \\
&&-2m\big(f-\frac{3g} 2\big) b -  (f^2 - 3fg + g^2) 
\end{eqnarray*}
Expressing the right hand side as a polynomial in $b$, we have 
\begin{eqnarray*}
2Nb^2 - \frac 5 4   N &=& 4f m b^3+(g^2/2+2f^2+3m^2/2) b^2\\
&& +(g m/2-2f m) b \\
&& -3g^2/8+f g/2-f^2-5m^2/8
\end{eqnarray*}
Equating the coefficients of $b^3$ we find $fm = 0$.
Hence either $m=0$ or $f=0$.  Suppose, for proof by contradiction, 
 that $f=0$.  Then the coefficient of $b$ on the right is $gm/2$, but on the left it is zero.  Hence 
 $g=0$ or $m=0$.  If $g=0$ then   equating the constant terms we 
 have 
 $$ \frac 5 4 N = \frac {5m^2} 8.$$
 Then $2N = m^2$. Equating the quadratic terms we have 
 \begin{eqnarray*}
 2N &=&  \frac {3m^2} 2 \\
 &=& 3N \mbox{\qquad since $m^2 = 2N$}
 \end{eqnarray*}
 This is a contradiction, so $g \neq 0$.  Still assuming $f=0$, we now have $m=0$. 
 Equating the constant terms we have 
 \begin{eqnarray*}
 \frac {5N} 4 &=& \frac{3g^2} 8 \\
 10 N &=& 3g^2 
 \end{eqnarray*}
 Equating the quadratic terms we have 
 \begin{eqnarray*}
 2N &=& \frac {g^2} 2 \\
 10 N &=& \frac 5 2 g^2
 \end{eqnarray*}
 Now $10N $ is equal to both $3g^2$ and $(5/2) g^2$, which is impossible.  This contradiction shows 
 that $f \neq 0$.  Then from the cubic terms we have $m=0$.  From the quadratic and constant terms
 we have 
 \begin{eqnarray*}
 2N &=& g^2/2 + 2f^2 \\
 \frac 5 4 N &=& f^2 - \frac {fg} 2 + \frac{3g^2} 8
 \end{eqnarray*}
 Multiplying the first of these equations by 5 and the second by 8 we have 
 \begin{eqnarray*}
 10 N &=& \frac 5 2 g^2 + 10 f^2 \\
 10 N &=& 8f^2 - 4fg + 3g^2 
 \end{eqnarray*}
 Subtracting the second equation from the first we have 
 \begin{eqnarray*}
 0 &=& 2f^2 + 4fg - \frac 1 2 g^2 
 \end{eqnarray*}
 Let $\xi = g/f$, which is legal since $f \neq 0$.  Dividing the previous equation by $f^2$ we have 
 \begin{eqnarray*}
 0 &=& 2 + 4\xi - \frac 1 2 \xi^2
 \end{eqnarray*}
 The discriminant of this quadratic equation is $\sqrt{16-4} = \sqrt{12} = 2\sqrt 3$, which is 
 irrational.  Hence there are no integer solutions $g$ and $f$ of the equations above.  That
 completes the proof of the lemma.

\begin{theorem} \label{theorem:gammaright}
Suppose that the triangle $ABC$ is $N$-tiled by a non-isosceles right triangle $T$. 
Then either 

(i)  $ABC$ is similar to $T$, or 
\smallskip

(ii) $ABC$ is isosceles, and $T$ is similar to half of $ABC$.  In this case exactly one 
vertex angle splits, and it splits into exactly two pieces. 
\smallskip

(iii) $ABC$ is equilateral, and $T$ is a 30-60-90 triangle, and $N$ is six times a square. 
\end{theorem}

\noindent{\em Proof}.     Let 
$P$, $Q$, and $R$ be the total number of $\alpha$, $\beta$,
and $\gamma$ angles of tiles at the vertices of $ABC$. 
   Since $\gamma = \pi/2$ we have $R \le 1$. 
We begin by showing that if $R=1$, then the specific value $\alpha = \pi/8$ implies conclusion (i) or (ii) of
the theorem.   Assume $R=1$. Then by the definition of $R$,
 one tile has a right angle at a vertex of $ABC$; it must be at angle $C$, the largest angle of triangle $ABC$.
  Assume that $\alpha = \pi/8$. 
  If  vertex $C$ is split, and one $\alpha$ angle occurs at $C$ in addition to the right angle, then  angle $C$ is $5\pi/8$,
and one of the two remaining angles is $\alpha$.  That case has been ruled out in Lemma \ref{lemma:special2}.   If
the angle at vertex $C$ is split, and two additional $\alpha$ angles occur there, then triangle $ABC$ is isosceles 
with base angles $\alpha$.  So $T$ is similar to half of $ABC$, which is conclusion (ii) of the 
theorem.  Therefore we may assume that the angle at vertex $C$
is not split (so only the right angle occurs there). If each of the other two vertices is split into two $\alpha$ angles,
then angles $A$ and $B$ are equal,  and $ABC$ is a right isosceles triangle.  This case has been ruled out in 
Lemma \ref{lemma:piover8not}.   If the four $\alpha$ angles are distributed one and three, then $ABC$ is similar to 
$T$, which is conclusion (i) of the theorem.  That disposes of the 
case $\alpha = \pi/8$ and $R=1$.   
     
We will prove that either one of the conclusions  of the theorem holds, or $P+Q+R=4$.   
First we take up the case $R=0$.  Then we have
\begin{eqnarray*}
 P \alpha + Q \beta &=& \pi \\
\alpha + \beta = \pi/2
\end{eqnarray*}
We note that  $P=Q=2$ does provide a solution of those two equations.  We shall show that in every other 
case, the conclusion of the theorem holds (or a contradiction ensues).
We begin by considering the possibilities for $Q$.
We have $\beta < \pi/2$ and $\alpha \le \beta$, and $\alpha + \beta = \pi/2$.  Therefore $\beta \ge \pi/4$.  From 
$P\alpha  + Q\beta = \pi$ we then obtain $Q \le 4$.   If $Q=4$ then $\beta = \alpha = \pi/4$, contradicting the assumption that $T$
is not isosceles. Hence $Q < 4$.   

Next we consider the case $Q=3$.  
Then subtracting the equation $\alpha + \beta = \pi/2$ from the equation $P \alpha + Q\beta = \pi$ we obtain
$$(P-1) \alpha + 2 \beta = \pi/2.$$
   If $P > 1$ then the left side exceeds $\alpha + \beta = \pi/2$, contradiction.   If $P=1$ then $\beta = \pi/4$, but then $\alpha = \pi/4$,
contradicting $\alpha < \beta$.    Hence we must have $P=0$.  Then $-\alpha + 2\beta = \pi/2$.  Together with $\alpha + \beta = \pi/2$
this implies that $\alpha = \pi/6$ and $\beta = \pi/3$.   With $R=0$ and $P=0$, all the angles of $ABC$ must be composed of three 
$\beta$ angles, so $ABC$ is equilateral, and there is no vertex splitting.  Then case (ii) holds, since an equilateral triangle 
is isoceles and $T$ is similar to half of $ABC$.  
That disposes of the case $Q=3$. 

Next we claim that if $Q=1$, then $P \ge 4$, and if $Q=0$ then $P > 4$.
First suppose $Q=1$. Then $P\alpha + \beta = \pi$.
Subtracting $\alpha + \beta = \pi/2$ we have $(P-1)\alpha = \pi/2$.  Since $\alpha < \pi/4$, we have $P-1 > 2$, i.e. $P \ge 4$, as claimed.
Now assume $Q=0$.  Then $P\alpha + Q \beta = \pi$ implies $\alpha  = \pi/P$; but since $\alpha < \pi/4$, we have $P > 4$.  This 
completes the proof of the claim.  

Now  we consider the cases $Q=0$ or $Q=1$ further (still assuming $R=0$).
In that case $ P \ge 4$, so in particular $P > Q$.  That is, there are more $\alpha$ angles
than $\beta$ angles at the vertices of $ABC$.  Therefore, there is another vertex $V$ of the tiling at which there are more $\beta$ angles 
than $\alpha$ angles  (since there are $N$ of each altogether).   Fix such a vertex $V$, and let $k$ be 1 if $V$ is a boundary or a non-strict
vertex, and $k=2$ if $V$ is a strict interior vertex; so the angle sum at $V$ is $k\pi$.   Let $n$, $m$, and $\ell$ be the numbers of 
$\alpha$, $\beta$, and $\gamma$ angles (respectively) at $V$.  Then we have $m > n$ and 
$$ \matrix n m \ell 1 1 1 P Q 0  \vector \alpha \beta \gamma = \vector {k\pi} \pi \pi $$
Solving for $\gamma$ by Cramer's rule (or Gaussian elimination) we find
\begin{equation}
[ (m-\ell) P - Q(n-\ell)] \gamma = [n(1-Q) + m(P-1) + k(Q-P)] \pi. \label{eq:3674}
\end{equation}
As long as we don't divide both sides by the coefficient of $\gamma$, the equation is valid whether or not 
that coefficient is zero.
Remembering that $\gamma = \pi/2$ we have 
\begin{eqnarray} \label{eq:lmnPQ}
(m-\ell) P - Q(n-\ell) &=& 2[n(1-Q)+m(P-1)+k(Q-P)]
\end{eqnarray}
Assume $Q=0$.  Then we will prove $\beta = \pi/3$ and $\alpha = \pi/6$.  
Since $Q=0$, we have $P > 4$ as proved above. We have 
\begin{eqnarray}
(m-\ell) P   &=& 2[n +m(P-1)-kP] \nonumber \\
P(m-\ell -2m + 2k) &=& 2(n-m)\nonumber \\
P(-\ell -m + 2k) &=& 2(n-m) \nonumber\\
P(\ell + m - 2k) &=& 2(m-n) \nonumber\\
P &=& \frac {2(m-n)}{\ell + m -2k} \label{eq:Pbound33}
 \end{eqnarray}
 The right side is positive since $m > n$, so the left side is also positive.  Since we 
 still are assuming $R=0$, the assumption $Q=0$ means that all the vertices of $ABC$ are composed
 only of $\alpha$ angles, so $\pi = P\alpha$, i.e. $\alpha = \pi/P$.
 We claim $k$ cannot be 1.    Suppose that $k=1$ and $\ell = 0$; then $\pi = n\alpha + m\beta$, and since $\beta > \pi/4$
 we have $m = 0,1,2$, or 3.  Since $2\beta + 2\alpha = \pi$ we must have $m=3$ in order that $m > n$.
 Then $n=1$ or $n=0$.   If $n=0$ then $\beta = \pi/3$ as desired. 
 Hence we may suppose $n=1$; then $n-m=2$, so $P = 2(m-n)/(\ell + m - 2k) = 4$,  contradicting $P > 4$.  On the other hand 
 if $k=1$ and $\ell = 1$ we have only an angle of $\pi/2$ to fill with $\alpha$ and $\beta$ angles, 
 and it cannot be done with $m > n$.  That proves that $k$ cannot be 1.  
 
 So $k=2$.   Then 
 $P = 2(m-n)/(\ell + m - 4)$.  We cannot have $\ell = 4$, as in that case there would be four right 
 angles at $V$, and hence there would be no $\alpha$ and $\beta$ angles at $V$, contradicting $m > n$. 
 If $\ell =3$ and $n > 0$ then there is one $\alpha$ and one $\beta$ at $V$,  contradicting $m > n$.  
 If $\ell = 3$ and $n = 0$ then $\pi/2$ is made of $m$  angles of measure $\beta$, i.e. $\beta = \pi/m$.
 If $m=3$  our claim that $\beta = \pi/3$ holds;  
 and since $\alpha < \beta$ and $\alpha + \beta = \pi/2$, we have $\beta > \pi/4$, and hence $\beta \neq \pi/m$
 for any $m$ (except possibly 3, which we are trying to prove must be the case). 
 If $\ell = 2$ and $n > 0$ 
 then there must be 3 $\beta$ angles and one $\alpha$ at $V$, i.e. $m=3$ and $n=1$ (since if $n = 2$ then 
 $m=2$ because $2\alpha + 2\beta + 2\gamma = 2\pi$, so $n < 2$).   
  \begin{eqnarray*}
 P &=& 2(m-n)/(\ell + m - 4)  \\
 &=& 2/(2+3-4) \\
 &=& 2
 \end{eqnarray*}
 contradicting $P > 4$.  If $\ell = 2$ and $n=0$ we have  $m\beta = \pi$, and the only possibility 
 with $\beta > \pi/4$ is $\beta = \pi/3$, which is the desired conclusion. 
 If $\ell = 1$ then there is $3\pi/2$ to make up from $\alpha$ and $\beta$ angles,  and since 
 $\alpha = \pi/P$,  we have $\beta = \pi/2 -\alpha = \pi/2 - \pi/(2P)$.
 Then $n\alpha + m \beta = 3\pi/2$ becomes
  $$ n \frac \pi P  + m \big( \frac \pi 2 - \frac \pi P\big) = \frac {3\pi} 2.$$
 Solving this equation for $n$ we find
 $$n = 3P/2 - m(P-2)/P.$$
 Then 
 \begin{eqnarray*}
 P &=& \frac{2(m-n)}{\ell + m - 4} \\
 &=& \frac {2(m- (3P/2 - m(P-2)/P))}{1+m-4} \\
 &=& \frac {2m -3P + 2m(P-2)/P}{m-3} \\
 P^2(m-3) &=&  2mP - 3P^2 + 2m(P-2) \\
 mP^2 - 4mP +4m &=& 0 \\
 P^2 - 4P + 4 &=& 0   \mbox{\qquad since $m > n$, so $m \neq 0$}\\
 (P-2)^2  &=& 0 \\
 P &=& 2
 \end{eqnarray*}
 But that contradicts $P > 4$, so $\ell = 1$ is impossible.

 Hence $\ell = 0$.  Then the angle of $2\pi$ at $V$ is made up entirely of $\alpha$ and $\beta$ angles.  Since
 we still have $R=Q=0$,  $\alpha = \pi/P$, so $n\alpha + m \beta = 2\pi$ becomes
 $$ n \big(\frac\pi P\big) +m \big( \frac \pi 2 - \frac \pi P\big) = 2\pi.$$
 Solving for $n$ we find 
 $$ n = 2P - m \frac { 2-P} P. $$
Then 
 \begin{eqnarray*}
  P &=& \frac{2(m-n)}{\ell + m - 4}  \\
  &=& \frac{2(m-n)}{m-4} \mbox{\qquad since $\ell = 0$} \\
  &=& \frac{2(m- (2P-m(P-2)/P))}{m-4} \\
 (m-4)P^2 &=& 2mP -4P^2 +2m(P-2)    \\
 mP^2 &=& 4mP - 4m  \\
P^2-4P + 4 &=& 0 \mbox{\qquad since $m \neq 0$} \\
(P-2)^2 &=& 0 \\
P &=& 2
 \end{eqnarray*}
 But that contradicts $P > 4$, so $\ell= 0$ is also impossible. 
 That exhausts all the possibilities for $\ell$ and $k$.  Hence $Q=0$ is possible only when $\beta = \pi/3$,
 as claimed.

Now assume, for proof by contradiction, that $Q=1$.  Then we have by (\ref{eq:3674}) 
\begin{eqnarray*}
(m-\ell) P - (n-\ell) &=& 2[m(P-1)+k(1-P)] \\
\ell(1-P) +mP - n  &=& 2(m-k)(P-1)  \\
(P-1)(m-\ell) + (m-n) &=& (P-1)(2m-2k)\\
m-n &=& (P-1)(m-2k+\ell) \\
2k-\ell &=& \frac { n + (P-2)m}{P-1} \\
&=& \frac{ n -m  + (P-1)m}{P-1} \\
&=&   m - \frac {m-n}{P-1}
\end{eqnarray*}
The numerator $m-n$ is positive, and the denominator $P-1$ is at least 3, so $P-1$ divides $m-n$.
Since $m - n$ is positive, we have $m-n \ge P-1 \ge 3$, so $m \ge 4$.   Since $\beta > \pi/4$ and 
$m \ge 4$, we have $m\beta > \pi$, so the $\beta$ angles at vertex $V$ occupy more than a straight angle.
Hence $k \neq 1$.  Hence $k=2$, and also $\ell = 0$ or $\ell = 1$, since more than one right angle cannot 
be accommodated at $V$ (since $m \ge 4$, and $4\beta > \pi$). 
 Substituting $2$ for $k$ we have
\begin{eqnarray*}
4-\ell &=& m - \frac {m-n}{P-1} \\
(P-1)(4-\ell) &=& (P-1)m - (m-n) \\
m(P-2) &=& (P-1)(4-\ell) - n \\
m &=& \frac {P-1}{P-2} (4-\ell) - \frac n {P-2}
\end{eqnarray*}
Since $P \ge 4$, the fraction $n/(P-2)$ is positive, and the fraction $(P-1)/(P-2)$ is between 0 and 1, so 
\begin{eqnarray*}
m &<& \frac{P-1}{P-2} (4-\ell) \\
  &<& 4-\ell \le 4
\end{eqnarray*}
But this contradicts the result $m \ge 4$ that we proved above.   This contradiction completes the
proof that $Q\neq 1$.

Now assume $Q=2$.
We  claim that $P=2$. Subtracting $\alpha + \beta = \pi/2$ from $P\alpha + Q\beta = \pi$, and using $Q=2$, we have
$$(P-1)\alpha + \beta = \pi/2$$
But since $\alpha + \beta = \pi/2$, it follows that $P=2$ as claimed.   That completes the proof that $P=2$, 
in case $R=0$ and $Q=2$.
 
We now consider the case $R=1$.  We claim $P+Q+R = 4$.
Assume first that no 
vertex splitting occurs, i.e. each angle of $ABC$ is equal to $\alpha$, $\beta$, or $\gamma$.
Since $R=1$, angle $C = \gamma = \pi/2$.   Then angles $A$ and $B$ are each equal to $\alpha$ or $\beta$; but 
since $T$ is not isosceles, we must have $B = \beta$ and $A = \alpha$, so $ABC$ is similar to $T$, and case (i)
of the theorem holds.  Hence we may assume from now on (under the assumption $R=1$) that some vertex splitting does occur;
and we proved already that when $R=0$, some vertex splitting occurs.  Therefore  $P+Q \ge 3$.

We have $P\alpha + Q\beta = \pi/2$.  Subtracting the equation 
$\alpha + \beta = \pi/2$, we have $(P-1)\alpha + (Q-1)\beta = 0$.  Since $P+Q \ge 3$, we cannot have $P=Q=1$;  hence 
either $P=0$ or $Q=0$.  If $P=0$ then $\alpha = (Q-1)\beta$, which is impossible since $\alpha < \beta$.   Hence $Q=0$
and $\beta = (P-1)\alpha$.  Since $\alpha + \beta = \pi/2$, we have $\alpha = \pi/(2P)$.   Since $Q=0$, there are 
more $\alpha$ angles than $\beta$ angles at the vertices of $ABC$.  Therefore there exists a vertex $V$ at which 
there are more $\beta$ angles than $\alpha$ angles, i.e. at which $m > n$.   We have
\begin{eqnarray*}
n \alpha + m \beta + \ell \gamma &=& k \pi   \mbox{\qquad where $k=1$ or $k=2$}\\
n \frac \pi {2P} + m \big(\frac \pi 2 - \frac \pi {2P}\big) + \ell \frac \pi 2 &=& k \pi 
 \mbox{\qquad since $\alpha = \frac \pi {2P} $ }\\
\end{eqnarray*}
Dividing by $\pi/2$ we have
\begin{equation} \label{eq:344}
\frac n  P + m \big(1- \frac 1 {P}\big) + \ell   = 2k    
\end{equation}
Now we consider the possible values of $\ell$.   We have $\ell \le 3$ since $ \ell \gamma < 2\pi$ and $\gamma = \pi/2$.
  Since $m > n$ there 
is at least one $\beta$ angle at $V$.  If $\ell = 3$ then we must have $m=n=1$, since $T$ is not isosceles and $\alpha + \beta = \pi/2$;
but since $m > n$, that is ruled out.  

Consider the case $\ell = 2$.  That means $\pi$ of the angle at $V$ is used up by two $\gamma$ angles.  Since $m > 0$, there is at least one $\beta$ angle at $V$, so the total angle at $V$ is more than $\pi$.  Therefore  $k=2$ and we must have $\alpha$ and $\beta$ angles 
adding up to $\pi$.  It cannot be two of each since $m > n$.  Hence there are at least three $\beta$ angles.  Since
$\beta > \pi/4$,  there cannot be as many as four $\beta$ angles, as that would make more than $2\pi$ total at $V$.
Hence there are exactly three $\beta$ angles, i.e. $m=3$.  Equation (\ref{eq:344}) becomes
$$ \frac n P + 3\big(1 - \frac 1 P \big) + 2 = 4 $$
Solving for $P$ we find $P=3-n$.  Since $P+Q \ge 3$ and $Q=0$, we have $P \ge 3$.  Hence $n=0$ and $P=3$.
Since $R=1$ and $Q=0$, we now have $P+Q+R = 4$, 
 which is what we are trying to prove.  That disposes 
of the case $\ell = 2$.

Now consider the case $\ell = 1$.
Then we have $k=2$, since if $k=1$ we must have $\alpha$ and $\beta$ angles adding up to $\pi/2$, with more $\beta$ than 
$\alpha$ angles, but that is not possible, since one $\beta$ angle leaves room for only one $\alpha$ angle.   With $k=2$ and $\ell = 1$, 
we must have $\alpha$ and $\beta$ angles adding up to $3\pi/2$, with more $\beta$ than $\alpha$ angles.  We cannot have 
$m > 4$ since $5\beta = 5\pi/2 - 5\alpha = 5\pi/2 - 5 \pi/(2P) > 3\pi/2$, since $P \ge 3$.   By
(\ref{eq:344}) we have 
\begin{eqnarray*}
\frac n  P + m \big(1- \frac 1 {P}\big) + 1   &=& 4     \\
\frac n  P + m \big(1- \frac 1 {P}\big)   &=& 3     \\
n &=& (3-m)P+m
\end{eqnarray*}
With $m\le 3$ we have $n = (3-m)P + m \ge m$, contradicting $m > n$.   With $m=4$ we have 
$n = 4-P$.  Since $n \ge 0$ and $P \ge 3$, the only possibilities are $n=0$ and $P=4$, or $n=1$ and $P=3$.
If $P=3$ then $P+Q+R = 4$, which is what we are trying to prove. Hence we may assume $n=0$, $m=4$, and $P=4$.
Then $\alpha = \pi/8$ and $\beta = 3\pi/8$.  But at the beginning of the proof we already 
disposed of the case $R=1$ and $\alpha = \pi/8$.

Finally we consider the case $\ell = 0$.   First consider $k=1$;  then the total angle at $V$ is $\pi$
and none of it is used by right angles, so $n\alpha + m\beta = \pi$.   Equation (\ref{eq:344}) becomes
$$ \frac n P + m\big( 1  - \frac 1 P\big) =2 $$  
which implies 
\begin{eqnarray*}
 m\big(1-\frac 1 P\big) &\le& 2 \\
 m & \le& \frac {2P} {P-1} \\
 m&<& 3
 \end{eqnarray*}
That implies $m \le 2$, since $m$ is an integer.   But since $m > n \ge 0$, we must have $m=1$ and $n=0$, or 
$m=2$ and $n=1$, or $m=2$ and $n=0$.   If $m=1$ and $n=0$ we have
\begin{eqnarray*}
2 &=& \frac n P + m\big(1 - \frac 1 P \big) \\
  &=& \big( 1 - \frac 1 P \big) \\
  &<& 1
\end{eqnarray*}
so the case $m=1$ and $n=0$ is impossible.  If $m=2$ we have 
\begin{eqnarray*}
2 &=& \frac n P + m\big(1 - \frac 1 P \big) \\
  &=& \frac n P + 2\big( 1 - \frac 1 P \big) \\
0 &=& \frac n P - \frac 2 P \\
&=& \frac {n-2} P \\
n &=& 0
\end{eqnarray*}
Then the equation $n\alpha + m\beta = \pi$ becomes $2\beta = \pi$; but since $\beta < \pi/2$, that is 
impossible.  This contradiction eliminates the case $\ell=0$, $k=1$.  

That leaves the case $\ell = 0$, $k=2$.
 Then we have from (\ref{eq:344})
\begin{eqnarray*}
\frac n  P + m \big(1- \frac 1 {P}\big)   &=& 4     \\
m-n &=& (m-4)P
\end{eqnarray*}
Since $m - n > 0$, we have $m \ge 5$.
Since $P \ge 3$ and $P = (m-n)/(m-4)$, we have 
\begin{eqnarray*}
3 &\le&  \frac {m-n} {m-4} \\
3m-12 &\le& m-n \\
m \le 6 - n/2
\end{eqnarray*}
Since $5 \le m$ we have only the possibilities $m=5$, $n= 0$, 1, or 2, and $m=6$, $n=0$.  
Since $P= (m-n)/(m-4)$,  if $m=6$ we have $P = 3$, and in that case $P+Q+R = 4$ as claimed.
So we may assume that $m=5$, in which case we have $P = 5$, $4$, or $3$, according as $n=0$, 1, or 2.
Then $\alpha = \pi/(2P)$ is $\pi/10$, $\pi/8$, or $\pi/6$.  In case $P=3$ we have $P+Q+R = 4$ as desired.
In case $P=4$ we have $\alpha = \pi/8$, and $R=1$, and this case was already treated at the beginning 
of the proof.  In case $m=P=5$, we have $\alpha = \pi/10$.  This case cannot occur, by Lemma \ref{lemma:pioverten}.
  
We have now proved that in every case in which there is vertex splitting, either 
 $P+Q+R=4$, or  $\alpha = \pi/6$ and $\beta = \pi/3$, or conclusion (i) or conclusion (ii)   of the theorem
holds.  Now assume $P+Q+R=4$.
Then exactly one vertex angle of $ABC$ is split, and that angle is split into only two pieces.
   The other two angles must therefore each
be equal to an angle of the tile $T$.  If two of them are equal to the same angle of $T$, then $ABC$ is 
isosceles.  We may assume $T$ is not equilateral, since in that case no vertex splitting could occur.
The two equal angles cannot be $\beta$ since $2\beta + \gamma  > \alpha + \beta + \gamma = \pi$.
Therefore the two equal angles are $\alpha$.  Since an $\alpha$ angle cannot split, the vertex angle
must split.  It could split into $\alpha + \beta$, or $\alpha + \gamma$, or $2\beta$, or $\beta + \gamma$. 
The latter is impossible, since then $2\alpha + \beta + \gamma = \pi$, which contradicts $\alpha + \beta + \gamma = \pi$
since $\alpha > 0$.  If the vertex angle splits into $\alpha + \beta$, then the angle sum of $ABC$ is 
$\pi = 3\alpha + \beta$;  subtracting $\alpha + \beta = \pi/2$ we have $\alpha = \pi/4 = \beta$, contradicting
the hypothesis that $T$ is not isosceles.  If the vertex angle splits into $\alpha + \gamma$, then then
angle sum of $ABC$ is $3\alpha + \gamma = \pi$.  Since $\gamma = \pi/2$, we have $\alpha = \pi/6$.  Hence
$\beta = \pi/3$ and $\alpha + \gamma = 2\pi/3 = 2\beta$.  Hence the vertex angle of $ABC$ is $2\beta$.
 If the vertex angle splits into $2\beta$, then conclusion (ii) of the 
theorem holds, i.e., the tile is similar to half of $ABC$. 

We have proved that if the two unsplit angles of $ABC$ are equal,  conclusion (ii) of the theorem holds.
Therefore we may assume that the two unsplit angles of $ABC$ are not equal; hence they are 
equal to two different angles of $T$.  Therefore triangle $ABC$ is similar to $T$, which is conclusion (i) 
of the theorem.   

We had proved that either $P+Q+R=4$, or $\alpha = \pi/6$, or one of the conclusions of the theorem holds. 
We have now disposed of the case $P+Q+R=4$.  
That leaves only the case when $\alpha = \pi/6$ to consider.  But in that case, 
by Theorem~\ref{theorem:isosceles}, conclusion (iii) of the theorem holds.  That completes the proof. 
\medskip

\section{Conclusions}
We have  classified the tilings of equilateral and isosceles triangles $ABC$,  and also 
 tilings in which the tile is a right triangle.  Our main results are summarized
in these theorems:

\begin{theorem}
If $ABC$ is a right triangle and $T$ is similar to 
$ABC$, then 
\smallskip

\noindent
(i) If \, $\tan \alpha $ is rational,  then  $ABC$ can be $N$-tiled by $T$
if and only if $N$ is a square, or $N = e^2 + f^2$ for integers $e$ and $f$ such that $\tan \alpha = e/f$. 
\smallskip

\noindent
(ii) if $\alpha = \pi/6$, then $ABC$ can be $N$-tiled by $T$ if $N$ is a square, or three times a square,
and for no other values of $N$. 

\smallskip
\noindent
(ii) If \,$\tan \alpha$ is not rational and $\alpha \neq \pi/6$,  then $ABC$ can be $N$-tiled by $T$ if and only if $N$ is a square.
\end{theorem}

\noindent{\em Proof}.  All the tilings mentioned in the theorem were exhibited and illustrated,
so it only remains to prove that there are no more.  
Suppose $ABC$ is a right triangle and $T$ is similar to $ABC$,  and $ABC$ is $N$-tiled by $T$.

{\em Ad} (i): Suppose $\tan \alpha$ is rational. 
Then $T$ is not a 30-60-90 triangle, since $\tan (\pi/6) = 1/\sqrt 3$ is not rational. 
Then by  Theorem~\ref{theorem:Similar2},  $N$ is either a square, or  a sum of squares $e^2 + f^2$ where $\tan \alpha = e/f$.

{\em Ad} (ii).  Suppose $\alpha = \pi/6$.  Then by Theorem~\ref{theorem:Similar2}, $N$ is either a square, or 
three times a square,
as claimed.

{\em Ad} (iii). Suppose $\tan \alpha$ is not rational and $\alpha \neq \pi/6$.  Then by Theorem~\ref{theorem:Similar2},
$N$ is a square.   That completes the proof. 

\begin{theorem} If $ABC$ is equilateral,  then $ABC$ can be $N$-tiled if $N = 3m^2$ and the tile is isoceles with 
base angles $\pi/6$,  or if $N = 6m^2$ or $2m^2$  and the tile is a right angle with $\alpha = \pi/6$,  or if $N$ is a square
and the tile is similar to $ABC$.   If $T$ is a right triangle with $\alpha \neq \pi/6$,  then there are
no $N$-tilings of $ABC$ by $T$ for any $N$.
\end{theorem}

\noindent{\em Remark}.  In \cite{beeson-nonexistence}, we prove that there are no tilings of an equilateral $ABC$ by 
{\em any} tile other than those mentioned here, except possibly by a non-isosceles tile with a $120^\circ$ angle.   The case of a tile with a $120^\circ$ angle is taken up in \cite{beeson120}. 
 But in {\em this}  paper, the only tilings of an equilateral 
$ABC$ that are ruled out are those by a right-angled tile.
\medskip

\noindent{\em Proof}.  The $3m^2$ tilings are formed by dividing $ABC$ into three tiles as in Fig.~\ref{figure:3-tilings},
that is, by three line segments from the center to the vertices, and then quadratically tiling each of those 
triangles.  The $2m^2$ tilings are formed by dividing $ABC$ in half along an altitude, and then 
quadratically tiling each half.   The $6m^2$ tilings are formed by dividing $ABC$ in half along an altitude, and then 
tiling each half into $3m^2$ tiles, as illustrated in Fig.~\ref{figure:27-tilings}.  If $N$ is a square, 
there is the quadratic tiling.  If the tile $T$ is a non-isosceles right triangle then, according to 
Theorem~\ref{theorem:gammaright}, $T$ is similar to half of $ABC$.  Then $T$ is a 30-60-90 triangle,
so $\alpha = \pi/6$.  By Theorem~\ref{theorem:isosceles}, the only possibile tilings by this tile 
are when $N = 2m^2$ or $N=6m^2$.  That completes the proof of the theorem.

\begin{theorem} If $ABC$ is isosceles but not equilateral, and $T$ is a right triangle or is 
similar to $ABC$,  then $ABC$ can be $N$-tiled  by tile $T$ if and 
only if one of the following conditions holds:
\smallskip

(i)  $N$ is a square and $T$ is similar to $ABC$, or 
\smallskip

(ii) $N$ is a square and $T$ is similar to half of $ABC$,  and $\alpha = \beta = \pi/4$, or
\smallskip

(iii) $N/2$ is a square, and $T$ is similar to half of $ABC$, or 
\smallskip

(iv) $N/6$ is a square,  $T$ is a 30-60-90 triangle, and the base angles of $ABC$
are $\pi/6$. 
\end{theorem}

\noindent{\em Proof}. 
Since $T$ is assumed to be a right triangle, if $T$ is not isosceles, then 
either $T$ is similar to $ABC$ or $T$ is similar to half of $ABC$, 
by Theorem~\ref{theorem:gammaright}. 
If $T$ is similar to half of $ABC$, 
then by Theorem~\ref{theorem:isosceles}, and the hypothesis that $ABC$ is not equilateral, one
of the listed alternatives holds.  If $T$ is a right isosceles triangle,  then $\alpha = \beta = \pi/4$,
so all the angles of $ABC$ are multiples of $\pi/4$.  Since $ABC$ is isosceles, then two of its angles 
are equal, so these two must be $\pi/4$, making $ABC$ and half of $ABC$ both similar to the tile $T$,
so conclusion (ii) holds.  (We can either tile $ABC$ quadratically, or tile each half of $ABC$ 
quadratically, using a smaller tile of the same shape.) 

Now suppose that $T$ is not a right triangle. 
If $T$ is similar to $ABC$
then by Lemma~\ref{lemma:penultimate}, either $N$ is a square, or 
$\gamma = \pi/2$.  Hence, under the assumption that $T$ is a not a right triangle, either 
(i) holds or  $T$ is not similar to $ABC$.
That completes the proof of the theorem.
 
\begin{theorem} If $ABC$ is similar to $T$, and  $ABC$ is not a right triangle, 
and not isosceles 
or equilateral,  then $ABC$ can be $N$-tiled by $T$ only when $N$ is a square.
\end{theorem}

\noindent{\em Proof}.  This follows immediately from Theorem~\ref{theorem:Similar2}.

Taken together these theorems provide a complete characterization of the triples $(ABC, N, T)$ such that 
$ABC$ can be $N$-tiled by $T$,  when $ABC$ is similar to $T$, or $T$ is a right triangle. 
The remaining possibilities for $ABC$ and $T$, except for the case when $T$ has a $120^\circ$ angle
and is not isosceles, 
are successfully treated in \cite{beeson-nonexistence} and \cite{beeson-triquadratics}, 
where it is shown that is exactly one more family of tilings not mentioned here.   The 
case when $T$ has a $120^\circ$ angle
and is not isosceles is taken up in \cite{beeson120}.

\bibliography{TriangleTiling}

\begin{thebibliography}{10}

\bibitem{beeson-nonexistence}
Michael Beeson.
\newblock Triangle tiling {II}: Some non-existence theorems.
 \begin{quotation}\noindent To appear, available on the author's website.
  \end{quotation}

\bibitem{beeson-triquadratics}
Michael Beeson.
\newblock Triangle tiling {III}: the triquadratic tilings.
 \begin{quotation}\noindent To appear, available on the author's website.
  \end{quotation}

\bibitem{beeson120}
Michael Beeson.
\newblock Triangle tiling {IV}: A non-isosceles tile with a 120 degree angle.
 \begin{quotation}\noindent To appear, available on the author's website.
  \end{quotation}

\bibitem{boltyanski2}
V.~G. Boltyanskii.
\newblock {\em Equivalent and Equidecomposable Figures}.
\newblock Heath, Boston, Massachusetts, 1963.
 \begin{quotation}\noindent Translated and adapted from the 1st Russian ed.
  (1956) by Alfred K. Henn and Charles E. Watts. \end{quotation}

\bibitem{boltyanski}
V.~G. Boltyanskii and I.~T. Gohberg.
\newblock {\em The decomposition of figures into smaller parts}.
\newblock University of Chicago Press, Chicago, Illinois, 1980.
 \begin{quotation}\noindent Translated and adapted from the Russian edition by
  Henry Christoffers and Thomas P. Branson. \end{quotation}

\bibitem{cohen}
Henri Cohen.
\newblock {\em Number Theory, Volume 1: Tools and Diophantine Equations}.
\newblock Springer, 2007.


\bibitem{goldberg}
M.~Goldberg and B.~M. Stewart.
\newblock A dissection problem for sets of polygons.
\newblock {\em American Mathematical Monthly}, 71:1077--1095, 1964.


\bibitem{golomb}
Solomon~W. Golomb.
\newblock Replicating figures in the plane.
\newblock {\em The Mathematical Gazette}, 48:403--412, 1964.


\bibitem{macmahon}
Percy~Alexander MacMahon.
\newblock {\em New Mathematical Pastimes}.
\newblock Cambridge University Press, Cambridge, England, 1921.


\bibitem{soifer}
Alexander Soifer.
\newblock {\em How Does One Cut a Triangle?}
\newblock Springer, 2009.


\end{thebibliography}
\bibliographystyle{plain-annote}

\end{document}